%% file: icml2024.tex
\documentclass[dvipsnames]{article}

\usepackage{hyperref}
\usepackage{url}
\usepackage{graphicx}
\usepackage{subfigure}
\usepackage{xcolor}
\usepackage{algorithm,algpseudocode}
\usepackage{mathtools}
\usepackage{soul} 
\usepackage{wrapfig, framed}
\usepackage[capitalize,noabbrev]{cleveref}

\usepackage[accepted]{icml2024}

\usepackage{amsmath}
\usepackage{amssymb}
\usepackage{mathtools}
\usepackage{amsthm}

\input{math_commands.tex}

\def\b{\vec{b}}
\def\u{\vec{u}}
\def\n{\hat{\vec{n}}}
\def\x{\vec{x}}
\def\y{\vec{y}}
\def\f{\vec{f}}
\def\d{\vec{d}}
\def\r{\vec{r}}
\def\fext{\f^\text{ext}}
\def\fluid{\Omega}
\def\solid{\Omega_s}
\def\air{\Omega_a}
\def\fullDomain{\mathcal{D}}
\def\neumannBoundary{\Gamma_n}
\def\dirichletBoundary{\Gamma_d}
\def\network{\bm{\mathcal{P}}^\textrm{net}}
\def\numOrtho{n_\textrm{ortho}}

\def\kernel{\bm{\mathcal{K}}}
\newcommand{\kernelComponent}[1]{\mathcal{K}_{#1}}

\def\image{\bm{\mathcal{I}}}
\newcommand{\imageComponent}[1]{\mathcal{I}_{#1}}

\def\levels{\mathcal{L}}
\def\level{\ell}

\def\weights{\bm{\mathcal{W}}}
\newcommand{\weight}[1]{\mathcal{W}_{#1}}
\def\biases{\bm{\mathcal{B}}}
\newcommand{\bias}[1]{\mathcal{B}_{#1}}

\def\A{A^{\image}}
\def\Areduced{{\tilde A}^{\image}}
\def\Adash{$A$\nobreakdash-}
\providecommand{\norm}[1]{\left\lVert#1\right\rVert}
\newcommand{\defeq}{\vcentcolon=}

\theoremstyle{plain}

\theoremstyle{definition}

\theoremstyle{remark}

\usepackage[textsize=tiny]{todonotes}

\icmltitlerunning{A Neural-Preconditioned Poisson Solver for Mixed Dirichlet and Neumann Boundary Conditions}

\begin{document}

\twocolumn[
\icmltitle{A Neural-Preconditioned Poisson Solver for \\
Mixed Dirichlet and Neumann Boundary Conditions}




\begin{icmlauthorlist}
\icmlauthor{Kai Weixian Lan}{ucd}
\icmlauthor{Elias Gueidon}{ucla}
\icmlauthor{Ayano Kaneda}{wsd}
\icmlauthor{Julian Panetta}{ucd}
\icmlauthor{Joseph Teran}{ucd}
\end{icmlauthorlist}

\icmlaffiliation{ucd}{University of California, Davis, USA}
\icmlaffiliation{ucla}{University of California, Los Angeles, USA}
\icmlaffiliation{wsd}{Waseda University, Tokyo, Japan}

\icmlcorrespondingauthor{Kai Weixian Lan}{kai.weixian.lan@gmail.com}

\icmlkeywords{Machine Learning, fluid simulation, numerical methods}

\vskip 0.3in
]



\printAffiliationsAndNotice{}  

\begin{abstract}
We introduce a neural-preconditioned iterative solver for Poisson equations with mixed boundary conditions.
Typical Poisson discretizations yield large, ill-conditioned linear systems.
Iterative solvers can be effective for these problems, but only when equipped with powerful preconditioners.
Unfortunately, effective preconditioners like multigrid \citep{brandt:1977:multi} require costly setup phases that must be re-executed every time domain shapes or boundary conditions change, forming a severe bottleneck for problems with evolving boundaries.
In contrast, we present a neural preconditioner trained to efficiently approximate the inverse of the discrete Laplacian in the presence of such changes.
Our approach generalizes to domain shapes, boundary conditions, and grid sizes outside the training set.
The key to our preconditioner's success is a novel, lightweight neural network architecture featuring spatially varying convolution kernels and supporting fast inference.
We demonstrate that our solver outperforms state-of-the-art methods like algebraic multigrid as well as recently proposed neural preconditioners on challenging test cases arising from incompressible fluid simulations.
\end{abstract}

\section{Introduction}

The Poisson equation is ubiquitous in scientific computing: it governs a wide array of physical phenomena, arises as a subproblem in many numerical algorithms, and serves as a model problem for the broader class of elliptic PDEs.
Solving its discretized form, a linear system of equations involving the discrete Laplacian matrix, is the bottleneck in many engineering and scientific applications.
These large, symmetric positive definite and sparse systems of equations are notoriously ill-conditioned.
Fast Fourier Transforms \citep{cooley:1965:algorithm} are optimal for these problems when discretized on trivial geometric domains, but they are not applicable for practical domain shapes.
Direct methods like Cholesky factorization \citep{golub2012matrix} resolve conditioning issues but suffer from loss of sparsity/fill-in and are prohibitively costly in practice when per-time-step refactoring is necessary (\emph{e.g.}, with changing domain shape or coefficients).
Iterative methods like preconditioned conjugate gradient (PCG) \citep{saad:2003:sparse} and multigrid \citep{brandt:1977:multi} can achieve good performance, but an optimal preconditioning strategy is not generally available.
Though multigrid preconditioners can guarantee modest iteration counts, computational overhead associated with solver creation and other per-iteration costs can dominate runtimes in practice.
This is especially true for problems posed on evolving domains, where multigrid hierarchies \emph{must be rebuilt at each time step}, and for nonlinear problems that require per-iteration rebuilds.
Unfortunately, there is no clear algorithmic solution.

Recently, machine learning techniques have shown promise for these problems, eliminating setup costs at runtime by training a general-purpose solver \emph{once} on a diverse set of systems offline.
\citet{tompson:2017:accelerating} showed that a network (FluidNet) can be used to generate an approximate inverse across domain shapes, albeit only with Neumann boundary conditions.
\citet{kanedo:2023:dcdm} developed the Deep Conjugate Direction Method (DCDM), which improves on FluidNet by applying a similar network structure as a \emph{preconditioner} for an orthogonalized gradient descent on the matrix norm of the error, enabling highly accurate solutions to be obtained.
While DCDM is similar to PCG, the nonlinearity of their approximate inverse required a generalization of the PCG method.
Also, their approach only supports Neumann pressure boundary conditions.
We build on the DCDM approach and generalize it to domains with mixed Dirichlet and Neumann boundary conditions.
Notably, these problems arise in simulating free-surface liquid flows.
DCDM fails on these cases, yet we show that a novel, lighter-weight network structure can be used effectively in its iterative formalism.
In contrast to DCDM, our approximate inverse is a \emph{linear} operator and can handle mixed boundary conditions over time-varying fluid domains.
Furthermore, we demonstrate that this structure drastically improves performance over DCDM.

We design our network architecture to represent the dense nature of the inverse of a discrete Laplacian matrix.
That is, the inverse matrix for a discrete Laplace operator has the property that local perturbations anywhere in the domain have non-negligible effects at all other domain points.
Our network structure uses a hierarchy of grid scales to improve the resolution of this non-local behavior over what is possible with DCDM's network architecture.
In effect, the process of transferring information across the hierarchy from fine grid to increasingly coarse grids and back again facilitates rapid propagation of information across the domain.
This structure is similar to multigrid but has some important differences.
We incorporate the effects of Dirichlet and Neumann conditions at irregular boundaries with a novel convolution design.
Specifically, we use stencils that learn spatially varying weights based on a voxel's proximity to the boundary and the boundary condition types encoded there.
We show that our network outperforms state-of-the-art preconditioning strategies, including DCDM, FluidNet, algebraic multigrid, and incomplete Cholesky,
performing our comparisons across a number of representative free-surface liquid and fluid flow problems.

In summary, our work makes the following contributions:
\begin{itemize}
    \item We introduce a novel light-weight architecture employing spatially varying convolutional kernels that is highly effective at approximating the inverse of structured-grid Laplacian matrices with arbitrary mixed Dirichlet and Neumann boundary conditions.
    \item We show that a simple loss function based on the residual of the linear system suffices for training in an unsupervised manner, producing a network that generalizes to systems not seen during training.
    \item We demonstrate through a comprehensive benchmark on challenging fluid-simulation test cases that, when paired with an appropriate iterative method, our neural-preconditioned solver dramatically outperforms state-of-the-art solvers like algebraic multigrid and incomplete Cholesky, as well as recent neural preconditioners like DCDM and FluidNet.
\end{itemize}
To promote reproducibility, we have released our full code and a link to our pretrained model at \url{https://github.com/kai-lan/MLPCG/tree/icml2024}.

\section{Related Work}

Many recent works leverage machine learning techniques to accelerate numerical linear algebra computations.

\citet{ackmann:2020:machine} use supervised learning to compute preconditioners from fully-connected feed-forward networks in semi-implicit time stepping for weather and climate models.
\citet{sappl:2019:deep} use convolutional neural networks (CNNs) \citep{lecunconv} to learn banded approximate inverses for discrete Poisson equations arising in incompressible flows discretized over voxelized spatial domains.
However, their loss function is the condition number of the preconditioned operator, which is prohibitively costly at high resolution.
\citet{ozbay:2021:poissoncnn} also use CNNs to approximate solutions to Poisson problems arising in incompressible flow discretized over voxelized domains, but they do not learn a preconditioner and their approach only supports two-dimensional square domains.
Our approach is most similar to those of \citet{tompson:2017:accelerating} and \citet{kanedo:2023:dcdm}, who also consider discrete Poisson equations over voxelized fluid domains, however our lighter-weight network outperforms them and generalizes to a wider class of boundary conditions.
\citet{li:2023:learningpcg} build on the approach of \citet{sappl:2019:deep}, but use a more practical loss function based on the supervised difference between the inverse of their preconditioner times a vector and its image under the matrix under consideration.
Their preconditioner is the product of easily invertible, sparse lower triangular matrices.
Notably, their approach works on discretizations over unstructured meshes.
\citet{gotz:2018:bjacobi} learn Block-Jacobi preconditioners using deep CNNs.
The choice of optimal blocking is unclear for unstructured discretizations, and they use machine learning techniques to improve upon the selection.
\citet{kopanivcakova:2024:deeponet} develop a general collection of neural DeepONet preconditioners.
The most effective of these is akin to a model reduction technique with a set of reduced basis functions taken from a trained DeepONet.
However, each application of the preconditioner itself requires solving a large (usually dense) linear system.\\
\\
Numerous recent works use neural network architectures that are analogous to a multigrid hierarchy.
\citet{he2019mgnet} analyzed the similarities between the structure of a convolutional network and that of the multigrid method, and proposed their novel multigrid-structured network MagNet.
The UNet \citep{ronneberger:2015:unet} and MSNet architectures \citep{mathieu:2016:deep} are similar to a multigrid V-cycle in terms of data flow, as noted by \citet{cheng:2021:using} and \citet{azulay:2023:multihelmoltz}.
\citet{cheng:2021:using} use the multi-scale network architecture MSNet to approximate the solution of Poisson equations arising in plasma flow problems.
However, they only consider flows over a square domain in 2D.
\citet{azulay:2023:multihelmoltz} note the similarity between the multi-scale UNet architecture and a multigrid V-cycle.
They use this structure to learn preconditioners for the solution of heterogeneous Helmholtz equations.
\citet{eliasof:2023:multigrid} also use a multigrid-like architecture for a general class of problems.
\citet{huang:2023:optimalmultigrid} use deep learning to generate multigrid smoothers at each grid resolution that effectively smooth high frequencies:
CNNs generate the smoothing stencils from matrix entries at each level in the multigrid hierarchy.
This is similar to our boundary-condition-dependent stencils, however we note that our network is lighter-weight and allowed to vary at a larger scale during learning.
Furthermore, optimal stencils are known for the problems considered in this work, and we provide evidence that our solver outperforms them.

\section{Motivation: Incompressible Fluid Simulation with Mixed B.C.s}
While our solver architecture can be applied to any Poisson equation
discretized on a structured grid, our original motivation was to accelerate
a popular method for incompressible inviscid fluid simulation based on the splitting scheme
introduced by \cite{chorin:1967:numerical}. The fluid's velocity
$\u(\x, t)$
is governed by the incompressible
Euler equations:
\begin{equation*}
\rho \left(
    \pder{\u}{t} + (\u \cdot \grad) \u
\right) + \grad p = \fext
\quad \text{s.t.} \quad
    \div \u = 0
\quad \text{in } \fluid,
\end{equation*}
where $\fluid$ is the domain occupied by fluid, pressure $p$ is the Lagrange multiplier for the incompressibility constraint $\div \u = 0$,
$\rho$ is the mass density, and $\fext$ accounts for external forces like gravity.
These equations are augmented with initial conditions
${\u(\x, 0) = \u^0(\x)}$ and ${\rho(\x, 0) = \rho^0}$ as well as the boundary conditions discussed in \pr{sec:bc}.
Incompressibility implies that 
mass density is conserved throughout the simulation ($\rho \equiv \rho^0$).

Chorin's scheme applies finite differences in time and splits the integration from time $t^n$ to ${t^{n + 1} = t^n + \Delta t}$ into
two steps. First, a provisional velocity field $\u^*$ is obtained by an
\emph{advection step} that neglects pressure and incompressibility:
\begin{equation}
    \label{eqn:advect}
    \frac{\u^* - \u^n}{\Delta t} + (\u^n \cdot \grad) \u^n = \frac{1}{\rho^0} \fext.
\end{equation}
Second, a \emph{projection step} obtains $\u^{n + 1}$ by eliminating divergence from $\u^*$:
\begin{align}
    - \div \frac{1}{\rho^0} \grad p^{n + 1} &= -\frac{1}{\Delta t}\div \u^*,
    \label{eqn:poisson}
    \\
    \frac{\u^{n + 1} - \u^*}{\Delta t} &= -\frac{1}{\rho^0} \grad p^{n + 1}.
    \label{eqn:project}
\end{align}
Equations~\ref{eqn:advect}-\ref{eqn:project} hold inside $\fluid$, and we
have deferred discussion of boundary conditions to \pr{sec:bc}.
The bottleneck of this full process is \pr{eqn:poisson}, which is a Poisson equation since $\rho^0$ is spatially constant.


\vfill\null

\subsection{Boundary Conditions}
\label{sec:bc}
\begin{wrapfigure}{r}{0.15\textwidth}
    \vspace{-3.05em}
\includegraphics[width=0.15\textwidth]{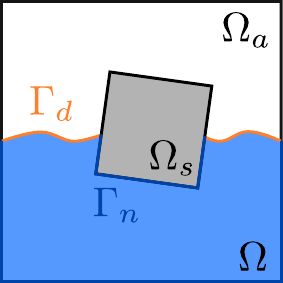}
    \vspace{-2.0em}
\end{wrapfigure}
Our primary contribution is our ability to handle both Neumann and
Dirichlet boundary conditions for the Poisson equation.
We assume the computational domain $\fullDomain$
is decomposed into $\fullDomain = \fluid \cup \air \cup \solid$,
as sketched in the inset,
where $\air$ denotes free space and $\solid$ the region filled with solid.
This decomposition induces a partition of the
fluid boundary
$\partial \Omega = \neumannBoundary \cup \dirichletBoundary$.
Boundary $\neumannBoundary$ contains the fluid-solid interface
and the intersection $\partial \Omega \cap \partial \fullDomain$
(\emph{i.e.}, the region outside $\fullDomain$ is treated as solid);
on it a free-slip boundary condition is imposed:
${\u(\x, t) \cdot \n(\x) = u^\neumannBoundary(\x, t)}$,
where $\n$ denotes the outward-pointing unit normal.
This condition on $\u$ translates via \pr{eqn:project} into a Neumann condition on \pr{eqn:poisson}:
\begin{equation}
    \label{eqn:bc_neumann}
    \n \cdot \grad p^{n + 1} = \frac{\rho_0}{\Delta t}(\n \cdot \u^* - u^\neumannBoundary) \quad \text{on } \neumannBoundary.
\end{equation}
Free-surface boundary $\dirichletBoundary$ represents the interface between the fluid
and free space. Ambient pressure $p_a$ then imposes on \pr{eqn:poisson} a Dirichlet condition
$p^{n + 1} = p_a \text{ on } \dirichletBoundary$.
In our examples, we set $p_a = 0$.

The Dirichlet conditions turn out to make solving \pr{eqn:poisson}
fundamentally more difficult: while the DCDM paper \cite{kanedo:2023:dcdm}
discovered that a preconditioner blind to the domain geometry and trained solely
on an empty box is highly effective for simulations featuring pure Neumann
conditions, the same is not true for Dirichlet (see \pr{fig:res_128}).

\subsection{Spatial Discretization}
We discretize the full domain $\fullDomain$ using a regular marker-and-cell (MAC) staggered grid with $n_c$ cubic elements \cite{harlow:1964:pic}.
The disjoint subdomains $\fluid$, $\air$, and $\solid$ are each represented by a per-cell rasterized indicator field;
these are collected into a 3-channel image, stored as a tensor $\image$. In the case of a 2D square with $n_c = N^2$, this tensor is of shape $(3, N, N)$,
and summing along the first index yields a single-channel image filled with ones.

Velocities and forces are represented at the \emph{corners} of this grid, and for smoke simulations the
advection step \pr{eqn:advect} is implemented using an explicit semi-Lagrangian
method \citep{stam:1999:stable, robert:1981:stable}.
For free-surface simulations, advection is performed by interpolating fluid velocities from the grid
onto particles responsible for tracking the fluid state, advecting those particles,
and then transferring their velocities back to the grid. We use a PIC/FLIP blend transfer scheme with a 0.99 ratio \citep{zhu:2005:sand-fluid}.

Pressure values are stored at element \emph{centers}, and the Laplace operator
of \pr{eqn:poisson} is discretized into a sparse symmetric matrix ${\A \in \R^{n_c \times n_c}}$ using the standard second-order accurate finite difference stencil
(with 5 points in 2D and 7 in 3D) but with modifications to account for Dirichlet and Neumann boundary conditions:
stencil points falling outside $\fluid$ are dropped,
\begin{wrapfigure}{r}{0.30\textwidth}
    \includegraphics[width=0.3\textwidth]{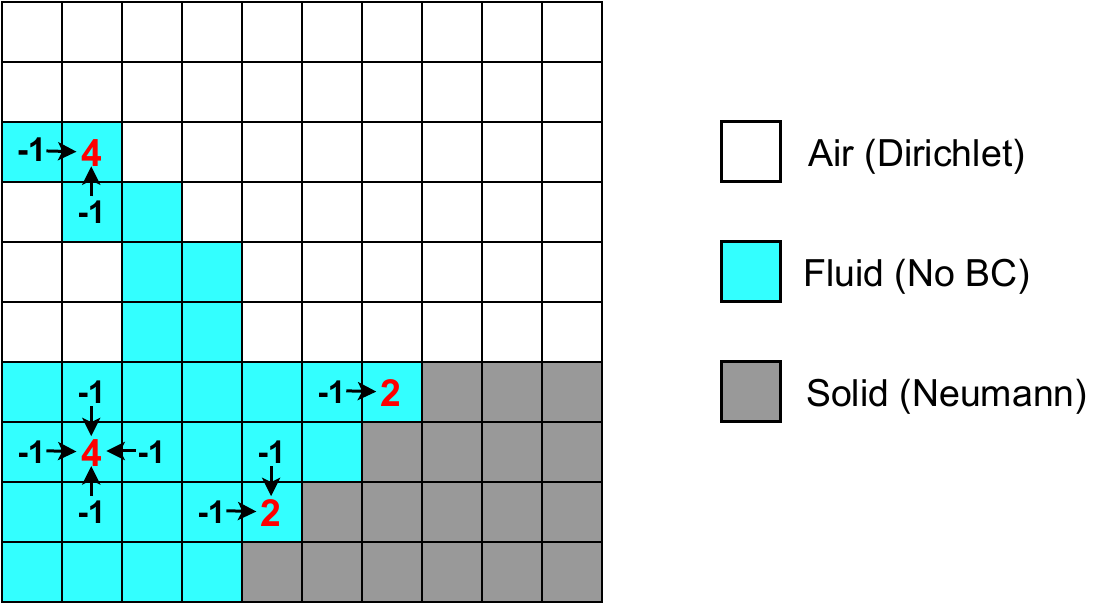}
\end{wrapfigure}
and the central value (\emph{i.e.}, the diagonal matrix entry)
is determined as the number of neighboring cells belonging to either $\fluid$ or $\air$.
Examples of these stencils are visualized in 2D in the inset.
Rows and columns corresponding to cells outside $\fluid$ are left empty, meaning $\A$
typically has a high-dimensional nullspace.
These empty rows and columns are removed before solving, obtaining a smaller positive definite matrix
${\Areduced \in \R^{n_f \times n_f}}$, where $n_f$ is the number of fluid cells.

The right-hand side of \pr{eqn:poisson} is discretized using the standard MAC
divergence finite difference stencil into a vector $\b \in \R^{n_c}$, which also
receives contributions from the Neumann boundary. Entries of this vector
corresponding to cells outside $\fluid$ are removed to form right-hand side vector $\tilde \b \in \R^{n_f}$
of the reduced linear system representing the discrete Poisson equation:
\begin{equation}
\label{eqn:linear_sys}
\Areduced \tilde \x = \tilde \b,
\end{equation}
where $\tilde \x \in \R^{n_f}$ collects the fluid cells' unknown pressure values (a discretization of $p^{n + 1}$).

The constantly changing domains and boundary conditions of a typical fluid simulation mean traditional
preconditioners for \pr{eqn:linear_sys} like multigrid or incomplete Cholesky, as well as direct sparse Cholesky factorizations,
must be \emph{rebuilt at every time step}. This prevents their high fixed costs from being amortized across frames
and means they struggle to outperform a highly tuned GPU implementation of unpreconditioned CG.
This motivates our neural-preconditioned solver which, after training, instantly adapts to arbitrary subdomain shapes
encoded in $\image$.

\section{Neural-preconditioned Steepest Descent with Orthogonalization}

Our neural-preconditioned solver combines a carefully chosen
iterative method (\pr{sec:algorithm}) with a preconditioner based on a novel neural network architecture (\pr{sec:architecture})
inspired by multigrid.


\subsection{Algorithm}
\label{sec:algorithm}

For symmetric positive definite matrices $A$ (like the discrete Laplacian $\Areduced$ from \pr{eqn:linear_sys}), the preconditioned
conjugate gradient (PCG) algorithm \citep{shewchuk:1994:cg} is by far the most
efficient iterative method for solving linear systems $A \x = \b$ when an effective preconditioner is available.
Unfortunately, its convergence rate is known to degrade when the preconditioner
itself fails to be symmetric, as is the case for our neural preconditioner (\pr{sec:symmetry}).
\citet{bouwmeester:2015:nonsym} have shown that good convergence can be recovered for nonsymmetric
multigrid preconditioners
using the ``flexible PCG'' variant at the expense of an additional dot product.
However, this variant turns out to perform suboptimally with our neural preconditioner, as shown
in \pr{apx:pcgs}. Instead, we adopt the preconditioned steepest descent with orthogonalization (PSDO)
method proposed in \cite{kanedo:2023:dcdm}, shown in \pr{alg:npsdo}, which performs well even
for their \emph{nonlinear} preconditioner.

\begin{algorithm}[h!]
    \caption{Neural-preconditioned Steepest Descent with $A$-Orthogonalization (NPSDO).}
    \label{alg:npsdo}
    \begin{algorithmic}[0]
       \State Given linear system ($A$, $\mathbf{b}$), image $\image$, and trained network $\network$
       \State $\r_0 = \b - A \x_0$
       \State $k = 0$
       \While{$\norm{\r_k} \ge \epsilon$}
       \State $k  =  k + 1$
       \State $\d_k = \network\left(\image,\frac{\r_{k-1}}{\norm{\r_{k-1}}}\right)$
       \For{$k - \numOrtho\leq i <k$}
       \State $\mathbf{d}_k = \mathbf{d}_k - \frac{\d_k^\top A \d_i}{\d_i^\top A \d_i}\mathbf{d}_i$
       \EndFor
       \State $\alpha_{k} = \frac{\r_{k-1}^\top \d_k}{\d_{k}^\top A \d_{k}} $
       \State $\x_{k} = \x_{k-1}+\alpha_k \d_{k}$
       \State $\r_k = \b - A \x_k$
       \EndWhile
       \State \Return $\x_k$
    \end{algorithmic}
 \end{algorithm}

The PSDO algorithm can be understood as a modification of standard CG that
replaces the residual with the preconditioned residual
as the starting point for generating search directions
and, consequently, cannot enjoy many of the simplifications baked into the traditional algorithm.
Most seriously, {\Adash}orthogonalizing against only the previous search direction no longer suffices
to achieve {\Adash}orthogonality to all past steps. Therefore, iteration $k$ of PSDO obtains its step direction $\d_k$ by explicitly {\Adash}orthogonalizing the preconditioned residual against the last $\numOrtho$
directions (where $\numOrtho$ is a tunable parameter) before determining step length $\alpha_k$ with an exact line search.
PSDO reduces to standard preconditioned steepest descent (PSD) when $\numOrtho = 0$, and it is
mathematically equivalent to unpreconditioned CG when
$\numOrtho \ge 1$ and the identity operator is used as the preconditioner.
In the case of a symmetric preconditioner $P = L L^\top$,
PSDO differs from PCG by taking steps that are {\Adash}orthogonal rather than $L A L^{\top}$\nobreakdash-orthogonal.
When combined with our neural preconditioner, we call this algorithm NPSDO, presented formally in \pr{alg:npsdo}.
We empirically determined $\numOrtho = 2$ to perform well (see \pr{apx:ortho}) and use this value in all reported experiments.

\subsection{Neural Preconditioner}
\label{sec:preconditioner}
The ideal preconditioner for all iterative methods described in
\pr{sec:algorithm} is the exact inverse $A^{-1}$; with it,
each method would converge to the exact solution in a single step.
Of course, the motivation for using an iterative solver is that inverting
or factorizing $A$ is too costly (\pr{fig:cholesky_comparison}),
so instead we must seek an inexpensive approximation of $A^{-1}$.
Examples are incomplete Cholesky, which does its best to factorize $A$
with a limited computational budget, and multigrid, which applies
iterations of a multigrid solver.

Our method approximates the map $\r \mapsto A^{-1} \r$
by our neural network $\network(\image, \r)$.
Departing from recent works like that of \citet{kanedo:2023:dcdm}, we use a novel architecture
that both substantially boosts performance on pure-Neumann problems
and generalizes to the broader class of Poisson
equations with mixed boundary conditions by considering geometric information from $\image$.
The network performs well on 2D or 3D Poisson equations of varying sizes,
but to simplify the exposition, our figures and notation describe the method on small square grids of size $N \times N$.

We note that \pr{alg:npsdo} runs on linear
system $\Areduced \tilde \x = \tilde \b$, featuring vectors of smaller size $n_f$, but the network always operates on
input vectors of full size $n_c$, reshaped into $(N, N)$ tensors. Therefore, to evaluate $\tilde \d = \network(\image, \tilde \r)$,
$\tilde \r$ is first padded by inserting zeros into locations corresponding to cells in $\air$ and $\solid$,
and then those locations of the output are removed to obtain $\tilde \d \in \R^{n_f}$.

\subsubsection{Architecture}

\label{sec:architecture}
\begin{figure}
    \subfigure[Full network architecture sketched for $\levels = 3$ levels.]{\label{fig:architecture}\includegraphics[width=\columnwidth]{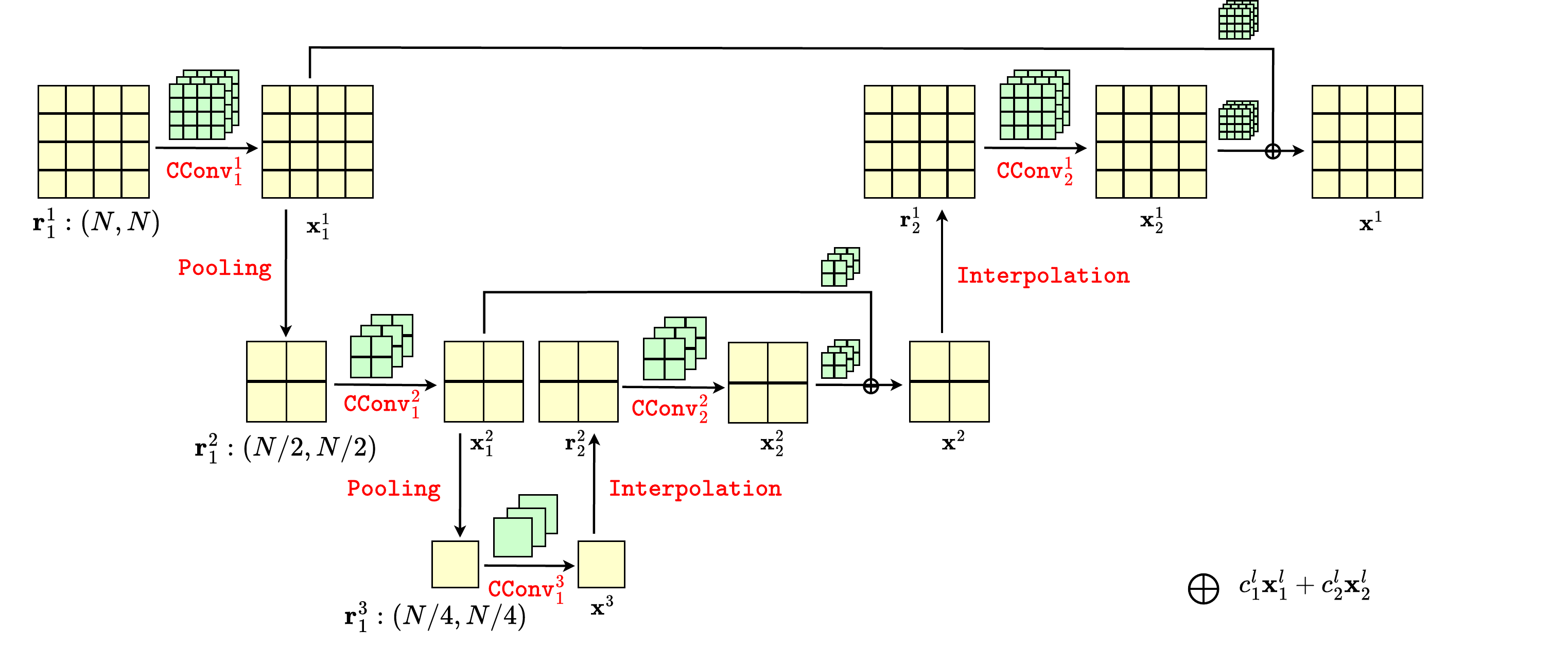}}
    \subfigure[Custom convolution block \\\texttt{CConv$^{\level}_i$} with $i=1,2$.]{\label{fig:network_conv}\includegraphics[width=0.64\columnwidth]{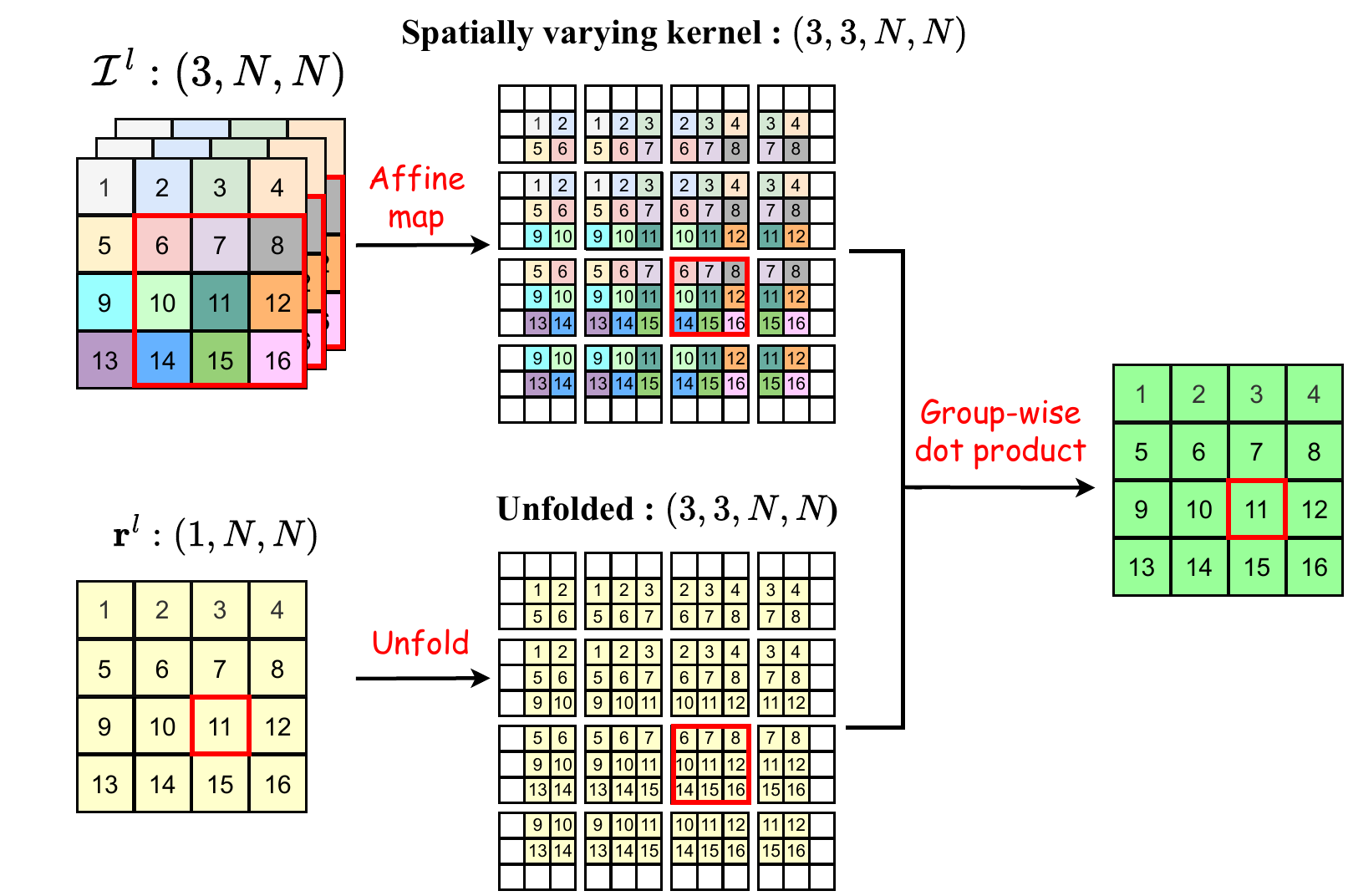}}
    \subfigure[Affine block \texttt{Aff$^{\level}_i$} with $i=1,2$.]{\label{fig:network_linear}\includegraphics[width=0.3\columnwidth]{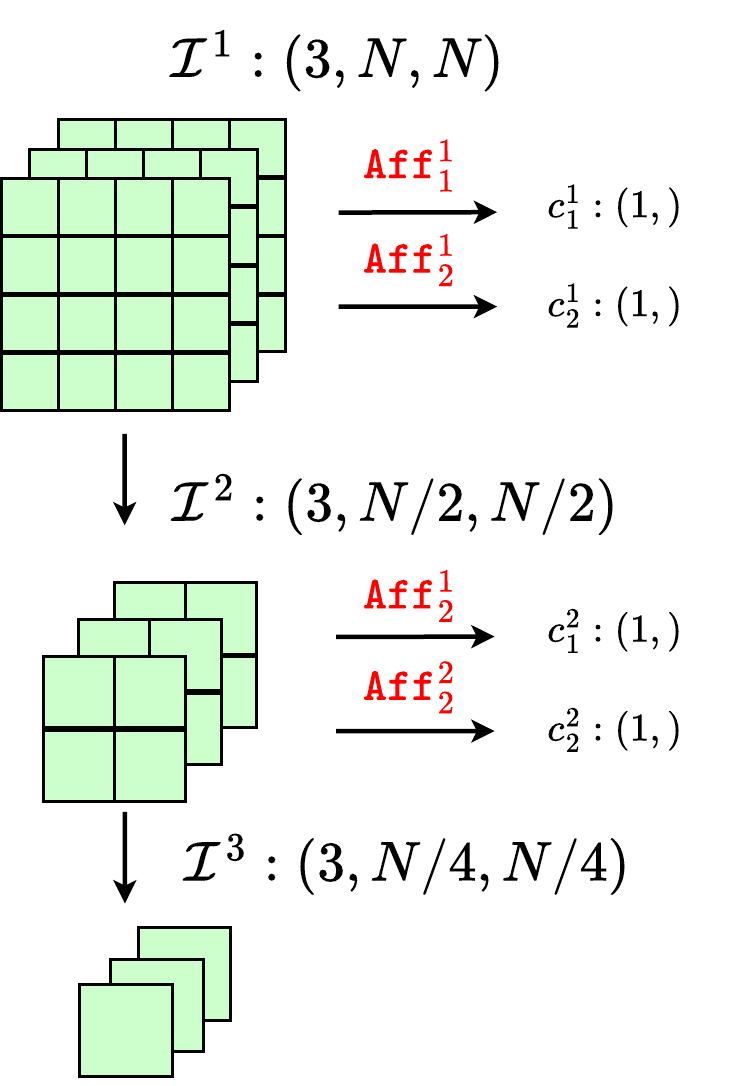}}
    \caption{Sketches of our network architecture.}
\end{figure}

Our multi-resolution network architecture (\pr{fig:architecture}) is inspired by UNet and geometric multigrid,
aiming to propagate information across the computational grid faster
than the one-cell-per-iteration of unpreconditioned CG.
The architecture is defined recursively, consisting of levels ${1 \le \level \le \levels}$ each executing the operations detailed in \pr{alg:network}.
A given level $\level$ operates on an input image $\image^{\level}$ of shape $(3, N^\level, N^\level)$ and an input
vector $\r^{\level}_1$. It performs a special image-dependent
convolution operation defined by our custom convolution block \texttt{CConv} and then downsamples the resulting
vector $\x^{\level}_1$, as well as $\image^{\level}$, to the next-coarser level
$\level + 1$ using average pooling with kernel size $2 \times 2$. 
The output of the level $\level + 1$ subnetwork is then upsampled with bilinear interpolation, 
run through another convolution stage, and finally linearly combined with
$\x^{\level}_1$ to obtain the output. At the finest level, $\image^{1} = \image$
and $\r^{1}_1 = \r$, while at the coarsest level only a single convolution
operation is performed.

\begin{algorithm}[h!]
    \caption{Operations executed by an $\levels$-level network in pseudocode form.}
    \label{alg:network}
    \begin{algorithmic}[0]
       \State Given input vector $\r$ and image $\image$
        \State $\r^1_1 = \r$
        \State $\image^1 = \image$
        \For {$\level = 1, \cdots, \levels - 1$}
            \State $\image^{\level+1} = \texttt{Pooling}(\image^\level)$

            \State $\x^\level_1 = \texttt{CConv}^\level_1 (\r^\level_1,\image^{\level})$
            \State $\r^{\level+1}_1 = \texttt{Pooling}(\x^\level_1)$
        \EndFor
        \State $\x^{\levels} = \texttt{CConv}^{\levels}_1 (\r^{\levels}_1, \image^{\levels})$
        \For {$\level = \levels-1, \ldots, 1$}
            \State $\r^\level_2 = \texttt{Interpolation} (\x^{\level+1})$
            \State $\x^{\level}_2 = \texttt{CConv}^{\level}_2(\r^\level_2, \image^\level)$
            \State $c^\level_1 = \texttt{Aff}^{\level}_1(\image^\level)$
            \State $c^\level_2 = \texttt{Aff}^{\level}_2(\image^\level)$
            \State $\x^{\level} = c^{\level}_1 \x^\level_1 + c^\level_2 \x^\level_2$
        \EndFor
       \State \Return $\x^1$
    \end{algorithmic}
 \end{algorithm}

One crucial difference between our network and existing neural solvers like
FluidNet \citep{tompson:2017:accelerating} is how geometric information from
$\image$ is incorporated. Past architectures treat this geometric data
on the same footing as input tensor $\r$, \emph{e.g.} feeding
both into standard multi-channel convolution blocks. However, we note that $\image$
determines the entries of $\A$, and so for the convolutions
to act analogously to the smoothing operations of multigrid,
this data should inform the weights of convolutions applied to $\r$.
This observation motivates our use of custom convolutional blocks $\texttt{CConv}$ whose \emph{spatially
varying kernels} depend on local information from $\image$. We note the
close connection between these varying kernels and attention \citep{bahdanau2014neural,vaswani2017attention}.

Each custom convolutional block (\pr{fig:network_conv}) at level $\level$ learns an affine map from a $3\times 3$ sliding window in $\image^{\level}$ to a $3\times 3$ kernel $\kernel^{i, j}$.
This affine map is parametrized by a weights tensor $\weights$ of shape $(3^2, 3, 3, 3)$ and a bias vector $\biases \in \R^{3^2}$.
Each entry of the block's output is computed as:
\begin{gather*}
    \left[\texttt{CConv}^\level(\r, \image^\level)\right]_{i, j} = \sum_{a,b=-1}^1 \kernelComponent{a, b}^{i, j} r_{i+a,j+b},
  \\
\kernelComponent{a, b}^{i, j} \defeq \sum_{c=1}^3\sum_{l,m=-1}^{1} \weight{3a+b,c,l,m}\imageComponent{c,i+l,j+m}^{\level} + \bias{3a + b}.
\end{gather*}
Out-of-bounds accesses in these formulas are avoided by padding $\image^{\level}$ with solid pixels
(\emph{i.e.}, the values assigned to cells in $\solid$)
and $\x$ with zeros.

In multigrid, solutions obtained on coarser grids
are \emph{corrections} that are added to the finer grids' solutions;
likewise, our network mixes in upsampled data from the lower level
using the linear combination $c^{\level}_1 \x^\level_1 + c^\level_2 \x^\level_2$.
The coefficients in this combination are defined by (i) convolving $\image^{\level}$ with a
(spatially constant) kernel $\overline \weights$ of shape $(3, 3, 3)$;
(ii) averaging to produce a scalar; and (iii) adding a scalar bias $\overline{\mathcal{B}}$.
For efficiency, these evaluation steps are fused into a custom affine block
(\pr{fig:network_linear}) that implements the formula:
   \begin{align*}
       \texttt{Aff}^{\level}(\image) = \overline{\mathcal{B}} + \frac{1}{3^2 n_c} \sum_{i,j=1}^{N^\level}\sum_{c = 1}^3\sum_{l,m=-1}^1 {\overline{\mathcal{W}}}_{c, l, m} \imageComponent{c, i + l, j + m}^{\level}.
   \end{align*}

Our custom network architecture has numerous advantages.
Its output is a linear function of the input vector (unlike the
nonlinear map learned by \cite{kanedo:2023:dcdm}), making it easier to interpret
as a preconditioner.
The architecture is also very lightweight: a model with $\levels=4$ coarsening levels has only $\sim\!25$k parameters.
Its simplicity accelerates network evaluations at solve time,
critical to make NPSDO competitive with the state-of-the-art solvers used in practice. We note that our solver is fully matrix free, with $\network$ relying only
on the image $\image$ of the simulation scene to infer information about $\A$.
Furthermore, since all network operations are formulated
in terms of local windows into $\image$ and $\r$, it can train and
run on \emph{problems of any size divisible by} $2^{\levels-1}$.

The 3D version of our architecture is a straightforward extension of the 2D
formulas above, simply using larger tensors with additional indices to account
for the extra dimension and extending the sums to run over these
indices.

\subsubsection{Symmetry and Positive Definiteness}
\label{sec:symmetry}
Our network enforces neither symmetry nor positive definiteness,
properties that would be needed to guarantee convergence of the standard PCG algorithm
(motiving our use of NPSDO).
However, we note that the additional operations in NPSDO vs PCG necessitated by
asymmetry (\emph{e.g.}, the extra $A$-orthogonalizations) do not add significant overhead
compared to the cost of network evaluation, which accounts for approximately 80\%
of the solver time. Furthermore, the strong, albeit suboptimal, performance of
the plain PCG algorithm (\pr{apx:pcgs}) suggests that the
deviations of our preconditioner from symmetry and positive definiteness are not
severe; this is confirmed by experiments reported in \pr{apx:symmetry_analysis}
that measure the magnitude of symmetry violation of our trained operator on
randomly generated vectors.

Nevertheless, we experimented with enforcing positive definiteness
by construction in two different ways: (i) concatenating our network with
a transposed copy of itself; and (ii) constraining ``post-smoothing'' convolution
blocks to be transposes of their corresponding ``pre-smoothing'' block and modifying the
shortcut connection.
Both architecture variants led to a significant degradation in preconditioner
quality when trained using the same methodology as our proposed network
(\pr{apx:symmetry_analysis}).

   \subsubsection{Training}
   We train our network $\network$ to approximate ${\A} \big \backslash \b$ when presented with image $\image$ and input vector $\b$.
   We calculate the loss for an example $(\image, \A, \b)$ from our training dataset as the residual norm:
   \begin{align*}
       Loss = \norm{\b - \A \network(\image, \b)}_2.
   \end{align*}
   We found the more involved loss function used in \cite{kanedo:2023:dcdm} not to benefit our network.

   Our training data set consists of 107 matrices collected from 11 different simulation scenes, some of domain shape (128, 128, 128) and others (256, 128, 128).
   For each matrix, we generate 800 right-hand side vectors using a similar approach to \cite{kanedo:2023:dcdm} but with far fewer Rayleigh-Ritz vectors.
   We first compute 1600 Ritz vectors using Lanczos iterations \citep{lanczos:1950:iteration} and then generate from them 800 random linear combinations.
   These linear combinations are finally normalized and added to the training set.
   To accelerate data generation, we create the right-hand sides for different matrices in parallel;
   it takes between 0.5 and 3 hours to generate the data for each scene.
   As Ritz vector calculation is expensive, we experimented with other approaches,
   like picking random vectors or constructing analytical eigenmodes for the Laplacian on $\fullDomain$ and masking out entries outside $\fluid$. Unfortunately these cheaper generation techniques led to degraded performance.

   In each epoch of training, we loop over all matrices in our dataset in shuffled order. For each matrix,
   we process all of its 800 right-hand sides in batches of 128, repeating five times.
   The full training process takes around 5 days on an AMD EPYC 9554P 64-Core Processor with an NVIDIA RTX 6000 GPU.
   We utilize the transfer learning technique \citep{transfer_learning}, training first a 5-level network and using those weights to initialize a 6-level network,
   which is subsequently fine-tuned and used for all experiments.

\begin{figure}

    \centering
    \includegraphics[width=0.3\columnwidth]{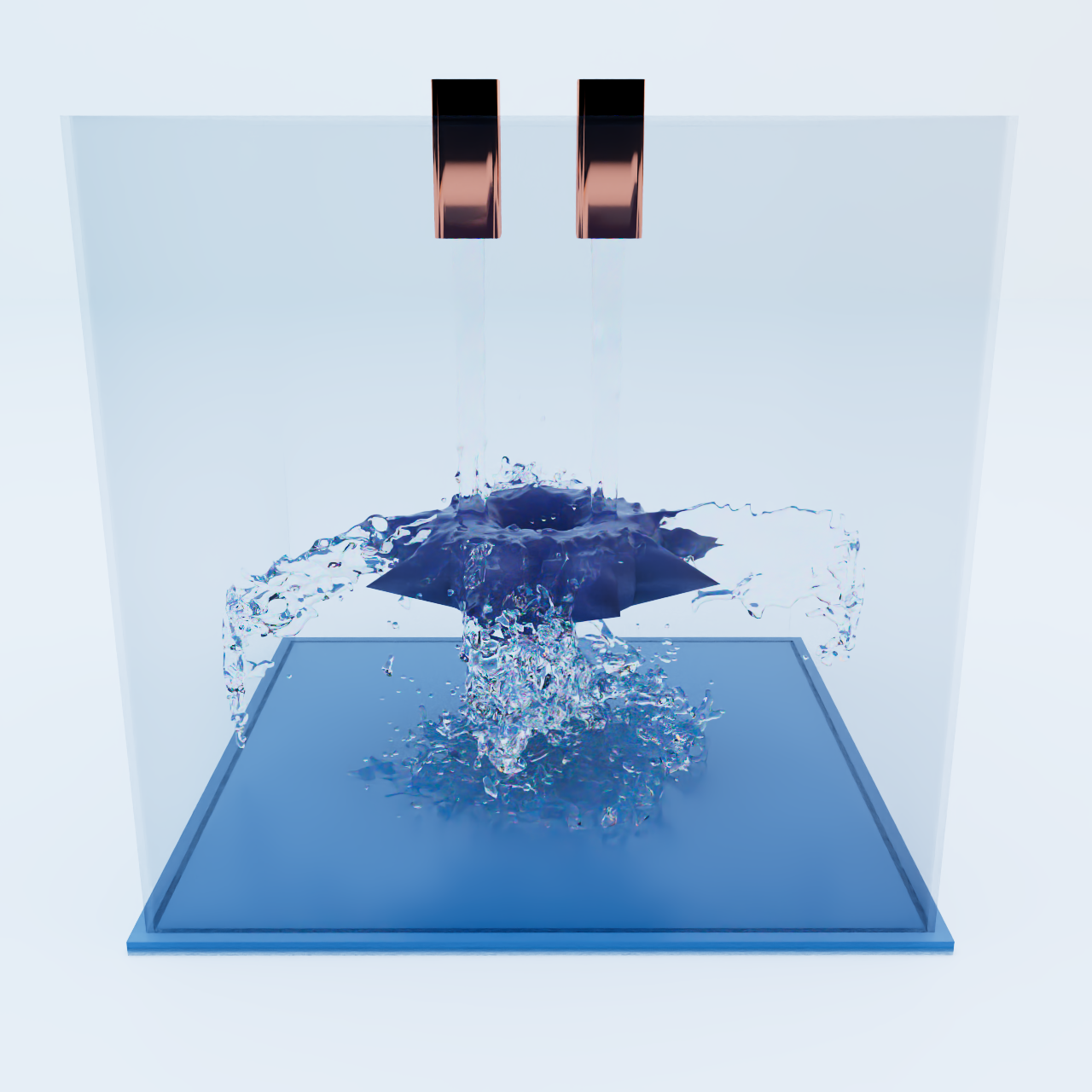}%
    \includegraphics[width=0.3\columnwidth]{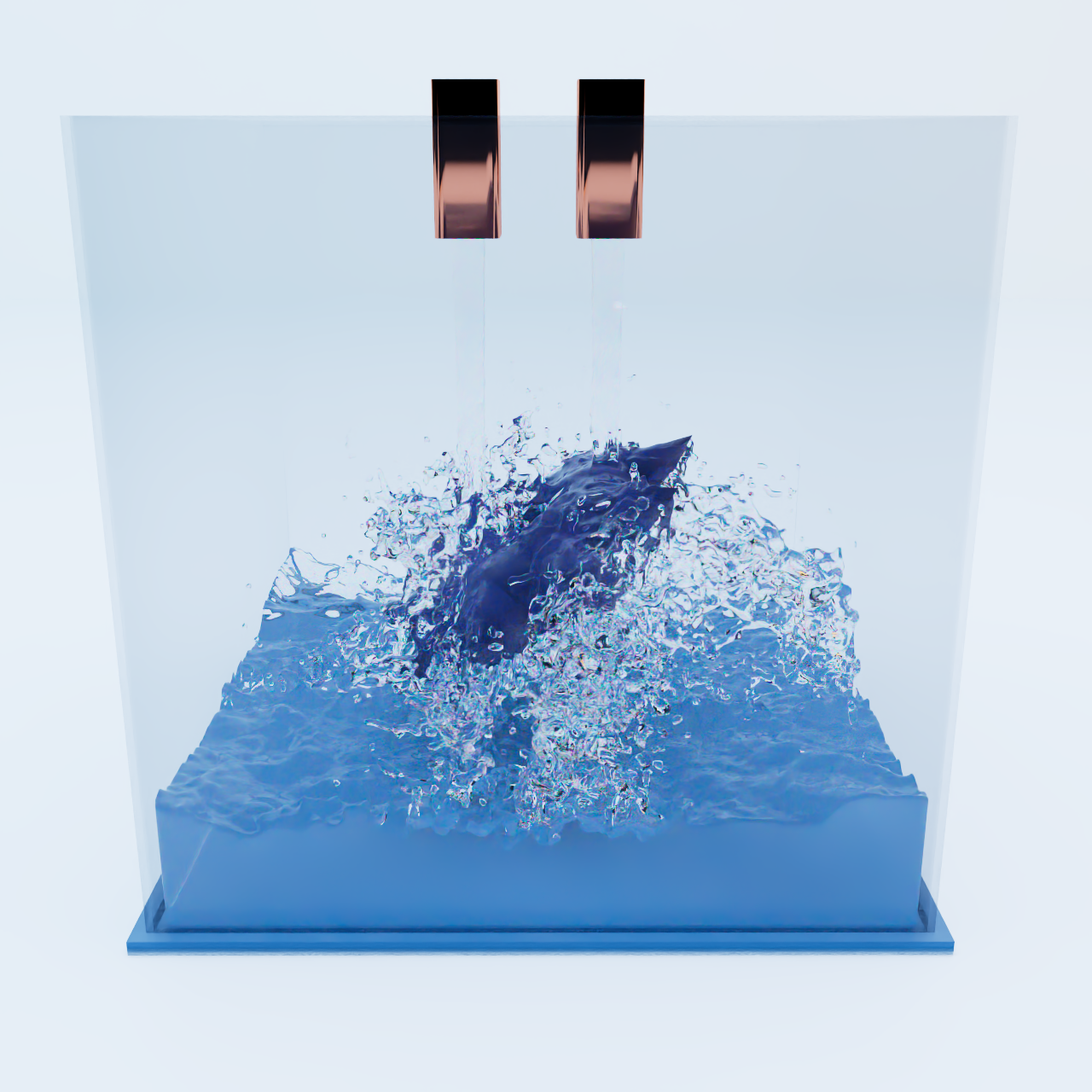}%
    \includegraphics[width=0.3\columnwidth]{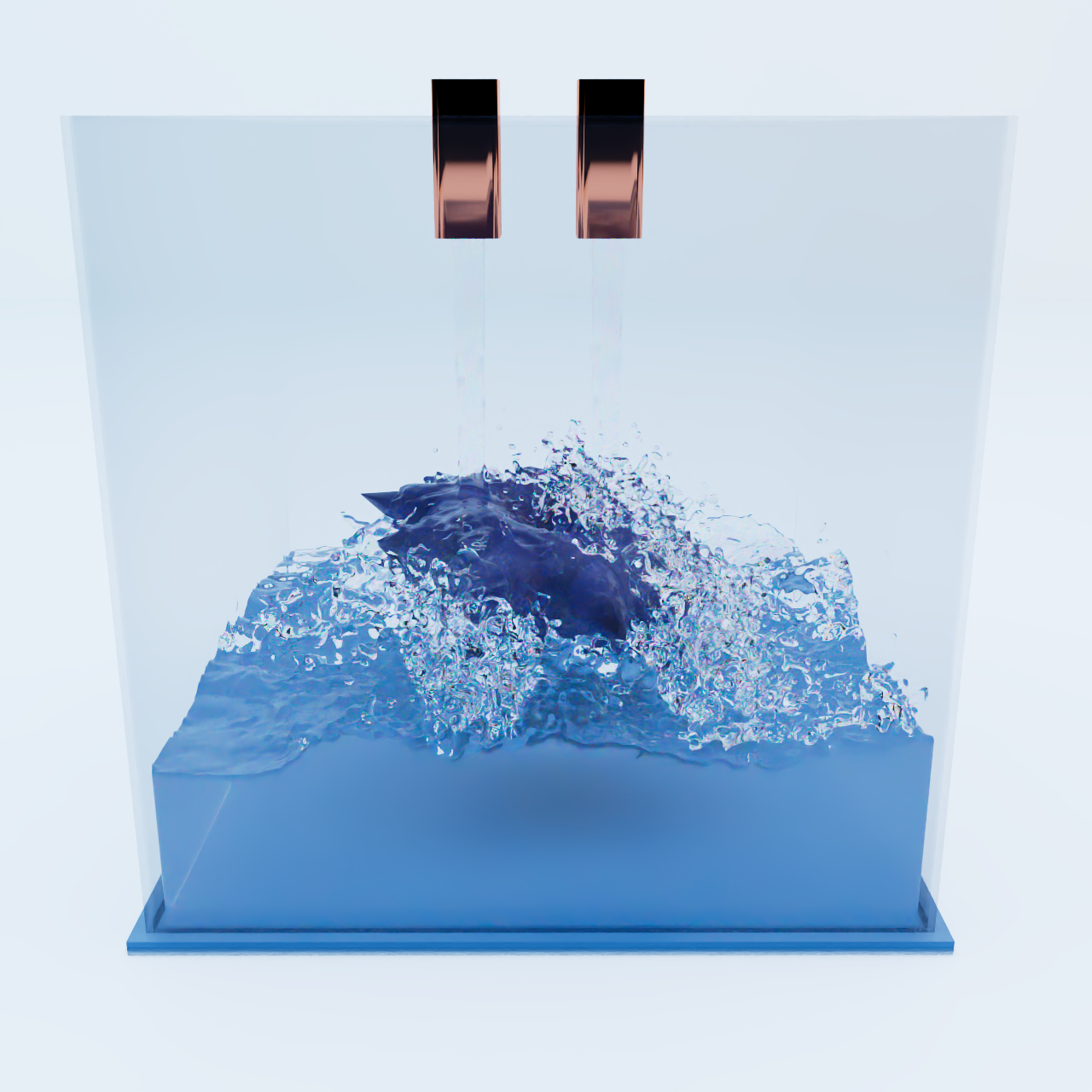} \\
    \includegraphics[width=0.3\columnwidth]{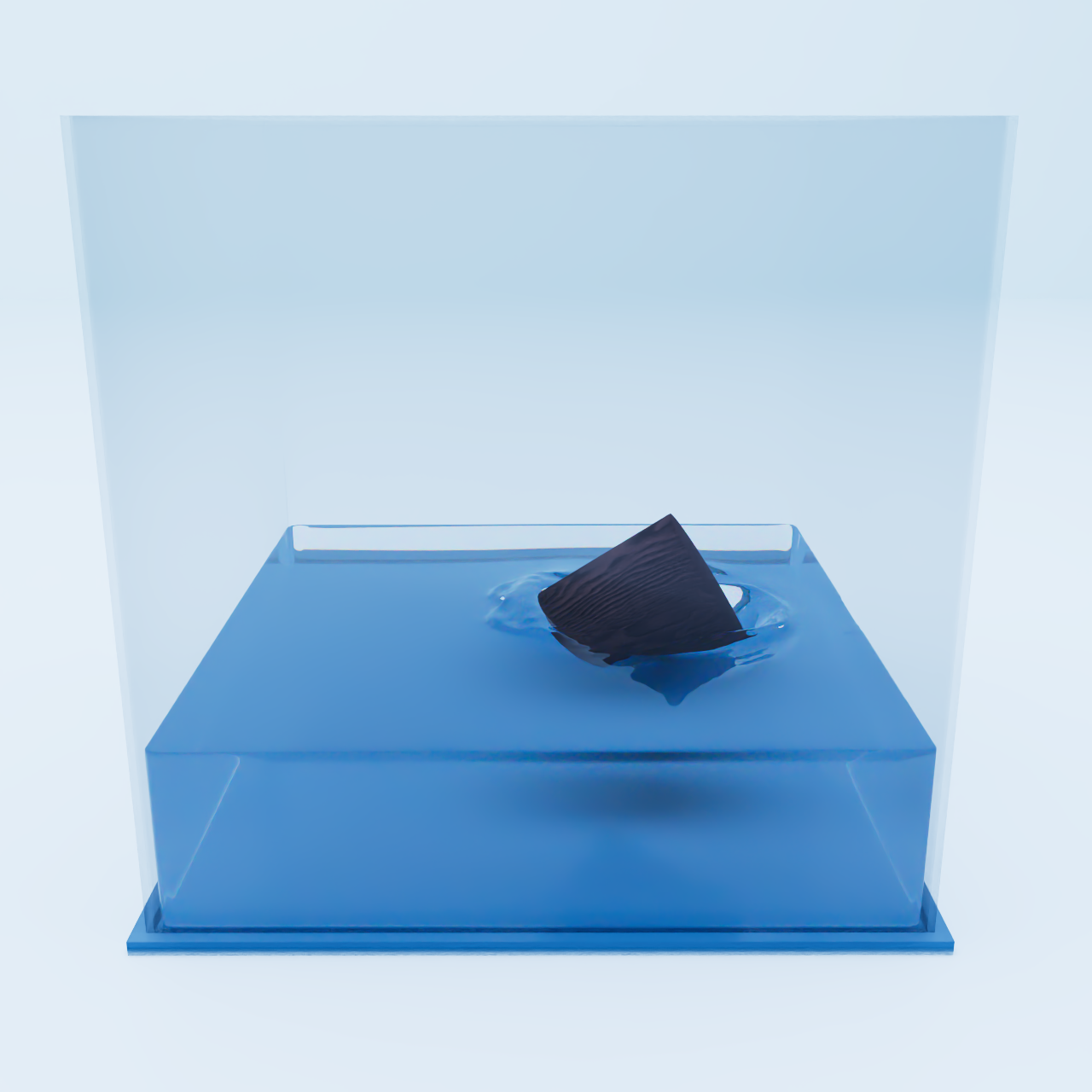}%
    \includegraphics[width=0.3\columnwidth]{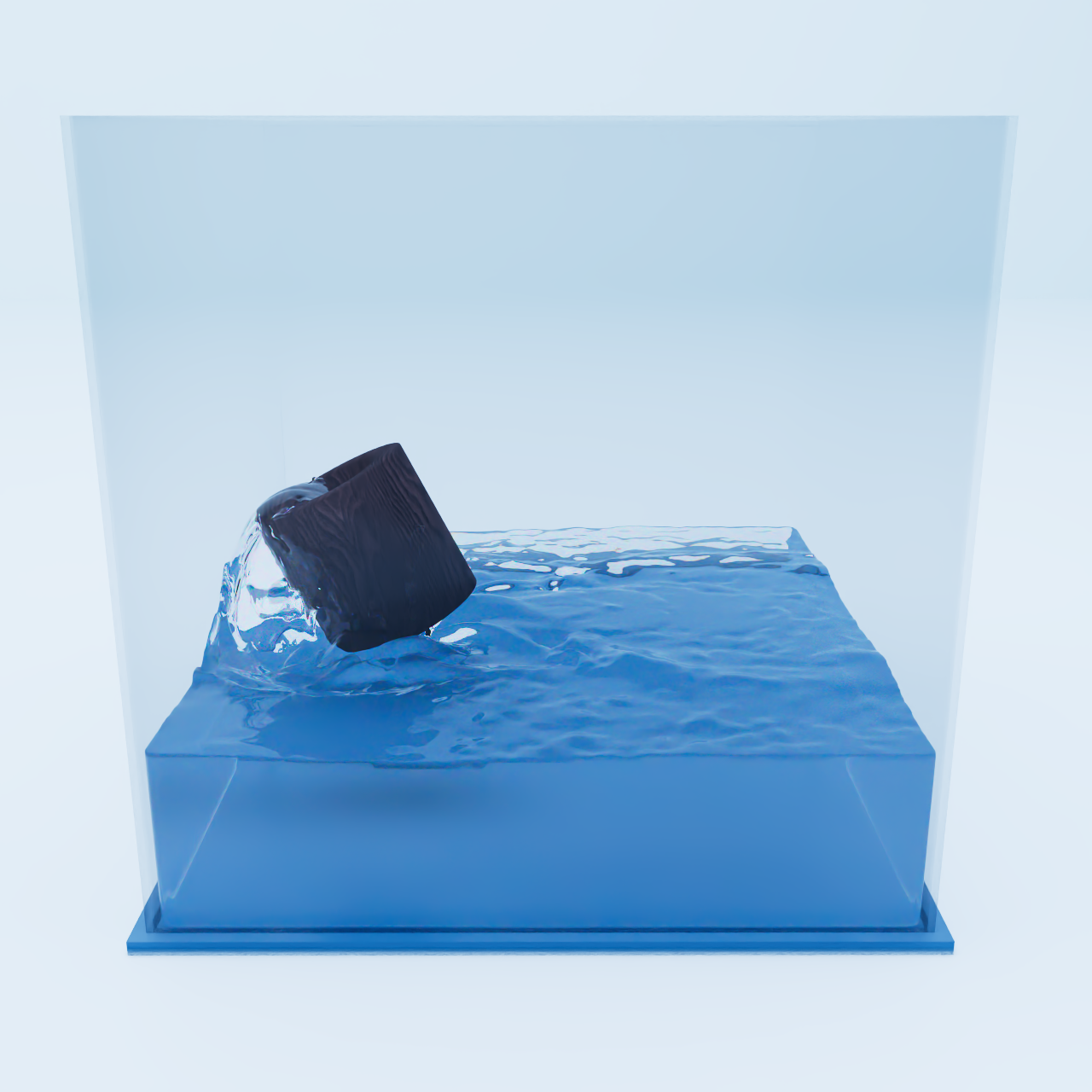}%
    \includegraphics[width=0.3\columnwidth]{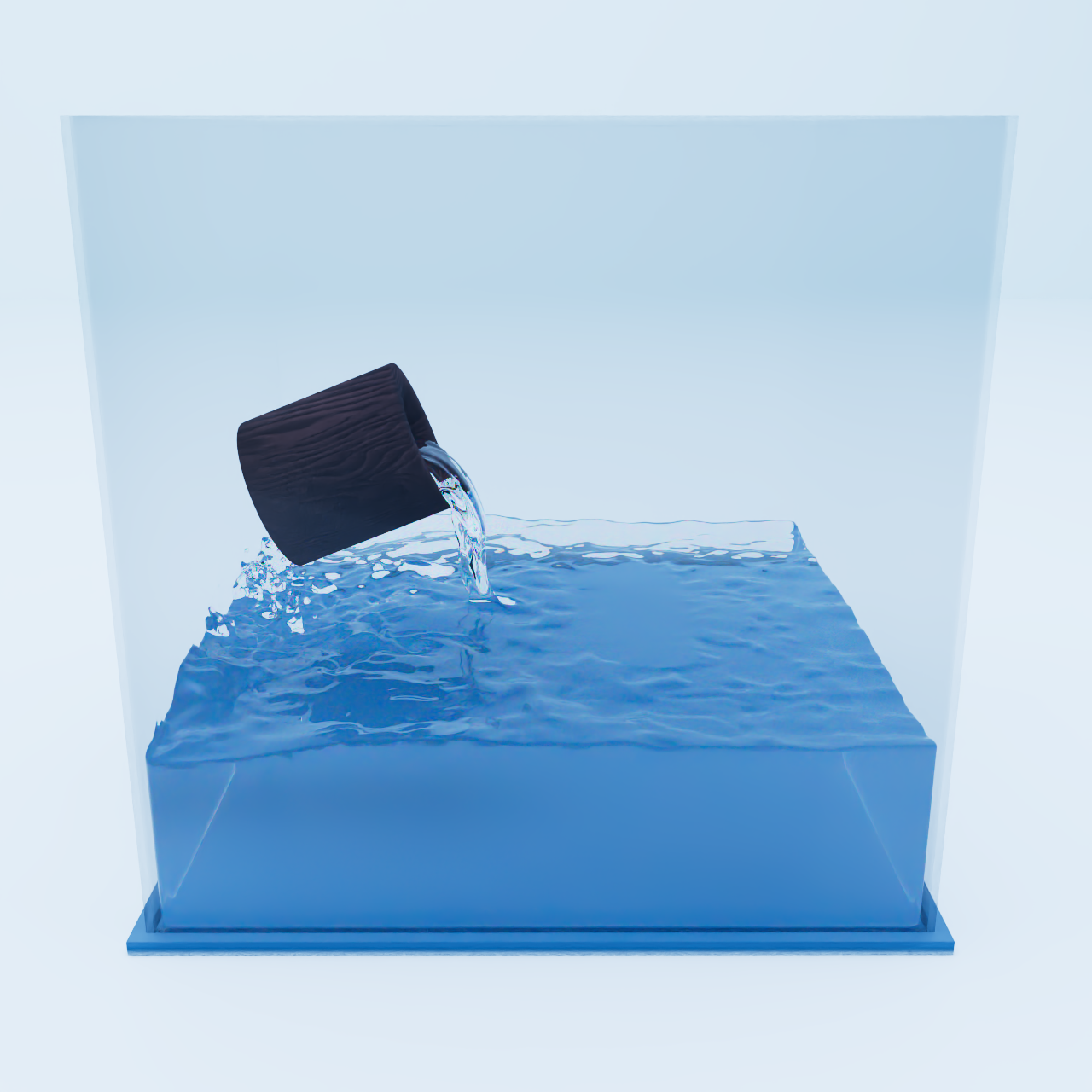} \\
    \includegraphics[width=0.3\columnwidth]{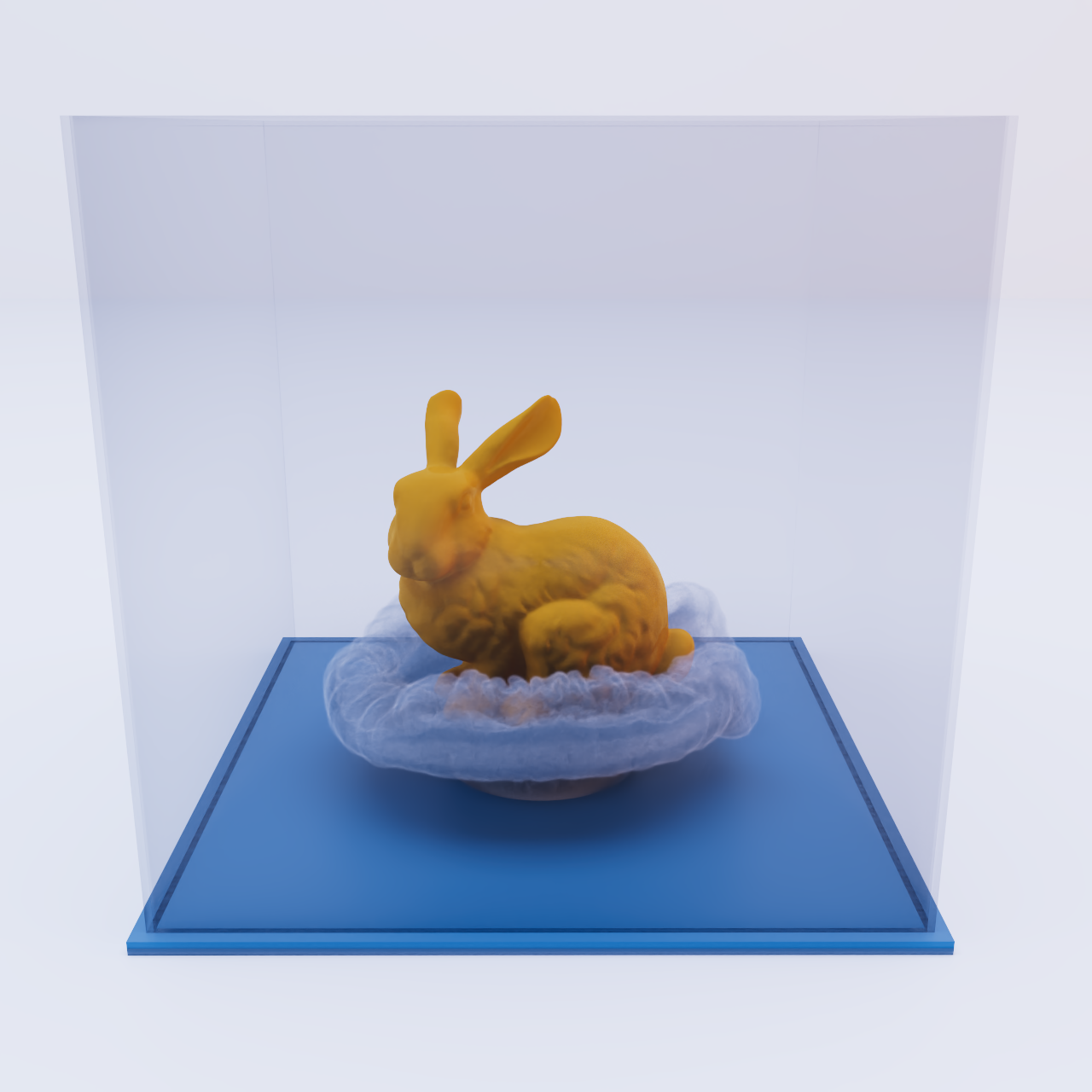}%
    \includegraphics[width=0.3\columnwidth]{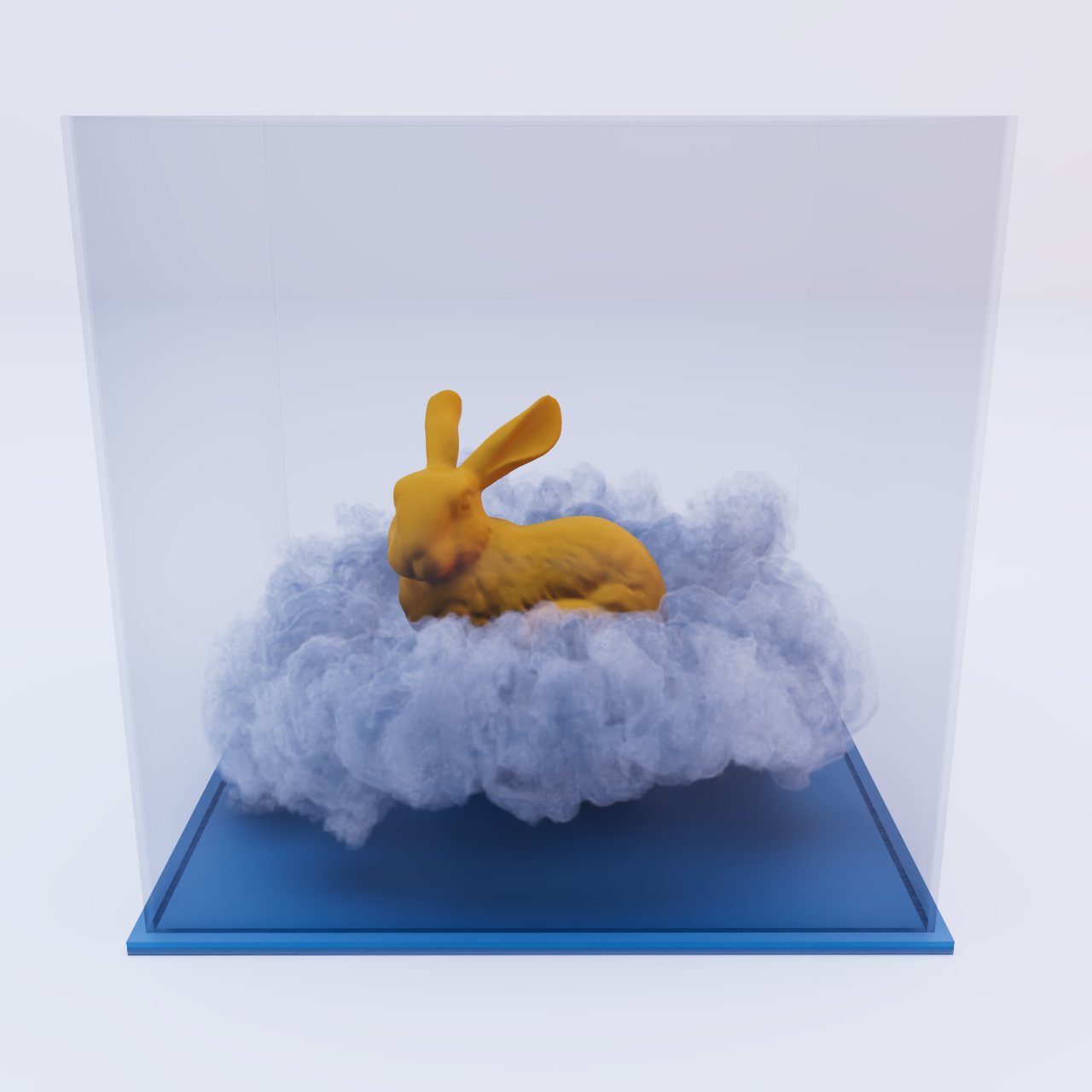}%
    \includegraphics[width=0.3\columnwidth]{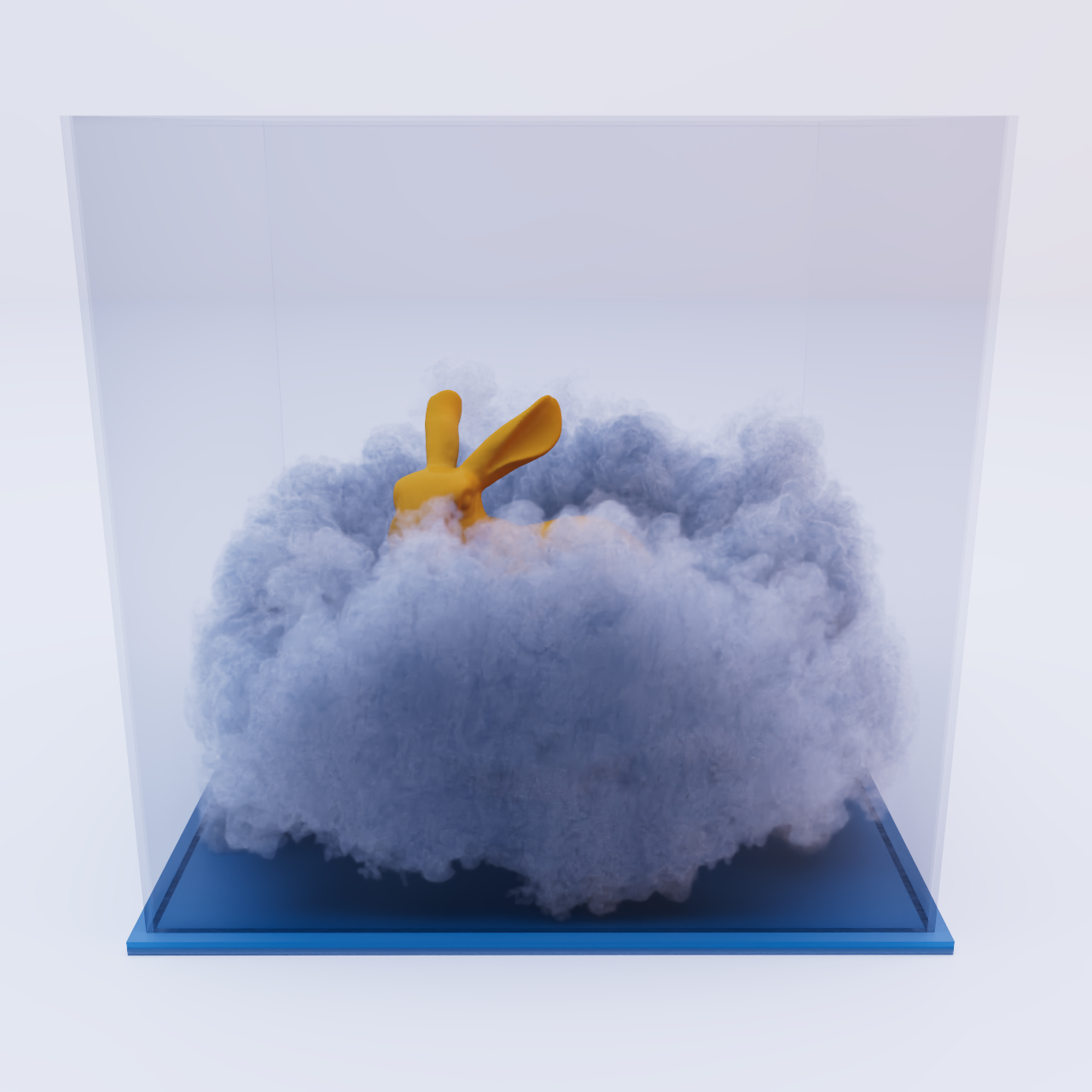}
   \caption{We demonstrate our solver with incompressible flow simulations requiring the solution of mixed Neumann/Dirichlet boundary conditions for the pressure Poisson equation.}
   \label{fig:examples}
\end{figure}

\subsubsection{Implementation}

   We built our network using PyTorch \citep{paszke:2019:pytorch},
   but implemented our convolutional and linear blocks as custom CUDA extensions.
   The neural network was trained using single-precision floating point.

\section{Results and Analysis}

We evaluate the effectiveness and efficiency of our neural preconditioned solver
by comparing it to high-performance state-of-the-art implementations
of several baseline methods: unpreconditioned CG provided by
the CuPy library \citep{okuta:2017:cupy}, as well as
CG preconditioned by the algebraic multigrid (AMG) and incomplete Cholesky (IC)
implementations from the AMGCL library \citep{demidov:2020:amgcl}.
We furthermore compare against NVIDIA's AmgX algebraic multigrid solver
\citep{naumov2015amgx}.
All of these baseline methods are accelerated by CUDA backends running on the GPU,
with the underlying IC implementation
coming from NVIDIA's cuSparse library.
Where appropriate, we also
compare against past neural preconditioners FluidNet \citep{tompson:2017:accelerating} and
DCDM \citep{kanedo:2023:dcdm}.
Finally, we include characteristic performance statistics of
a popular sparse Cholesky solver CHOLMOD \citep{chen:2008:cholmod}.
In all cases, our method outperforms these baselines, often dramatically.
We emphasize that a \emph{single} trained $\network$ instance is used throughout, demonstrating
its capacity to generalize to across simulation scenarios and grid shapes (see also \pr{apx:additional_domain_size}).

\begin{figure}
    \centering
      \subfigure[All systems]{\label{fig:hist_all}\includegraphics[width=0.45\columnwidth]{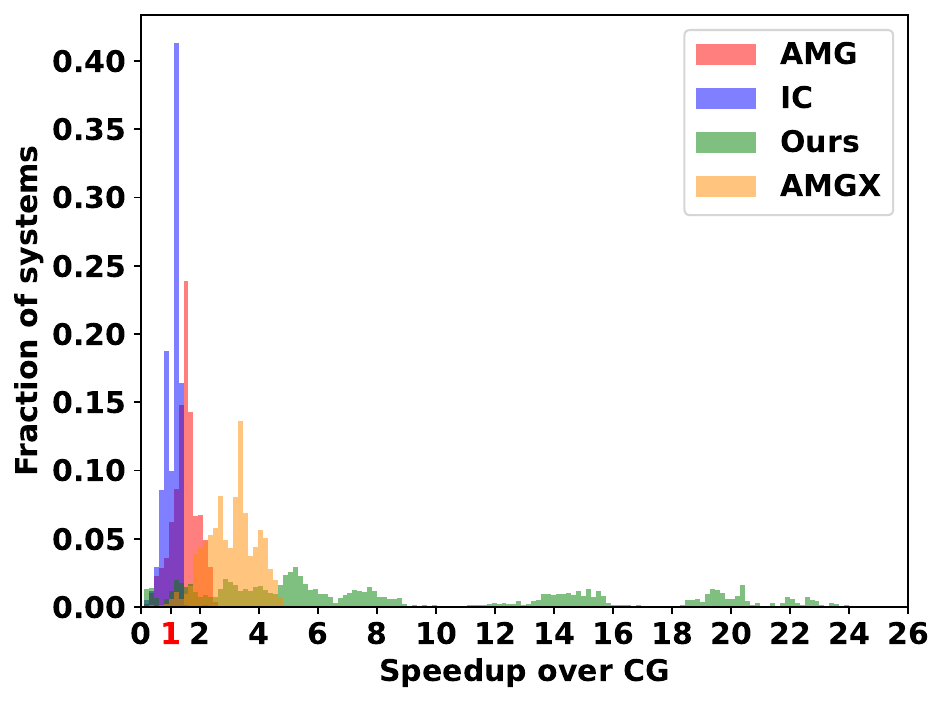}}
      \subfigure[Smallest 25\% of systems]{\label{fig:hist_all_low}\includegraphics[width=0.45\columnwidth]{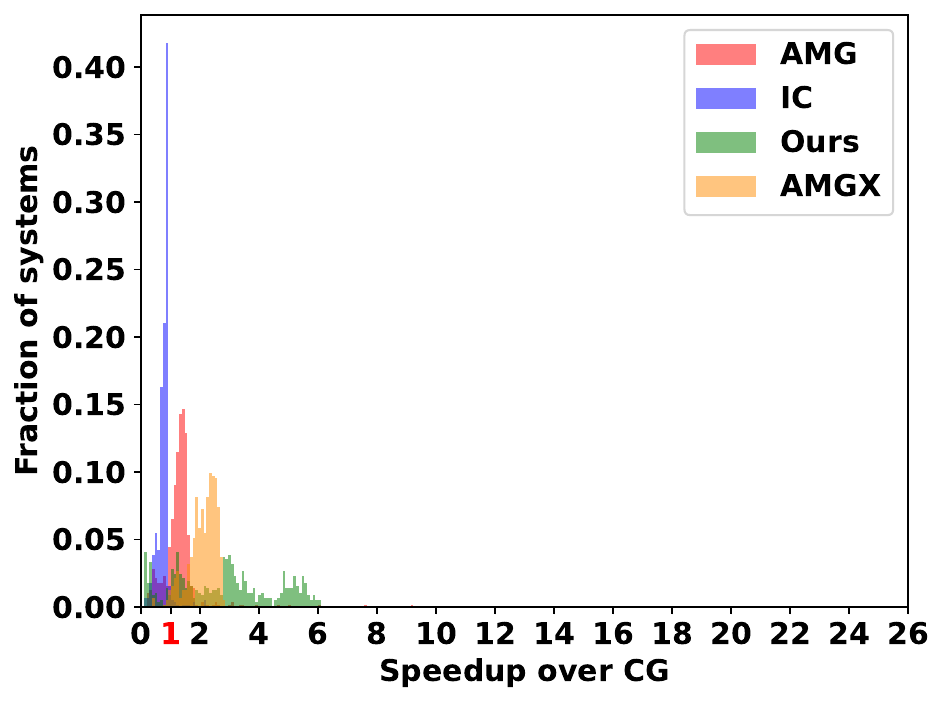}}
      \subfigure[Largest 75\% of systems]{\label{fig:hist_all_high}\includegraphics[width=0.45\columnwidth]{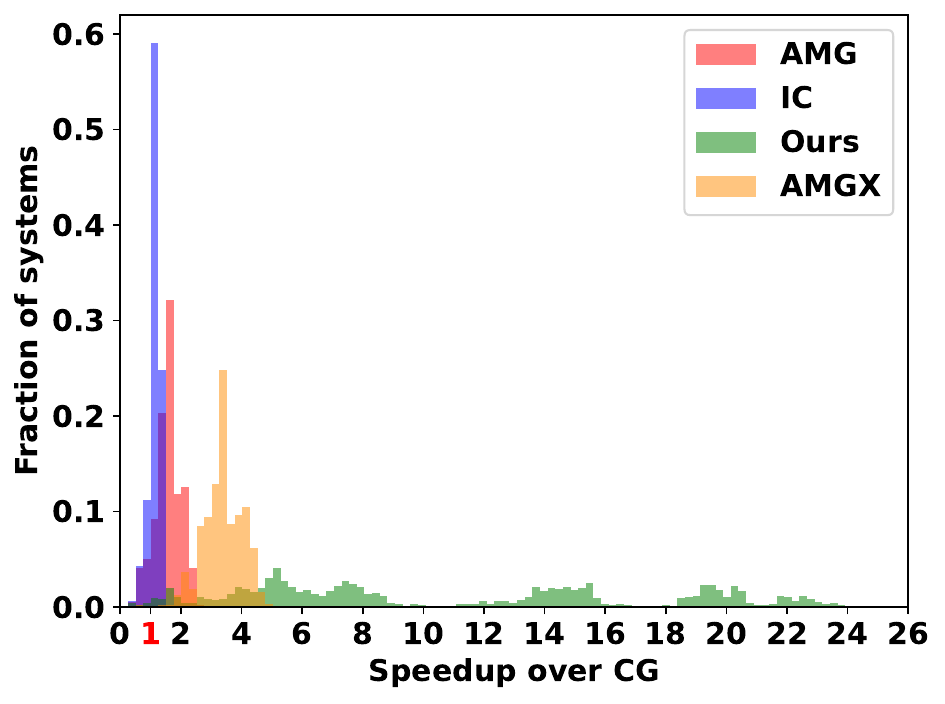}}
      \subfigure[Neumann only, $128^3$]{\label{fig:hist_smoke_128}\includegraphics[width=0.45\columnwidth]{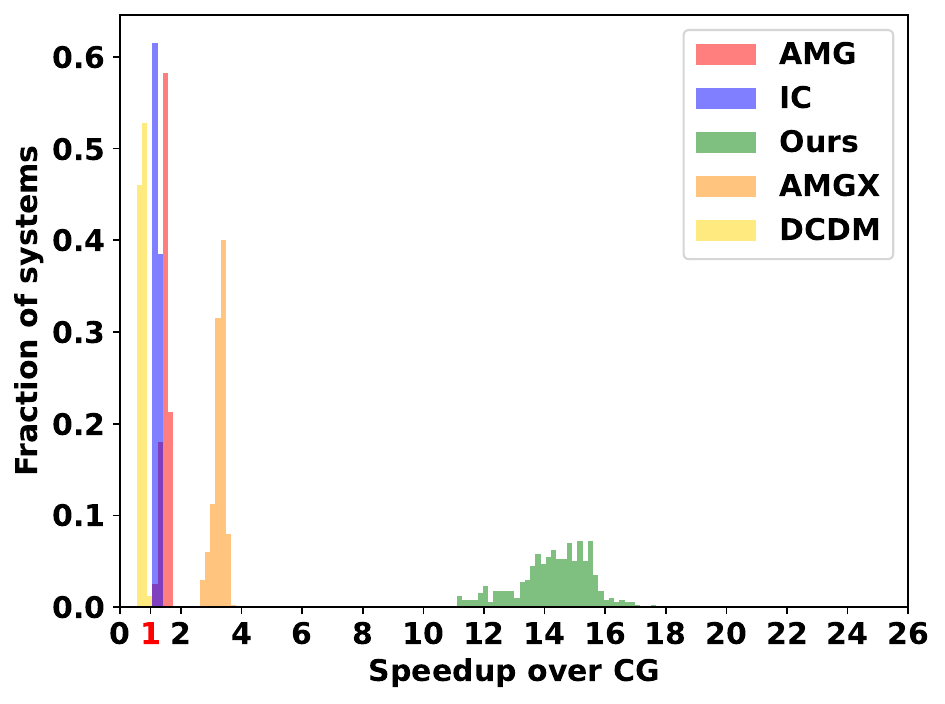}}
      \subfigure[Neumann only, $256^3$]{\label{fig:hist_smoke_256}\includegraphics[width=0.45\columnwidth]{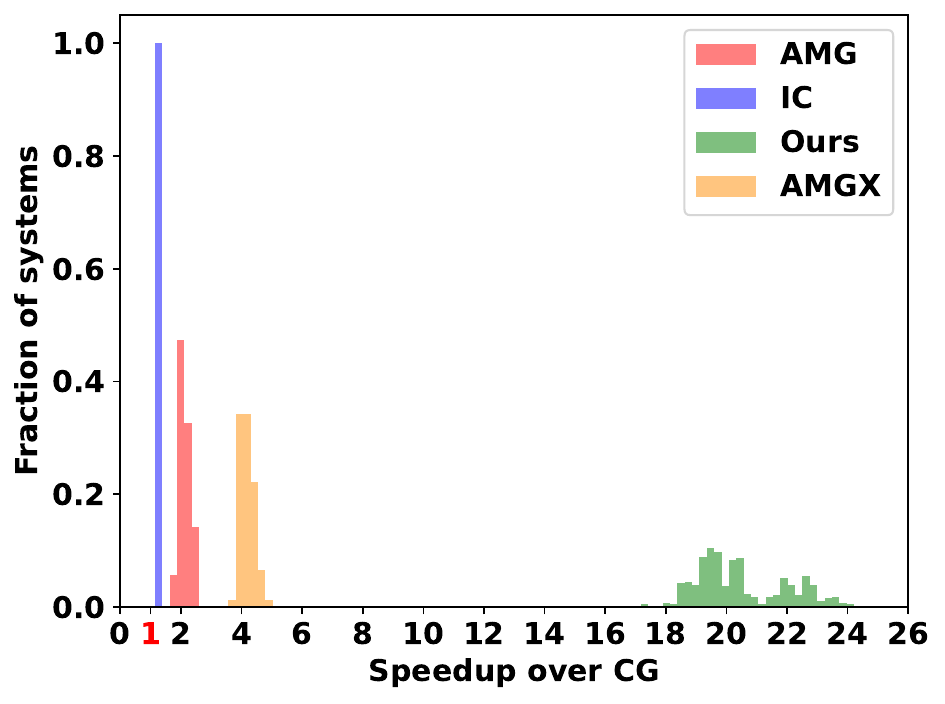}}
      \subfigure[Mixed, $128^3$ and $256^3$]{\label{fig:hist_water}\includegraphics[width=0.45\columnwidth]{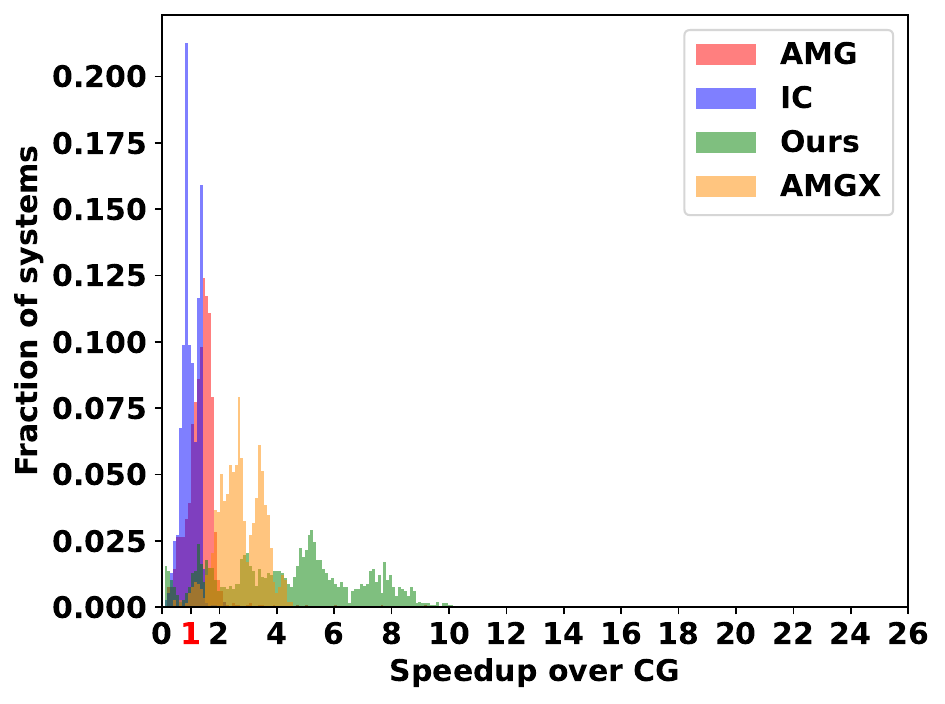}}
      \caption{Histograms of solution speedup vs. a baseline of unpreconditioned CG (a) for all solves;
      and (b-f) for certain subsets of the systems to help tease apart different modes of the distribution.}
      \label{fig:histograms}
  \end{figure}

We executed all benchmarks on a workstation featuring an AMD Ryzen 9 5950X 16-Core Processor and an NVIDIA GeForce RTX 3080 GPU.
We used as our convergence criterion for all methods a reduction of the residual norm by a factor of $10^6$, which is
sufficiently accurate to eliminate visible simulation artifacts.
We evaluate our neural preconditioner
in single precision floating point for efficiency
but implement the rest of the NPSDO algorithm in double precision
for numerical stability. We note that this use of single precision
does not limit the solver's overall accuracy, as demonstrated empirically in
\pr{apx:converge_test}.

\begin{figure}
      \centering
      \includegraphics[width=0.8\columnwidth]{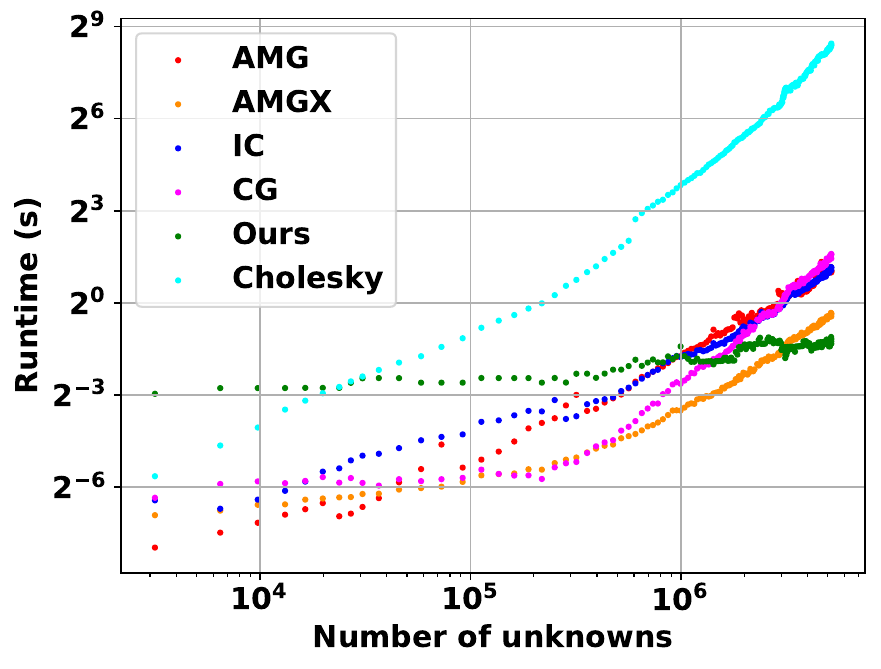}
      \caption{Solver scaling for mixed BC system matrices originating from a fixed-resolution domain $(n_c = 256^3)$; matrix row/col size $n_f$ is determined by the proportion of cells occupied by fluid. The vast majority of total solve time is contributed by the high-occupancy systems clustered to the right, where our method outperforms the rest.}
        \label{fig:cholesky_comparison}
      \end{figure}

   We benchmarked on twelve simulation scenes with various shapes---$(128,128,128)$, $(256,128,128)$, and $(256,256,256)$---each providing 200 linear systems
   to solve. For each solve, we recorded the number of iterations and runtime taken by each solver.
   These performance statistics are summarized visually in Figures~\ref{fig:histograms}-\ref{fig:cholesky_comparison} and \pr{apx:hist}, as well as in tabular form in \pr{apx:benchmarking}.

   \pr{fig:hist_all} summarizes timings from all solves in our benchmark suite:
   for each system, we divide the unpreconditioned CG solve time
   by the other methods' solve times to calculate their speedups
   and plot a histogram. We note that our method significantly outperforms
   the others on a majority of solves: ours is fastest
   on 95.6\% of the systems, which account for
   98.0\% of our total solve time.

   Our improvements are more substantial on larger problems
   (\Cref{fig:hist_all_low,fig:hist_all_high}) for two reasons.
   First, condition numbers increase with size,
   impeding solvers without effective preconditioners;
   this is seen clearly by comparing
   results from two different resolutions (\Cref{fig:hist_smoke_128,fig:hist_smoke_256}).
   Second, the small matrices $\Areduced$ correspond to
   simulation grids with mostly non-fluid cells.
   While CG, AMGCL, AmgX, and IC timings shrink significantly
   as fluid cells are removed, our network's evaluation
   cost does not:
   it always processes all of $\fullDomain$ regardless of occupancy.
   This scaling behavior is visible in \pr{fig:cholesky_comparison}.

   Our speedups are also greater for examples with $\dirichletBoundary = \emptyset$. DCDM is applicable for these, and so we included in it
   \pr{fig:hist_smoke_128} (but not in \pr{fig:hist_smoke_256} due to the network overspilling GPU RAM). DCDM's failure to outperform CG and IC in these
   results, contrary to \citep{kanedo:2023:dcdm}, is due to the
   higher-performance CUDA-accelerated implementations of those baselines used in this work. With
   Dirichlet conditions (\pr{fig:hist_water}), our preconditioner is less effective, and yet we
   still outperform the rest on 93.27\% of frames, which account for 96.36\%
   of our total solve time. Statistics are not reported in this setting for
   DCDM and FluidNet, which struggle to reduce the residual (\pr{fig:res_128}).

   \begin{figure}
      \centering
      \includegraphics[width=0.45\columnwidth]{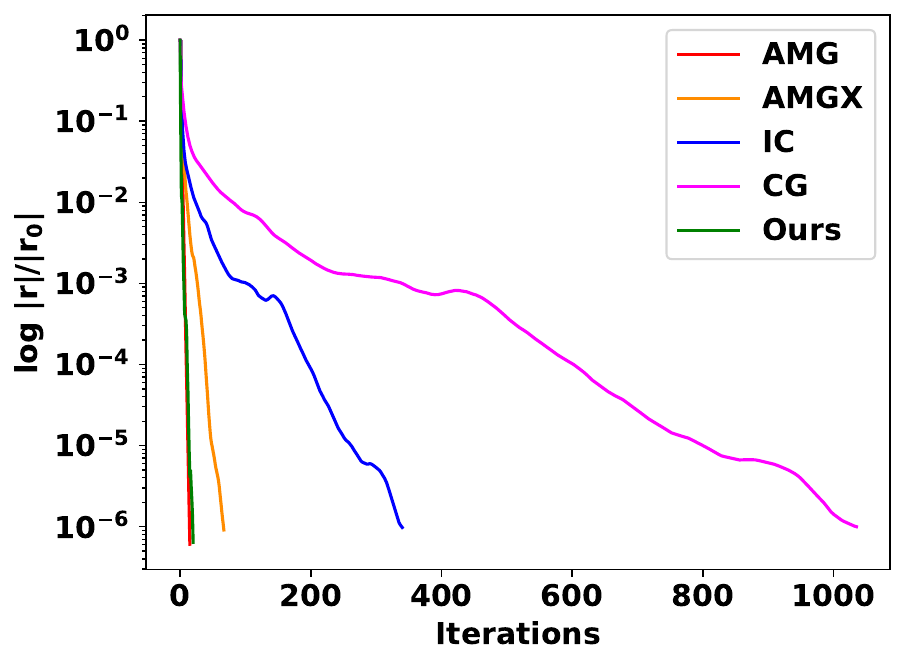}
      \includegraphics[width=0.45\columnwidth]{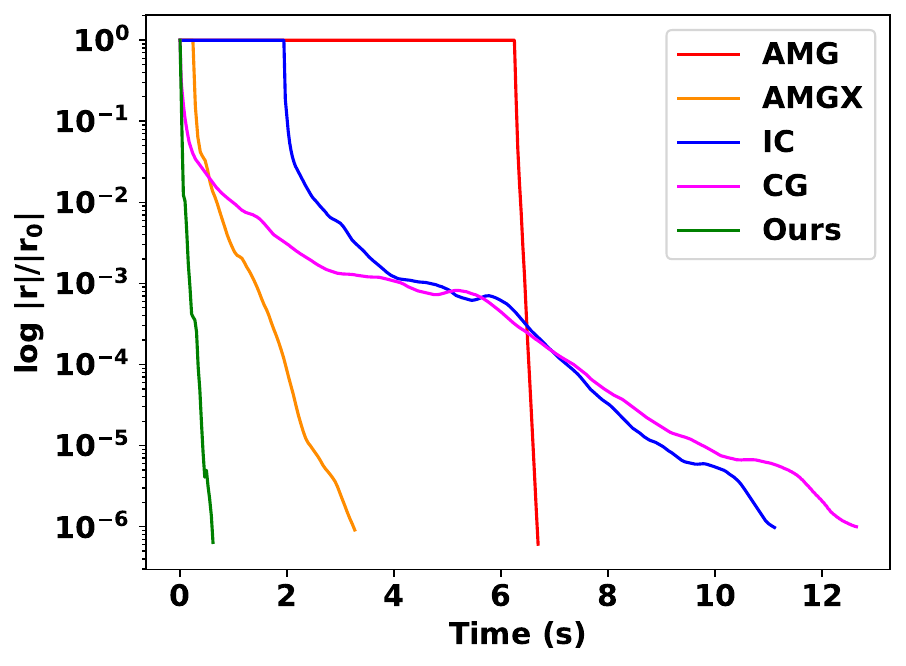}
      \includegraphics[width=0.45\columnwidth]{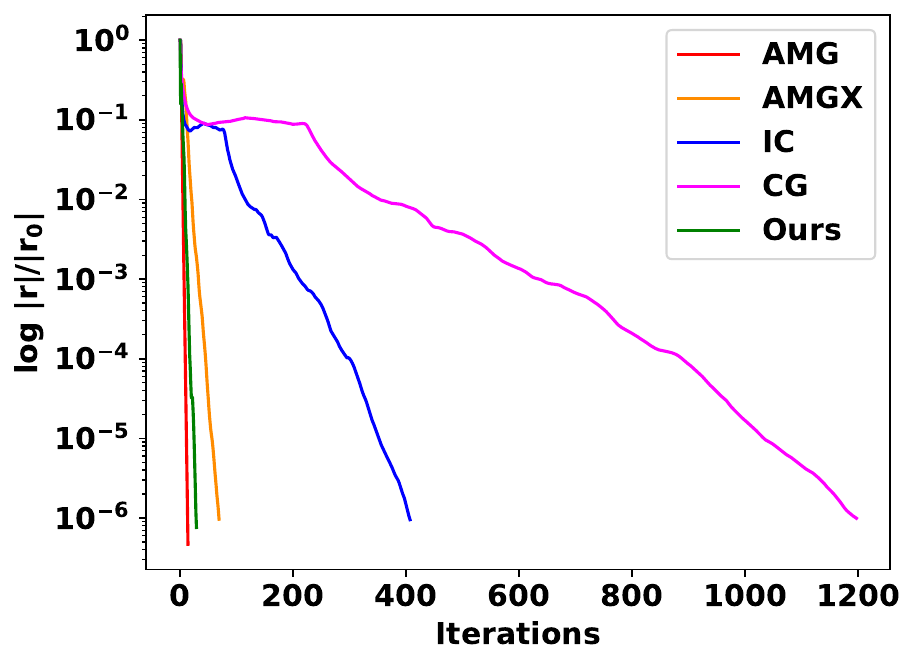}
      \includegraphics[width=0.45\columnwidth]{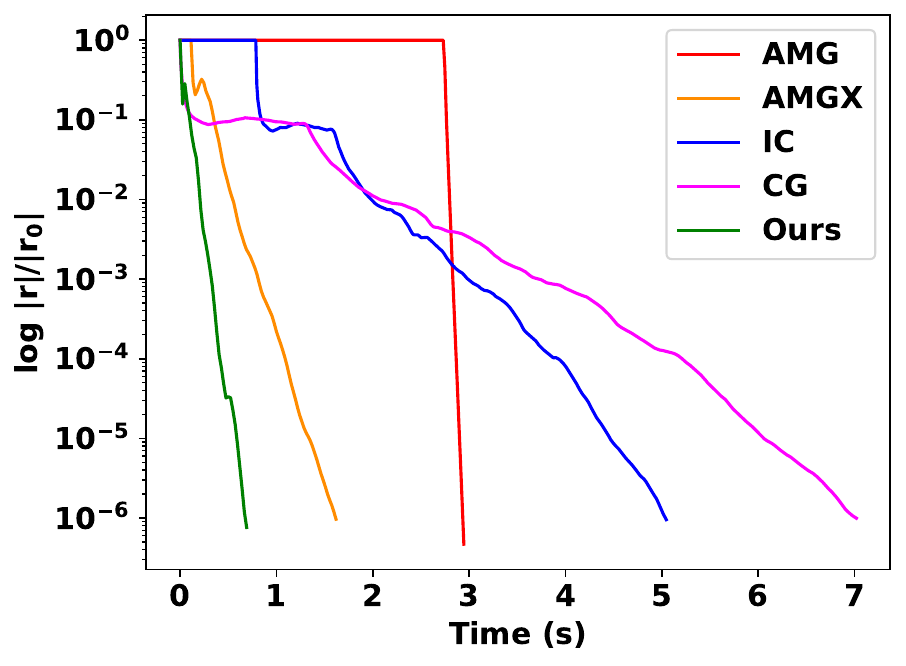}
    \caption{Comparisons among AMG, IC, CG and NSPDO (Ours) on a single frame at $256^3$ with Neumann only BC (top two) and mixed BC (bottom two).}
    \label{fig:res_256}

       \includegraphics[width=0.45\columnwidth]{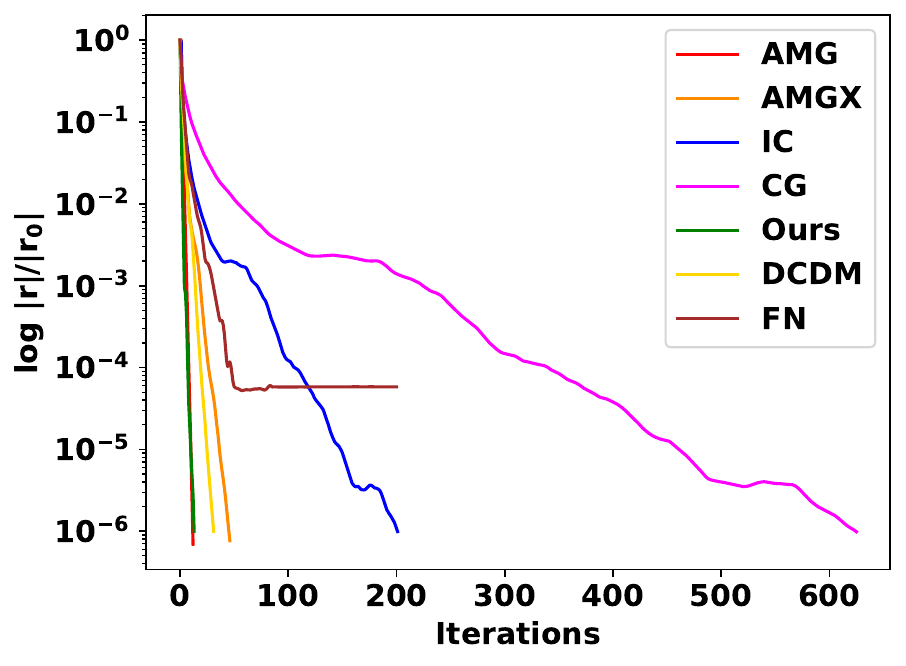}
       \includegraphics[width=0.45\columnwidth]{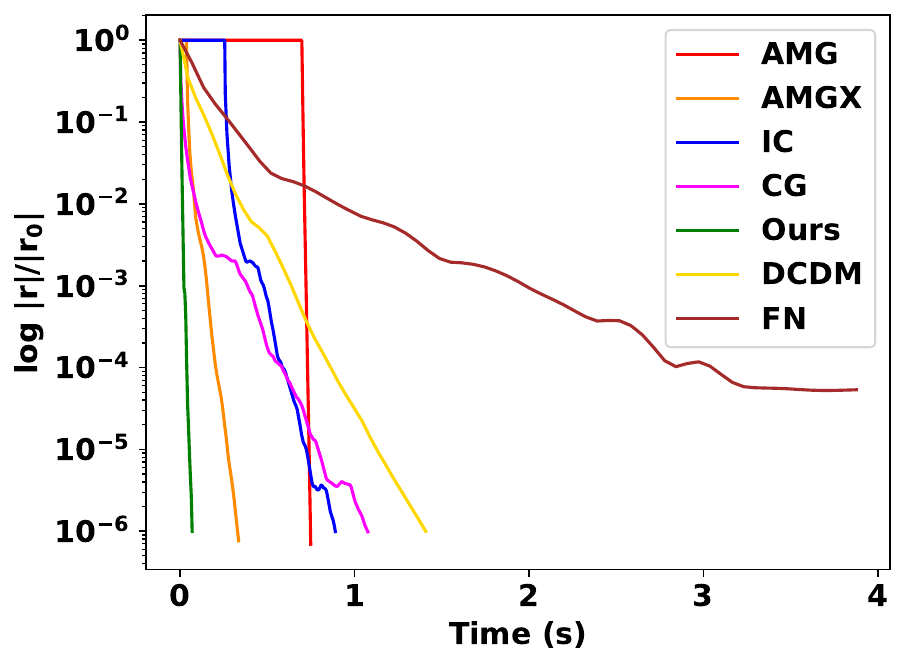}
       \includegraphics[width=0.45\columnwidth]{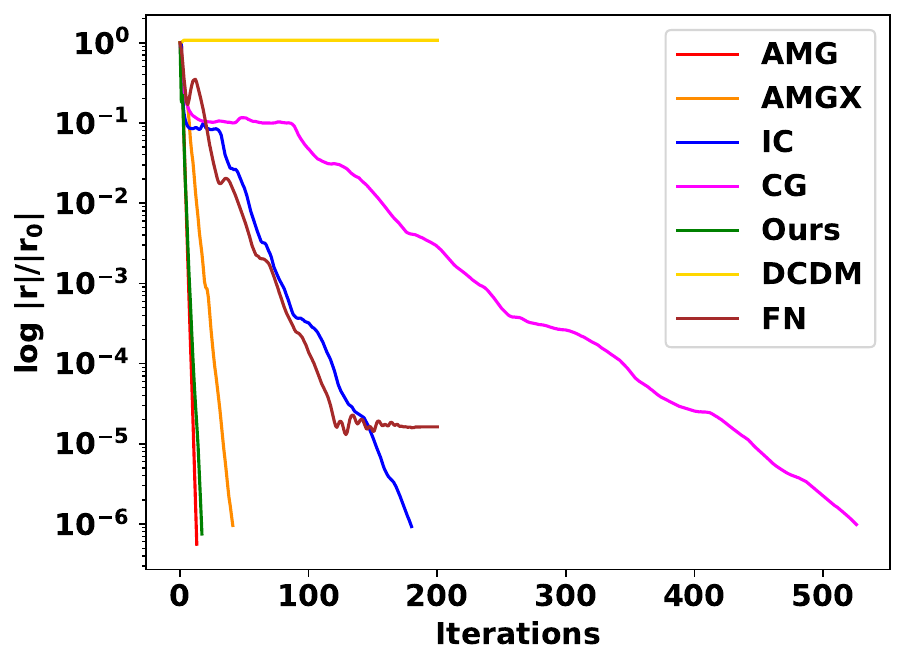}
       \includegraphics[width=0.45\columnwidth]{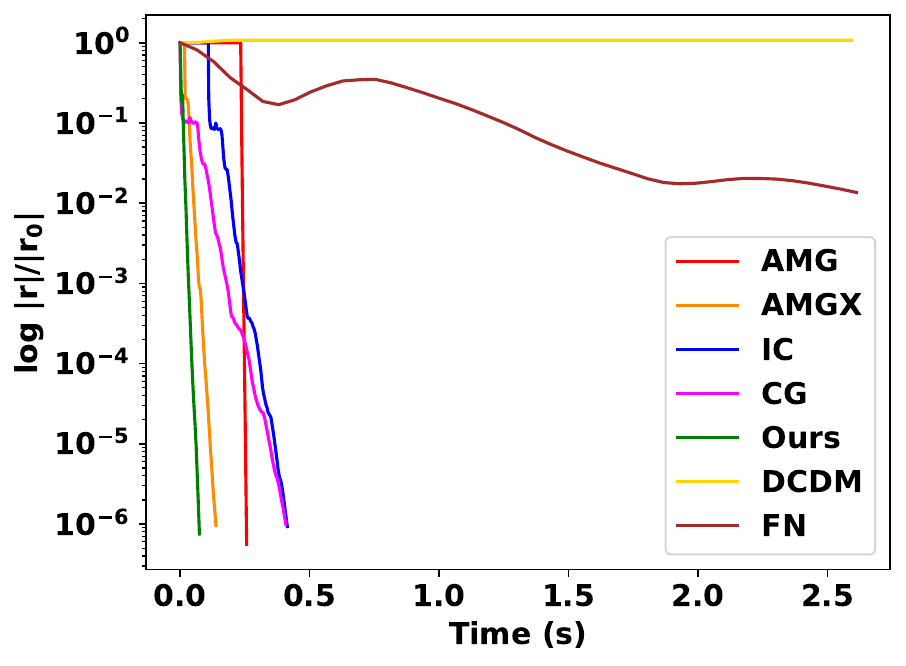}
    \caption{Comparisons among AMG, IC, CG, DCDM, FluidNet (FN) and NSPDO (Ours) on a single frame at $128^3$ with Neumann only BC (top two) and mixed BC (bottom two).}
    \label{fig:res_128}
  \end{figure}

   Further insights can be obtained by consulting \Cref{fig:res_128,fig:res_256},
   which show the convergence behavior of each iterative solver on characteristic example problems.
   AMG is clearly the most effective preconditioner, but this comes at the high cost
   of rebuilding the multigrid hierarchy before each solve:
   its iterations cannot even start until long after our solver already converged.
   Our preconditioner is the second most effective and, due to its lightweight architecture,
   achieves the fastest solves. DCDM is also quite effective at preconditioning Neumann-only problems
   but is slowed by costly network evaluations.
   IC's setup time is shorter than AMG but still
   substantial, and it is much less effective as a preconditioner; the same holds for AmgX.

   We note that the smoke example (\pr{fig:res_128}) also includes a comparison
   to FluidNet \emph{applied as a preconditioner} for PSDO.
   In the original paper, FluidNet was presented as a standalone solver, to be
   run just once per simulation frame.
   However, in this form it cannot produce highly accurate solutions.
   Incorporating it as a preconditioner as we do here in theory allows the system
   to be solved to controlled accuracy, but this solver
   ended up stalling before reaching a $10^6$ reduction in our experiments;
   it was omitted from \pr{fig:histograms} for this reason.

   On average,
   our solver spends 79.4\% of its time evaluating $\network$,
   4.4\% of its time in orthogonalization, and the remaining 16.2\%
   in other CG operations.
   In contrast, AMG takes a full 90\% of its time in its setup stage.
   IC's quicker construction and slower convergence mean it takes only 23\%
   in setup.
   Our architecture also confers GPU memory usage benefits: for $128^3$ grids,
   our solver uses 1.5GiB of RAM, while FluidNet and DCDM consume
   5.2GiB and 8.5GiB, respectively (\pr{apx:memory}).


   \subsection{Ablation Studies}
    To demonstrate the necessity of each component of our solver, we performed several ablation studies.
    Regarding the network architecture, we evaluated a simple CNN with fixed kernels and masked output, concatenating image and input vectors into separate channels. However, this approach performed poorly compared to our network. Specifically, the residual obtained by the masked CNN solver increased after the first few iterations and remained high even after reaching the $100$ iteration cap applied in this experiment. \pr{tbl:cnn_mask} summarizes this comparison.
We also investigated the effectiveness of several PCG variants, discussed in \pr{apx:pcgs}.
   \subsection{Hyperparameter Selection and More Experiments}
    The ideal number of network levels $\levels$ and orthogonalizations $\numOrtho$ were determined empirically as detailed in \pr{apx:levels} and \pr{apx:ortho}, respectively.

    To ensure fairness of our comparisons, we ran additional experiments on
    \href{https://www.dropbox.com/sh/5f3t9abmzu8fbfx/AAAkzW9JkkDshyzuFV0fAIL3a/bunny.capped.obj?e=1}{examples from FluidNet's original training dataset}
    processed according to the instructions provided with their \href{https://github.com/google/FluidNet}{open-source code release}.
    Naturally, these examples contain only pure-Neumann boundary conditions.
    We compared our model to DCDM and FluidNet, both as a standalone solver and
    as a preconditioner. For all iterative methods, except unpreconditioned
    Conjugate Gradient (CG), we capped the maximum iteration count at $100$.
    Notably, FluidNet, when used as a preconditioner, failed to converge to a
    residual tolerance of $10^{-6}$. As a standalone solver, FluidNet runs
    quickly but is highly ineffective at reducing the residual. These results
    are reported in \pr{tbl:fluidnet_data}.
    \begin{table}
        \caption{Comparison among our current network and CNN with masked output and CG} \label{tbl:cnn_mask}
        \begin{center}
            \scriptsize
            \begin{tabular}{l |*{3}{*{1}{c}}}
                \textbf{Methods} & \textbf{Iteration count} & \textbf{Time} & \textbf{Residual}\\\hline\\
                CG & $390$ & $0.1125$ & $9.72 \times10^{-7}$ \\\\
                Standard CNN & $100$ & $0.1025$ & $3.23\times 10^{-2}$\\\\
                Ours & $12$ & $0.027$ & $4.61\times 10^{-7}$
            \end{tabular}
        \end{center}
        \caption{Comparison among FluidNet (as a solver and as a preconditioner), DCDM, ours and CG}\label{tbl:fluidnet_data}
        \begin{center}
        \scriptsize
        \begin{tabular}{l |*{3}{*{1}{c}}}
            \textbf{Methods} & \textbf{Iteration count} & \textbf{Time} & \textbf{Final residual}\\\hline\\
            CG & $790$ & $0.5615$ & $1.00 \times 10^{-6}$\\\\
            Ours & $19$ & $0.05237$ & $8.77 \times 10^{-7}$\\\\
            FluidNet & $1$ & $0.03084$ & $2.035$\\\\
            FluidNet PSDO & $100$ & $3.116$ & $2.10 \times 10^{-5}$\\\\
            DCDM & $39$ & $12.35$ & $9.81 \times 10^{-7}$
        \end{tabular}
        \end{center}
    \end{table}

   \section{Conclusions}
   \label{sec:conclusions}
   The neural-preconditioned solver we propose not only addresses more
   general boundary conditions than past machine learning approaches for the
   Poisson equation \citep{tompson:2017:accelerating,kanedo:2023:dcdm} but also
   dramatically outperforms these solvers. It even surpasses state-of-the art
   high-performance implementations of standard methods like algebraic multigrid
   and incomplete Cholesky. It achieves this through a combination of its strong
   efficacy as a preconditioner and its fast evaluations enabled by our novel
   lightweight architecture.

   Nevertheless, we see several opportunities to improve and extend our solver in
   future work. Although we implemented our spatially-varying convolution
   block in CUDA, it remains the computational bottleneck of the network evaluation
   and is not yet fully optimized. We are excited to try porting our
   architecture to special-purpose acceleration hardware like Apple's Neural Engine;
   not only could this offer further speedups, but also it would
   free up GPU cycles for rendering the results in real-time applications like
   visual effects and games.
   For applications where fluid occupies only a small portion
   of the computational domain, we would like to develop techniques
   to exploit sparsity for better scaling (\pr{fig:cholesky_comparison}).
   Finally, we look forward to extending our ideas to achieve
   competitive performance for problems posed on unstructured
   grids as well as equations with non-constant coefficients, vector-valued unknowns
   (\emph{e.g.}, elasticity), and nonlinearities.


\clearpage
\section*{Acknowledgements}
Kai Lan was supported on a UC Multiple Campus Award (MCA) MCA-001639-M23PL6076.
Any opinions, findings, conclusions, or recommendations expressed are those of the authors and do not necessarily reflect the views of the sponsor.

\section*{Impact Statement}
This paper aims to advance the field of computational fluid simulation. We do not foresee any potential societal consequences resulting from our work that warrant discussion.

\bibliography{references.bib}
\bibliographystyle{icml2024}

\newpage
\appendix
\onecolumn
\section{Appendix}
\subsection{Additional Histograms}\label{apx:hist}

The following histograms offer additional views into the data presented in \pr{fig:histograms},
focusing on the linear systems arising from simulations with mixed boundary conditions (\emph{i.e.}, featuring
both Dirichlet and Neumann conditions).

\begin{figure}[h!]
   \centering
   \subfigure[Smoke solid, $(128,128,128)$]{\label{}\includegraphics[width=0.3\textwidth]{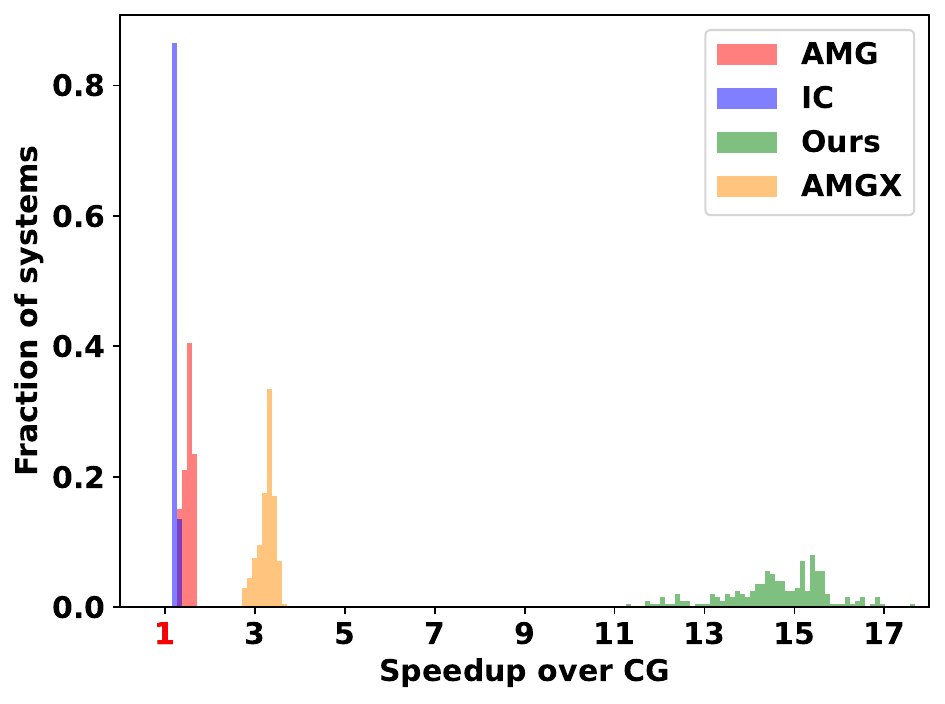}}
    \subfigure[smoke bunny, $(128,128,128)$]{\label{}\includegraphics[width=0.3\textwidth]{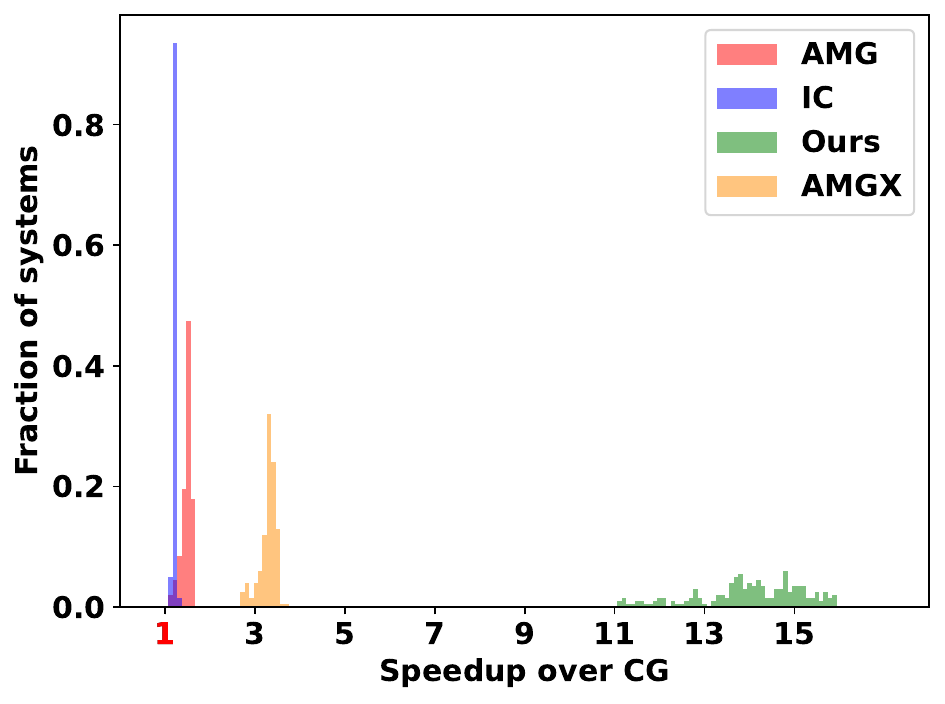}}
    \subfigure[Scooping, $(128,128,128)$]{\label{}\includegraphics[width=0.3\textwidth]{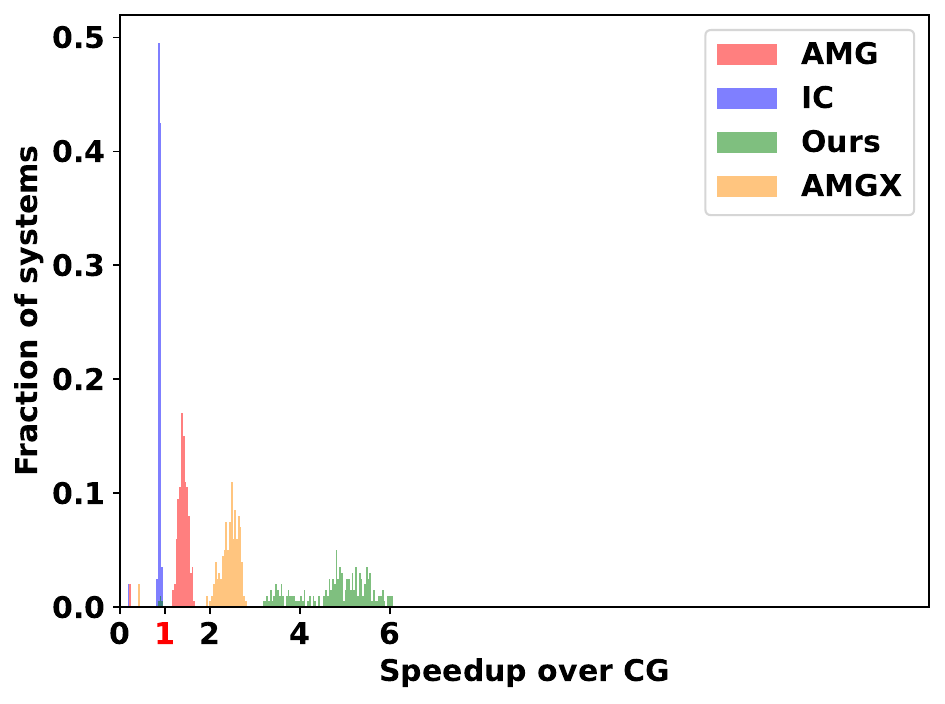}}
    \subfigure[Waterflow torus, $(128,128,128)$]{\label{}\includegraphics[width=0.3\textwidth]{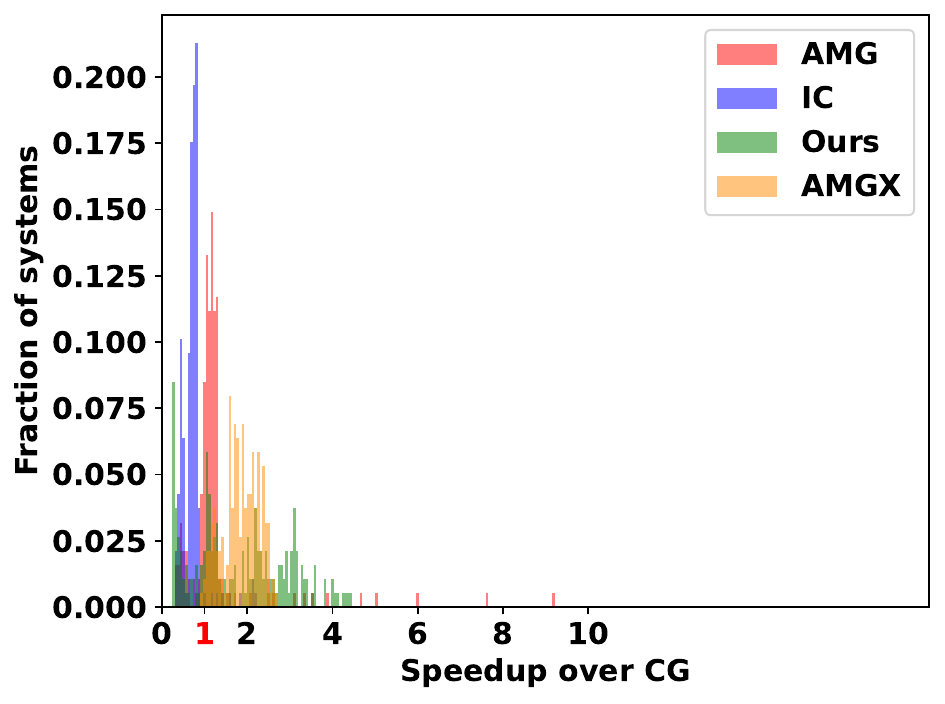}}
    \subfigure[Waterflow ball, $(128,128,128)$]{\label{}\includegraphics[width=0.3\textwidth]{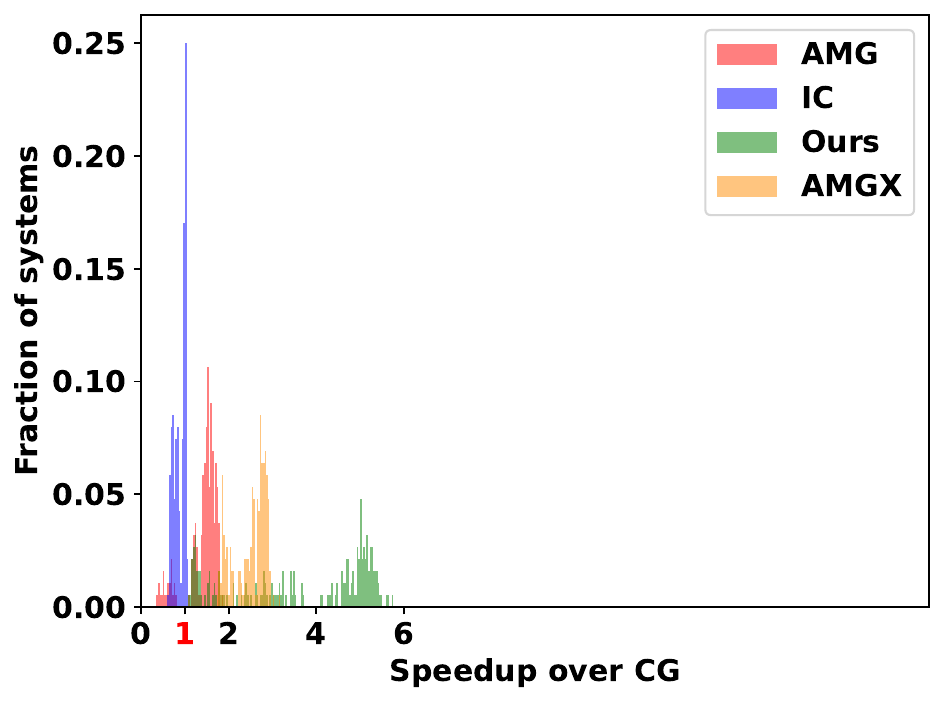}}
    \subfigure[Dambreak pillars, $(256,128,128)$]{\label{}\includegraphics[width=0.3\textwidth]{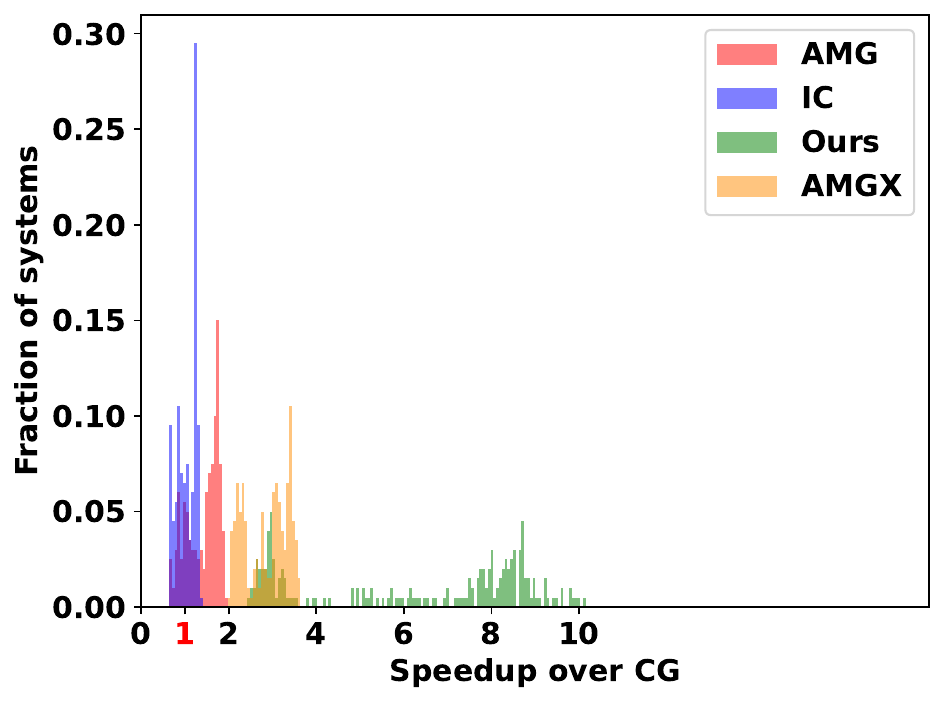}}
    \subfigure[Dambreak bunny, $(256,128,128)$]{\label{}\includegraphics[width=0.3\textwidth]{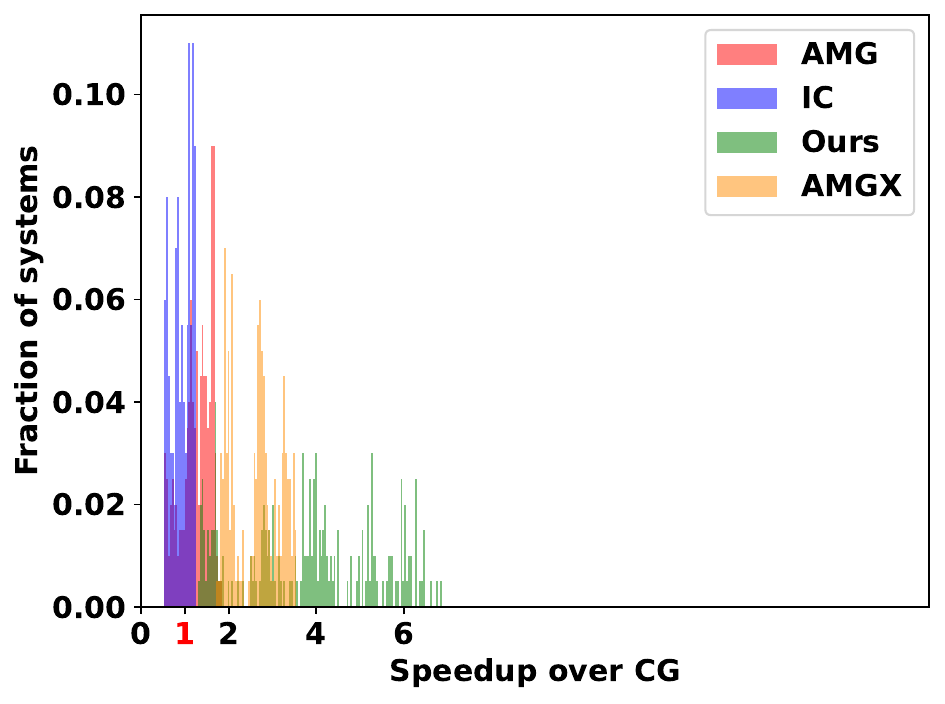}}
    \subfigure[Smoke solid, $(256,256,256)$]{\label{}\includegraphics[width=0.3\textwidth]{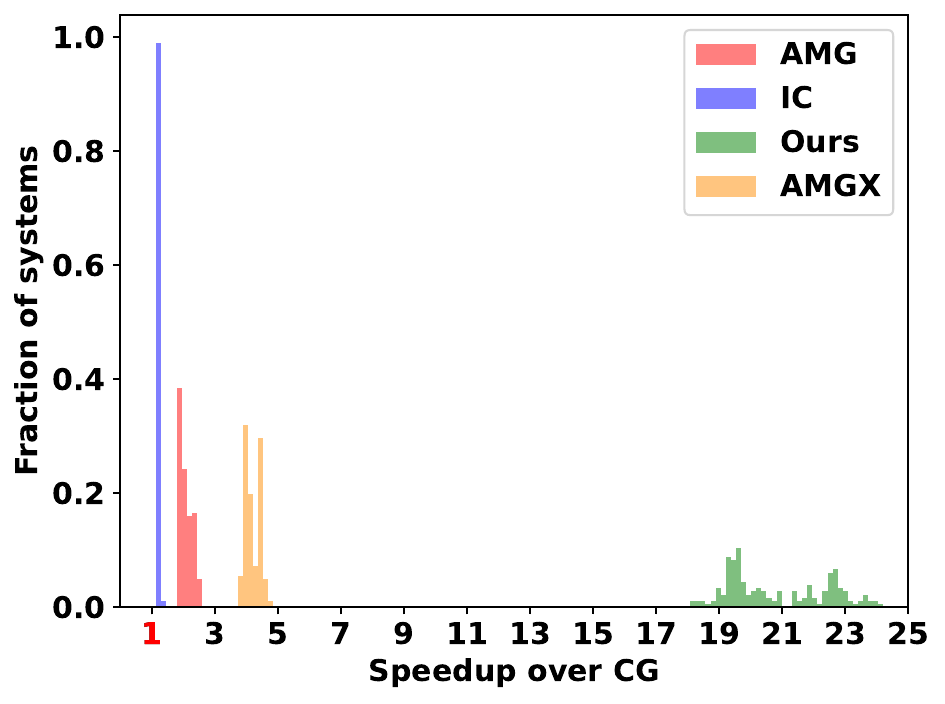}}
    \subfigure[Smoke bunny, $(256,256,256)$]{\label{}\includegraphics[width=0.3\textwidth]{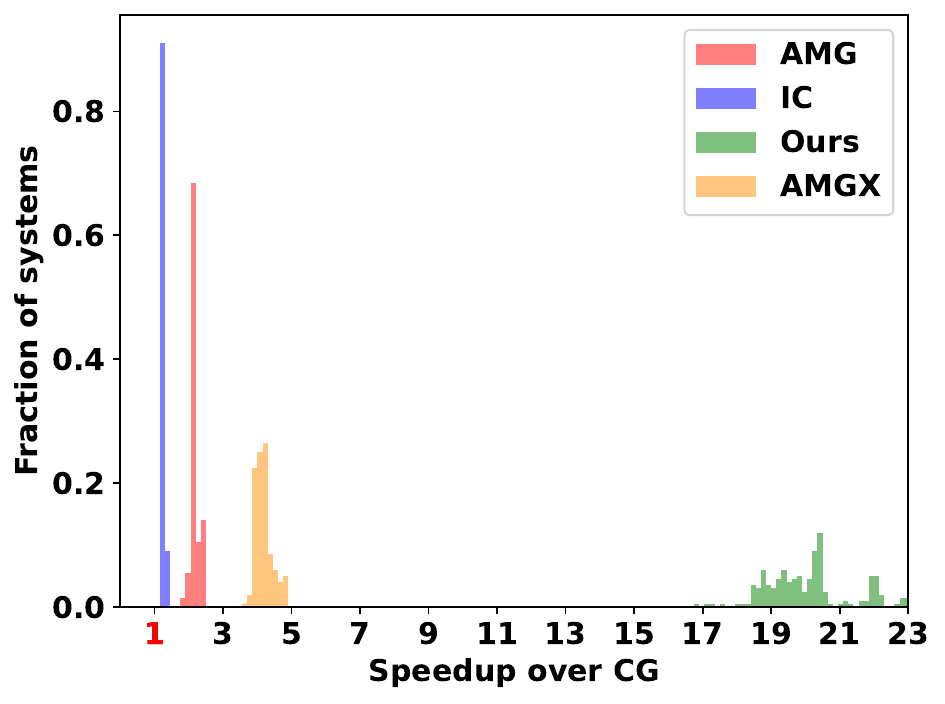}}
    \subfigure[Scooping, $(256,256,256)$]{\label{}\includegraphics[width=0.3\textwidth]{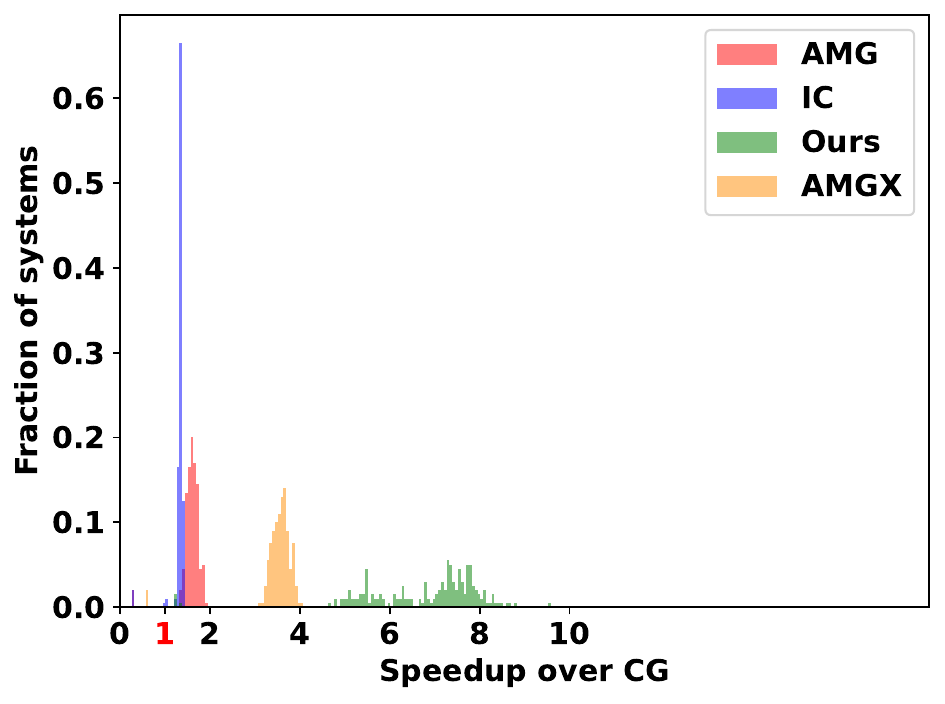}}
    \subfigure[Waterflow torus, $(256,256,256)$]{\label{}\includegraphics[width=0.3\textwidth]{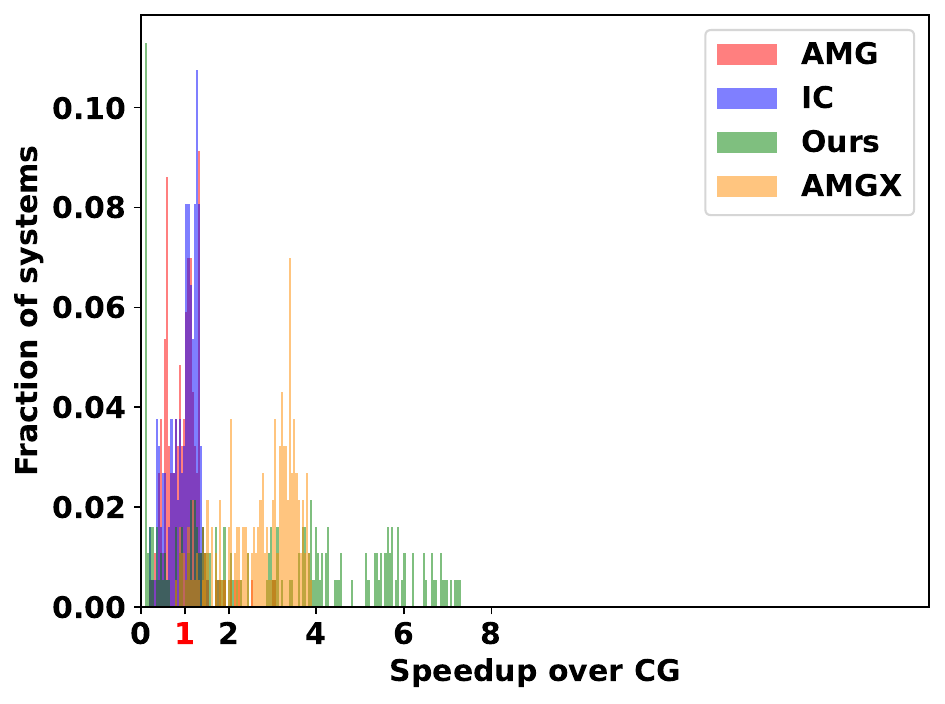}}
    \subfigure[Waterflow ball, $(256,256,256)$]{\label{}\includegraphics[width=0.3\textwidth]{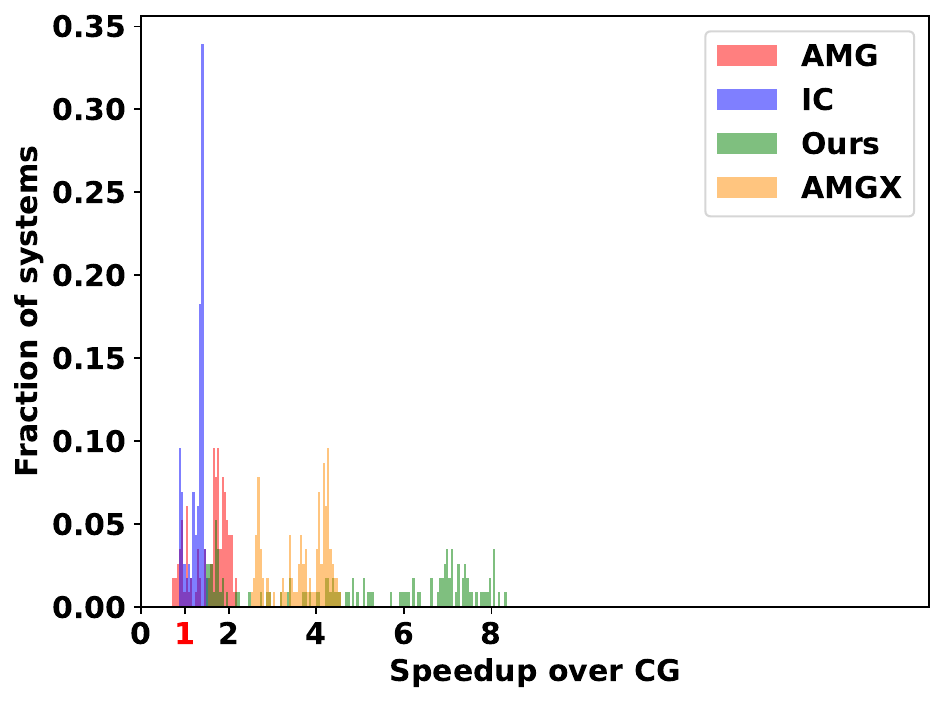}}
   \caption{Mixed BC examples split by linear system size into smallest $25\%$ and largest $75\%$.}
\end{figure}

\newpage
\subsection{Analysis on depth of the network}\label{apx:levels}

The depth $\levels$ of our network model is a hyperparameter. To study its impact on the performace of our solver, we compared the models for $\levels=3,4,5,6$ on a few examples.

\begin{figure}[h!]
    \centering
    \subfigure[Smoke solid, $(128,128,128)$]{\label{}\includegraphics[width=0.3\textwidth]{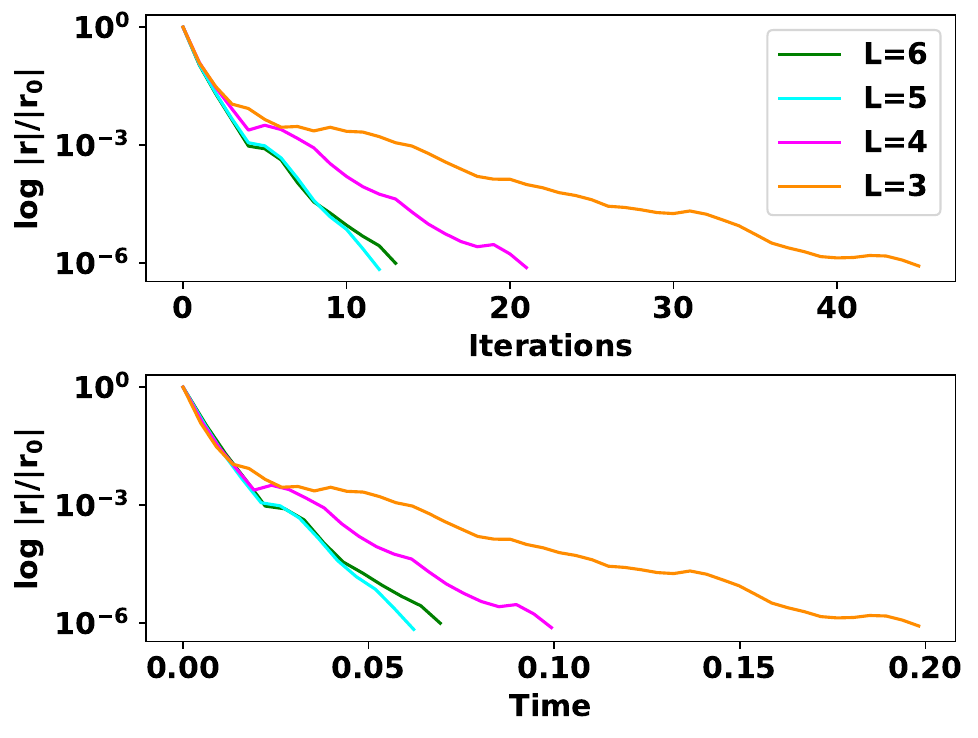}}
    \subfigure[smoke bunny, $(128,128,128)$]{\label{}\includegraphics[width=0.3\textwidth]{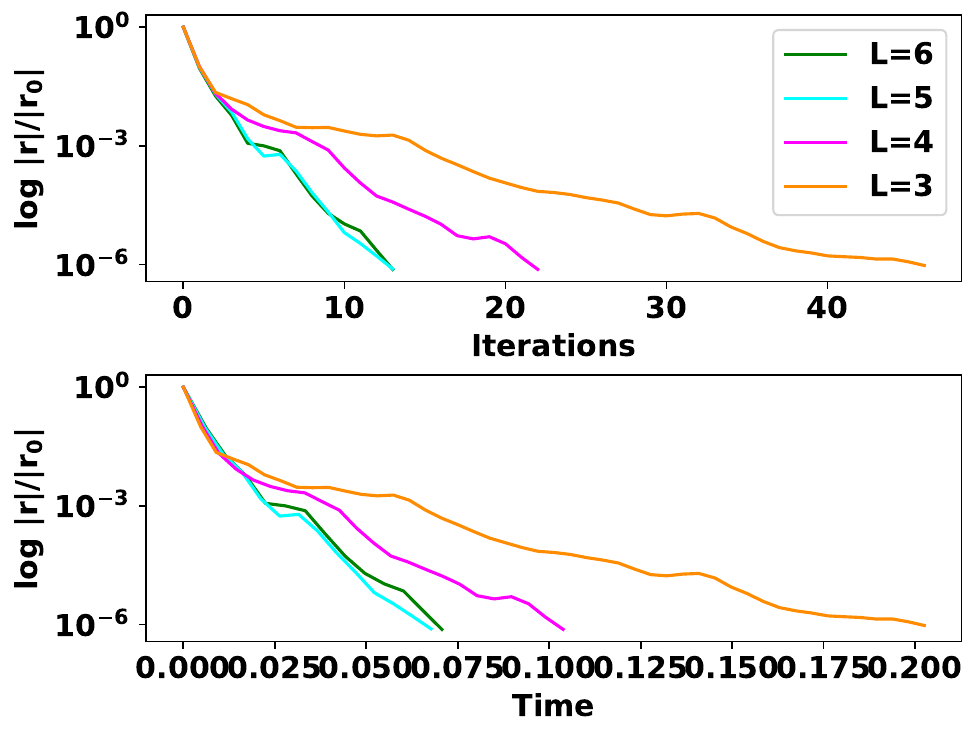}}
    \subfigure[Scooping, $(128,128,128)$]{\label{}\includegraphics[width=0.3\textwidth]{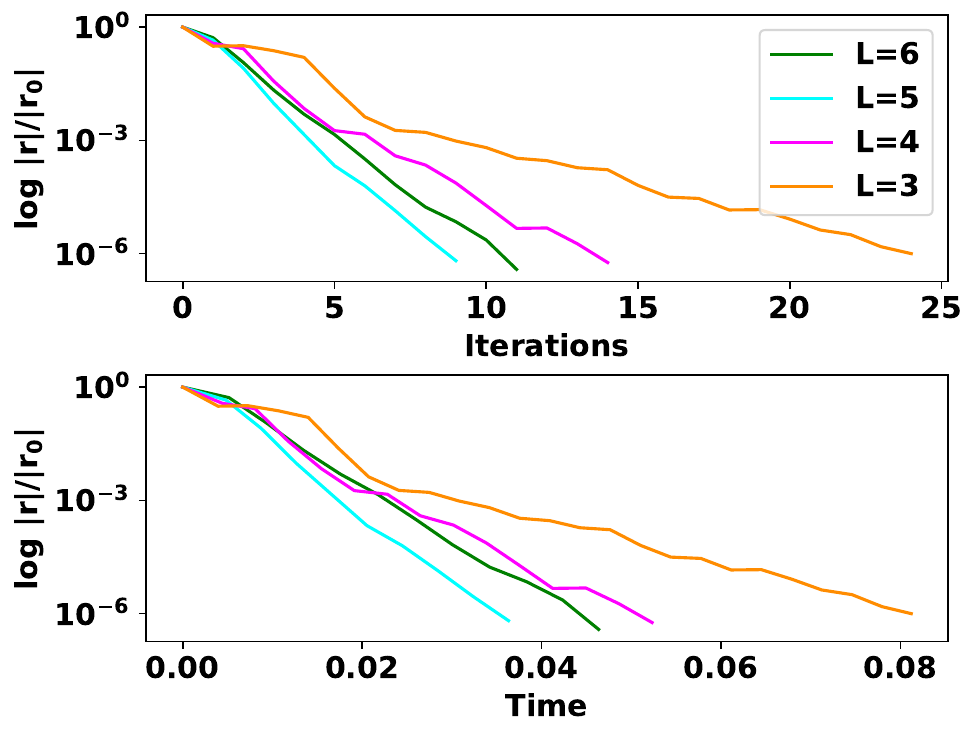}}
    \subfigure[Waterflow torus, $(128,128,128)$]{\label{}\includegraphics[width=0.3\textwidth]{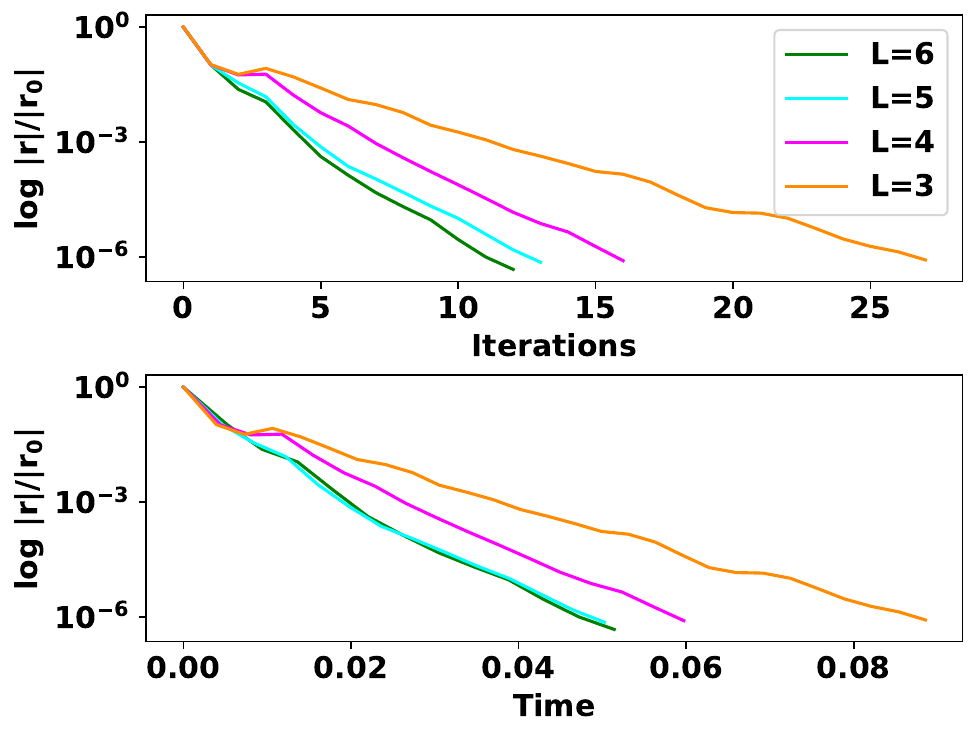}}
    \subfigure[Waterflow ball, $(128,128,128)$]{\label{}\includegraphics[width=0.3\textwidth]{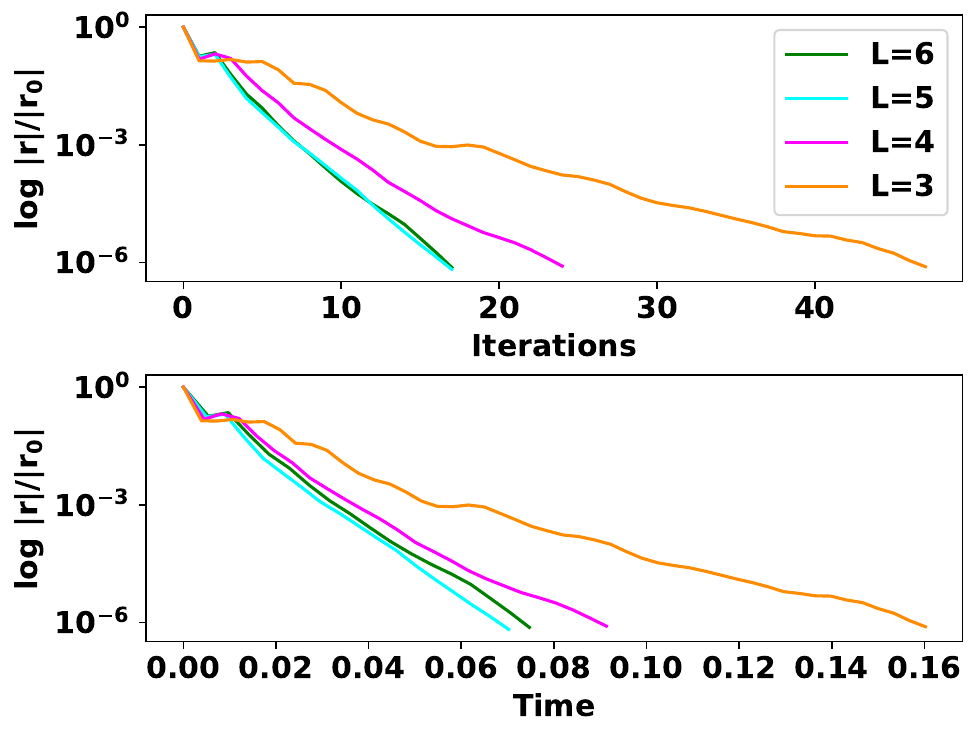}}
    \subfigure[Dambreak pillars, $(256,128,128)$]{\label{}\includegraphics[width=0.3\textwidth]{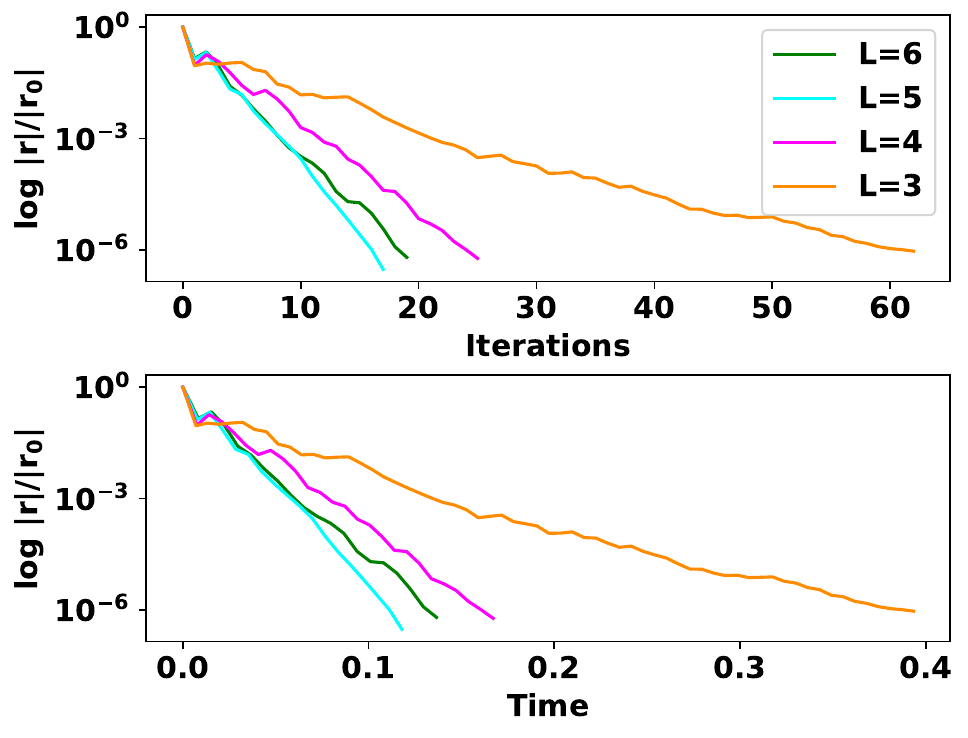}}
    \subfigure[Dambreak bunny, $(256,128,128)$]{\label{}\includegraphics[width=0.3\textwidth]{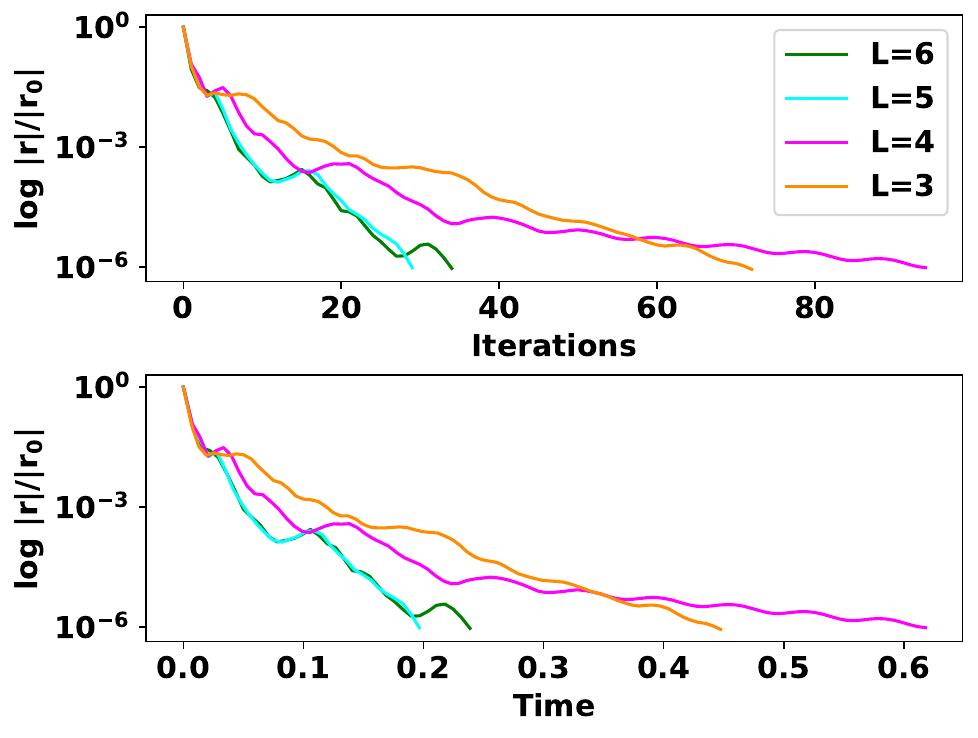}}
    \subfigure[Smoke solid, $(256,256,256)$]{\label{}\includegraphics[width=0.3\textwidth]{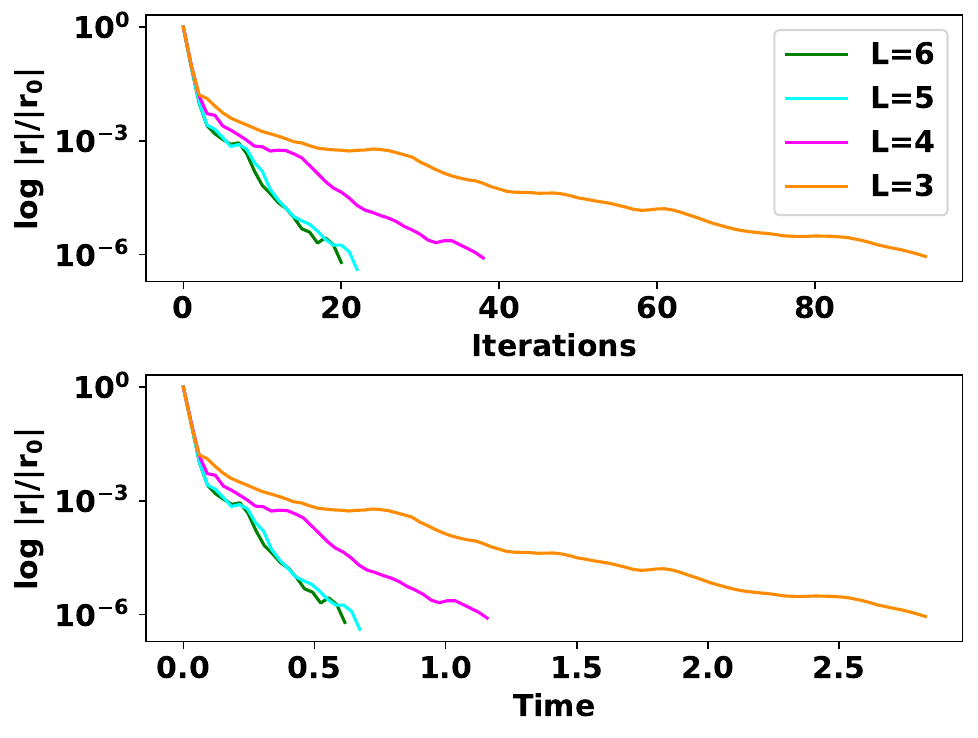}}
    \subfigure[Smoke bunny, $(256,256,256)$]{\label{}\includegraphics[width=0.3\textwidth]{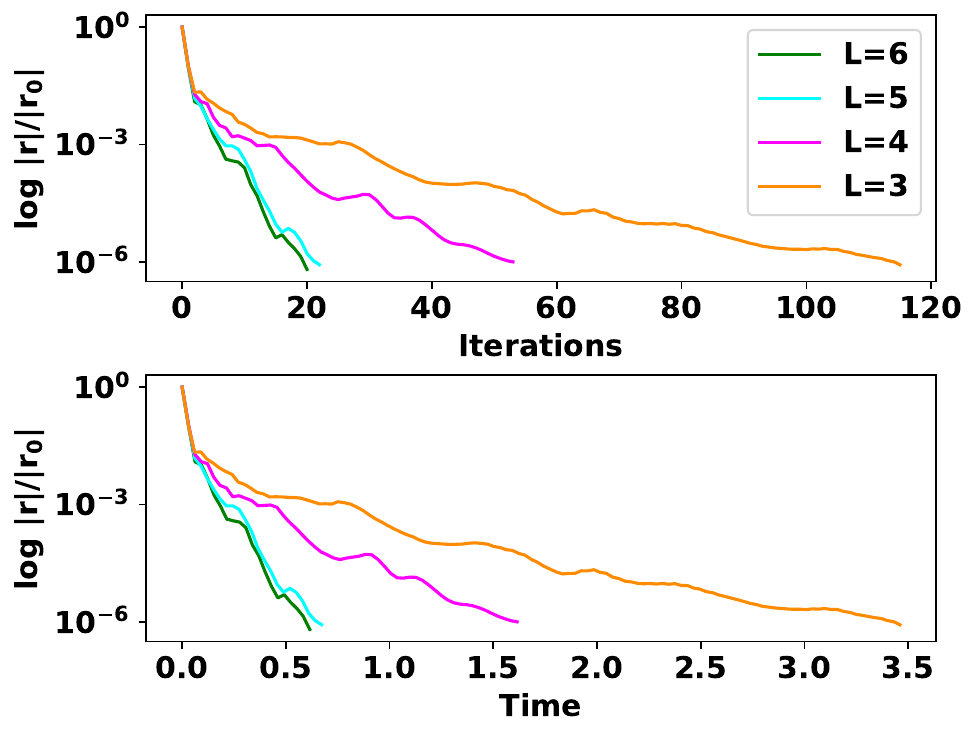}}
    \subfigure[Scooping, $(256,256,256)$]{\label{}\includegraphics[width=0.3\textwidth]{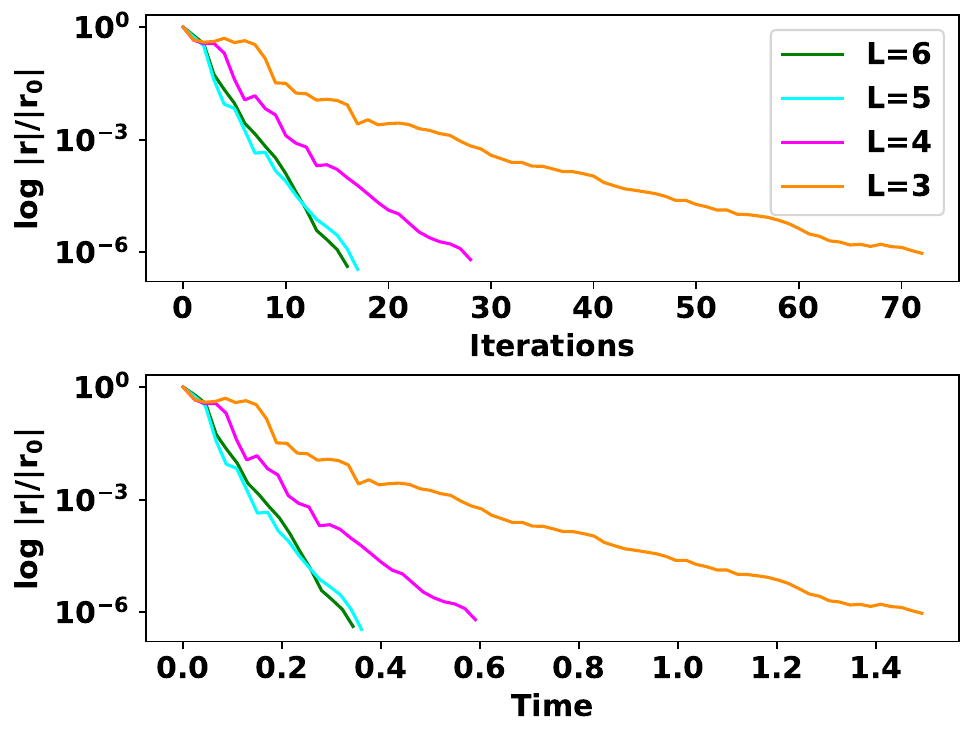}}
    \subfigure[Waterflow torus, $(256,256,256)$]{\label{}\includegraphics[width=0.3\textwidth]{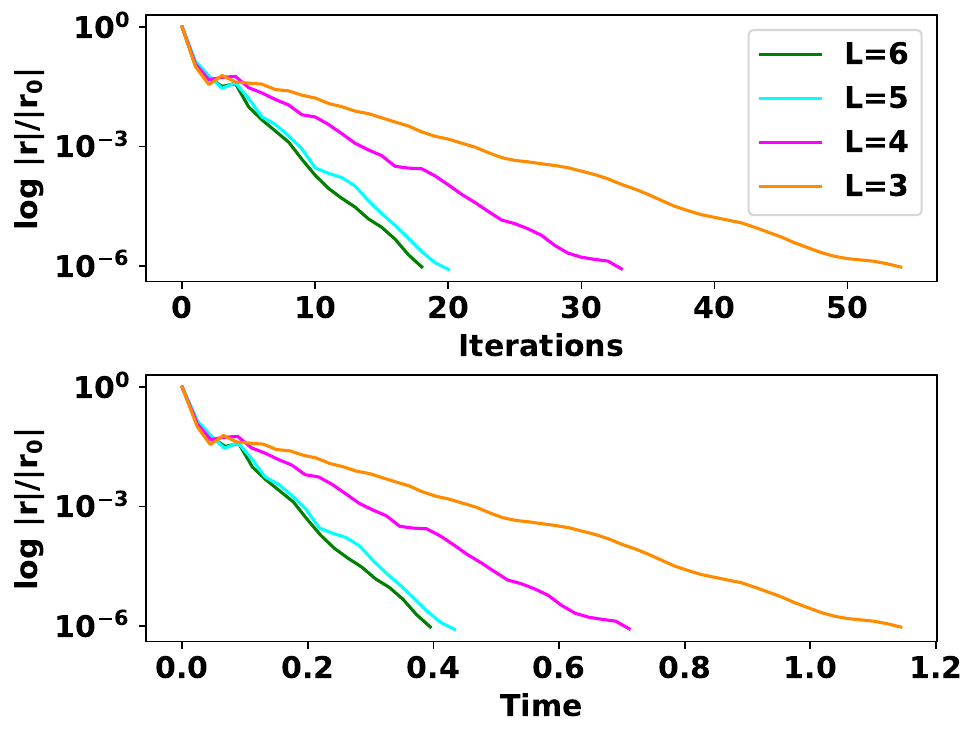}}
    \subfigure[Waterflow ball, $(256,256,256)$]{\label{}\includegraphics[width=0.3\textwidth]{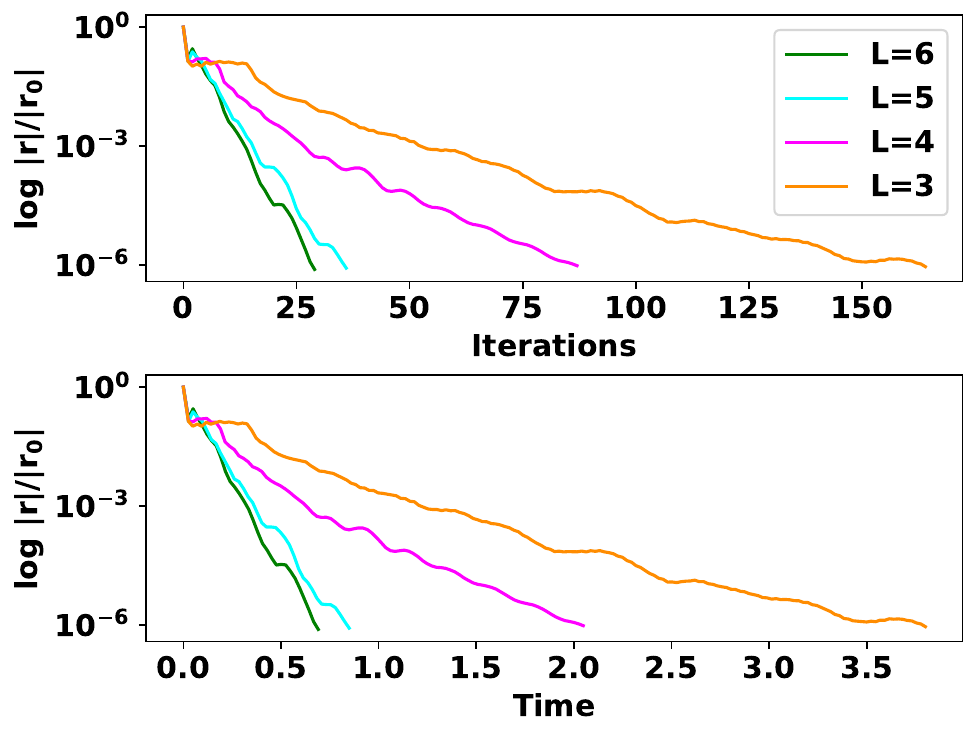}}
  \label{fig:comparing_levels}
  \caption{Comparisons among models with different depth $\levels$ on 12 frames from all scenes at varying resolutions.}
\end{figure}

\newpage
\subsection{PCG Solver Variant Comparisons}
\label{apx:pcgs}

The following plots compare the performances of various iterative solvers
preconditioned by $\network$. Statistics for unpreconditioned CG are also
included for reference. While PCG and Flexible PCG both perform reasonably,
PSDO achieves a modest speedup over them.

\begin{figure}[h!]
    \centering
    \subfigure[Smoke solid, $(128,128,128)$]{\label{}\includegraphics[width=0.3\textwidth]{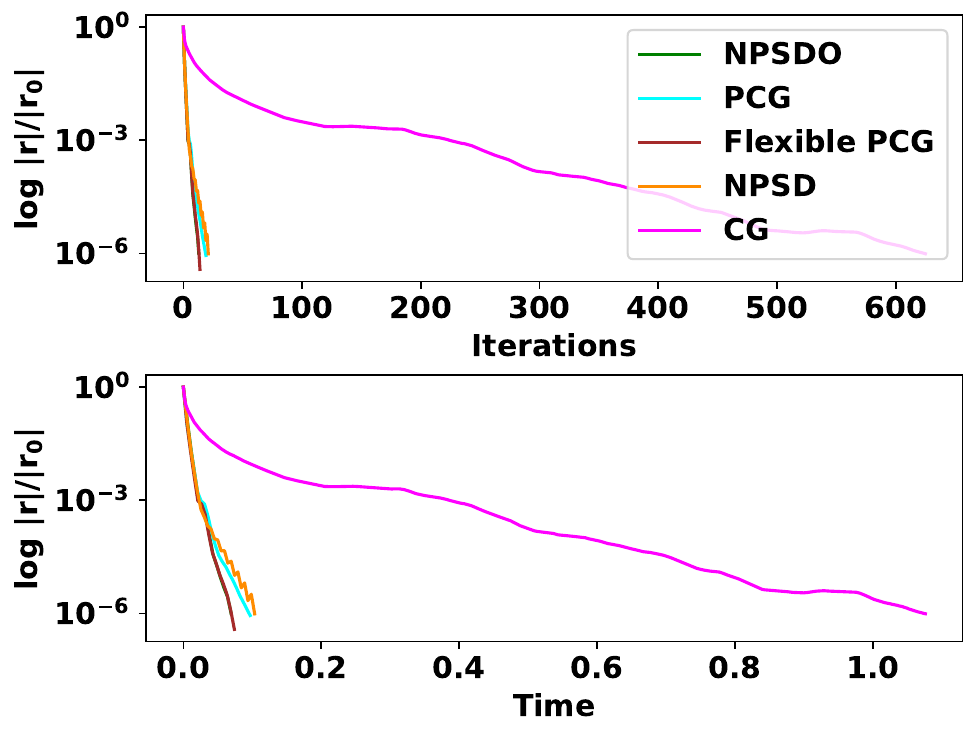}}
    \subfigure[smoke bunny, $(128,128,128)$]{\label{}\includegraphics[width=0.3\textwidth]{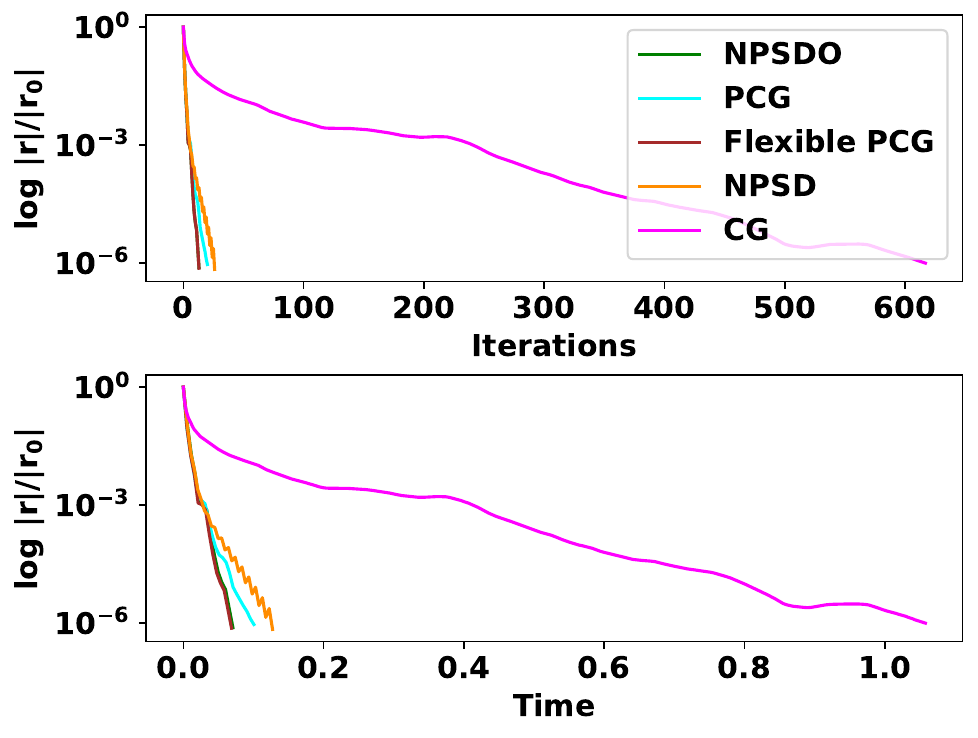}}
    \subfigure[Scooping, $(128,128,128)$]{\label{}\includegraphics[width=0.3\textwidth]{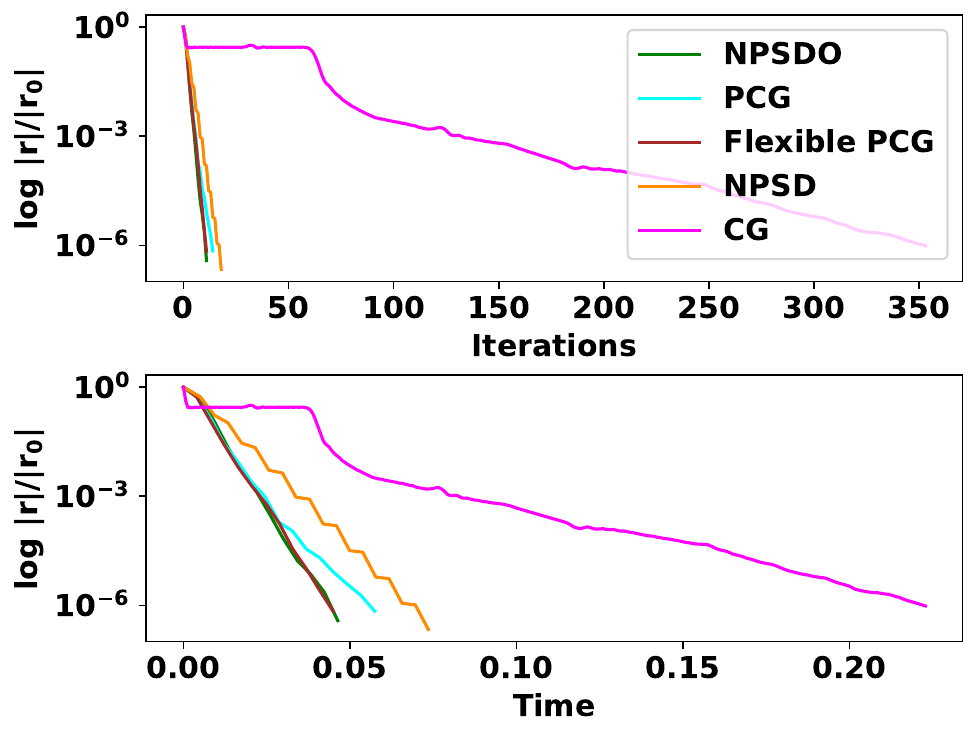}}
    \subfigure[Waterflow torus, $(128,128,128)$]{\label{}\includegraphics[width=0.3\textwidth]{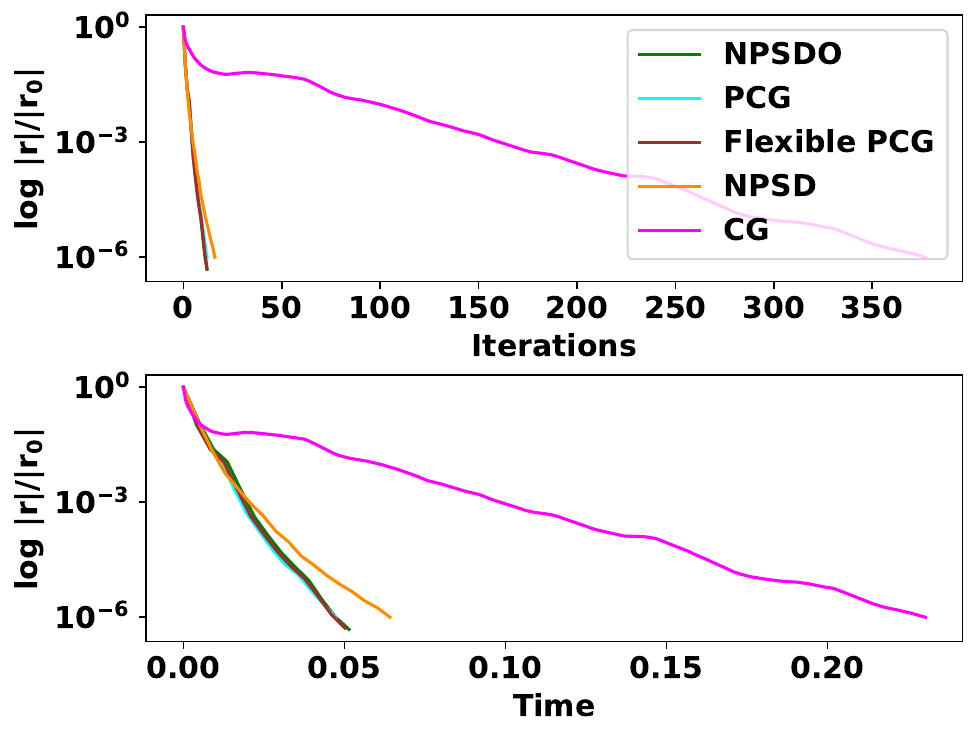}}
    \subfigure[Waterflow ball, $(128,128,128)$]{\label{}\includegraphics[width=0.3\textwidth]{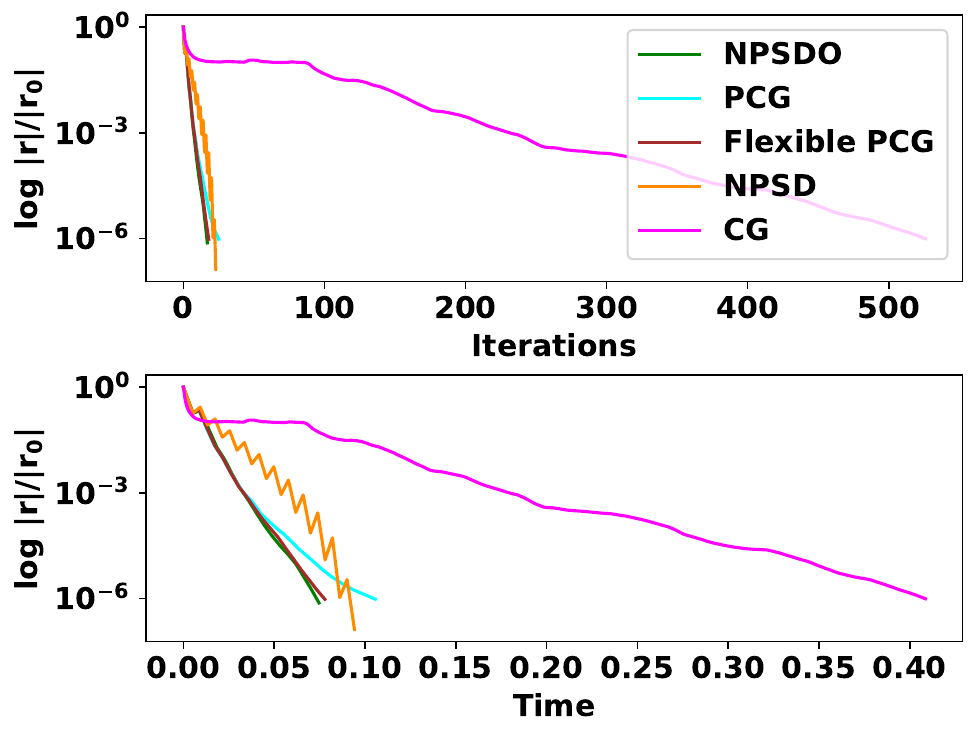}}
    \subfigure[Dambreak pillars, $(256,128,128)$]{\label{}\includegraphics[width=0.3\textwidth]{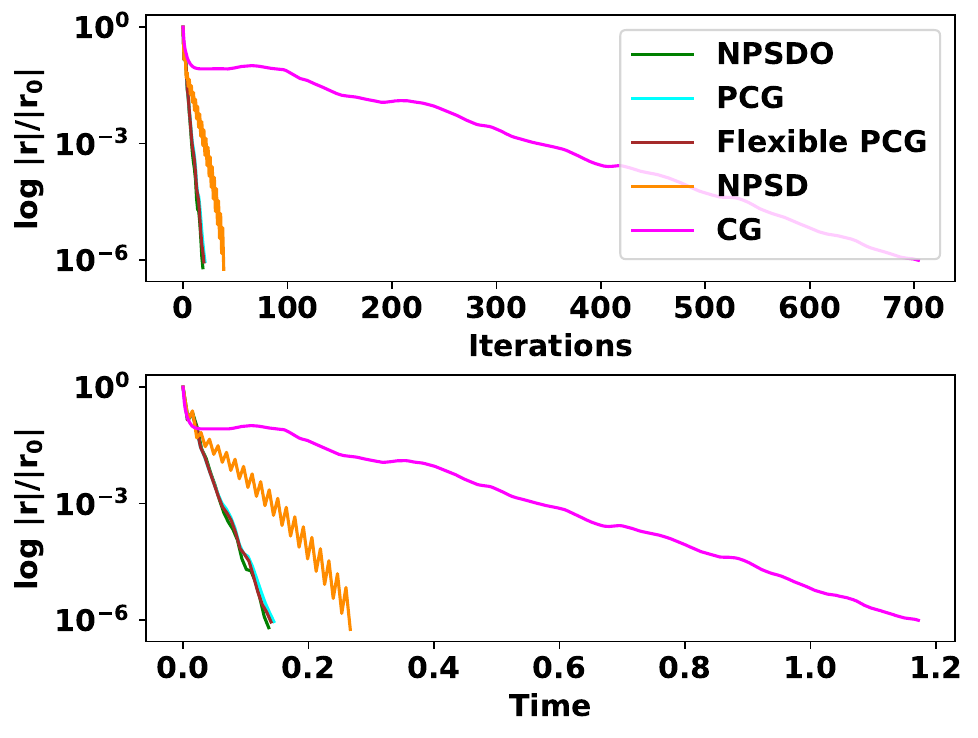}}
    \subfigure[Dambreak bunny, $(256,128,128)$]{\label{}\includegraphics[width=0.3\textwidth]{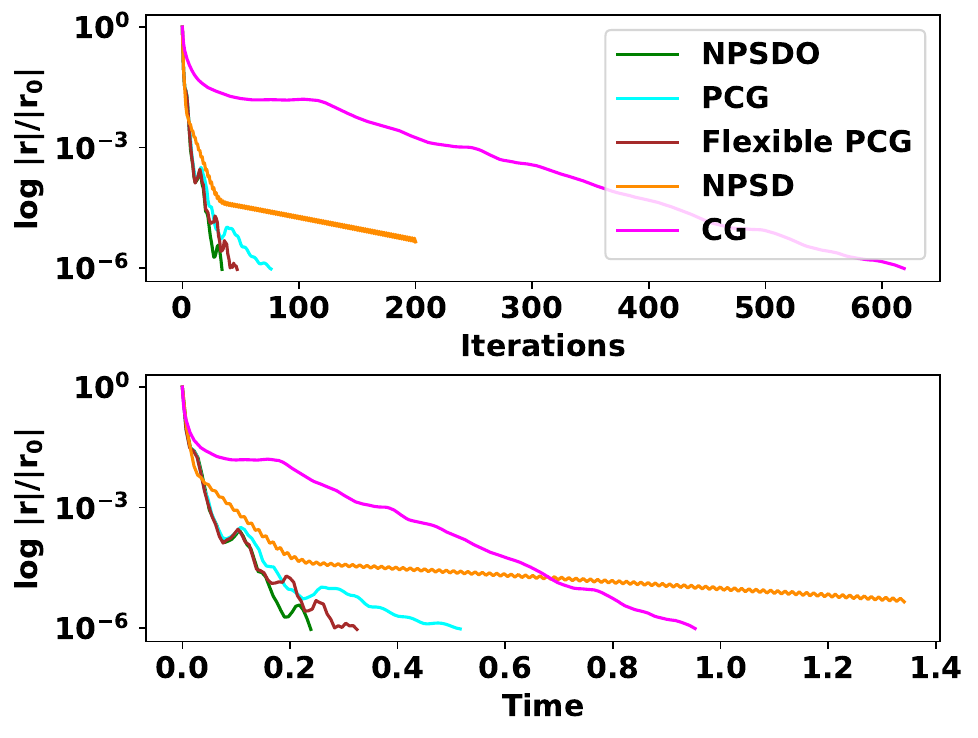}}
    \subfigure[Smoke solid, $(256,256,256)$]{\label{}\includegraphics[width=0.3\textwidth]{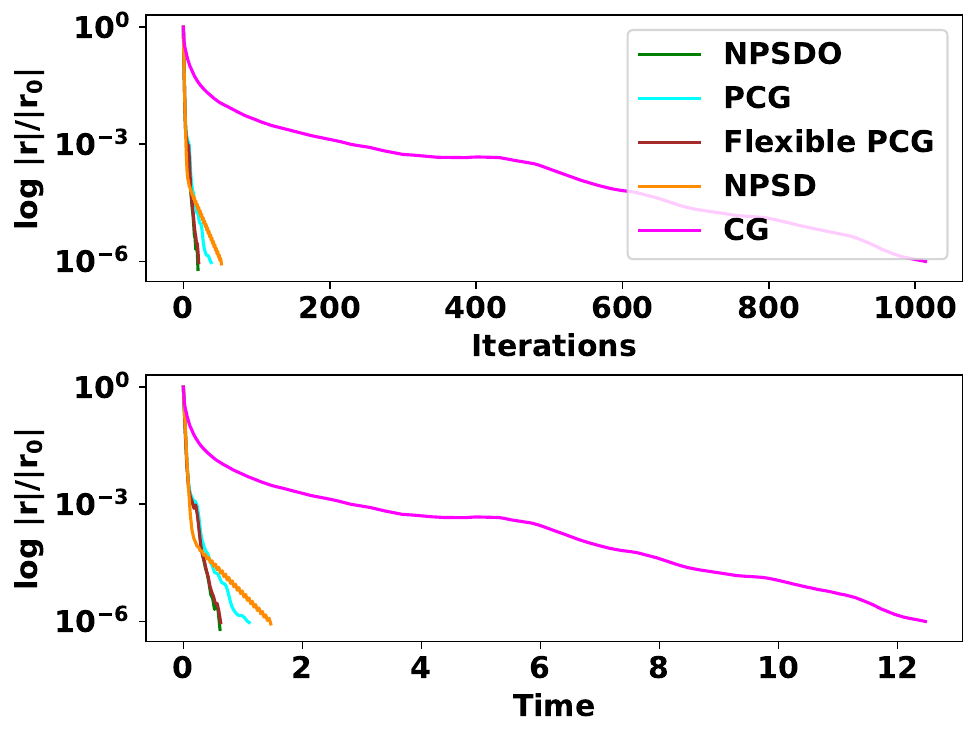}}
    \subfigure[Smoke bunny, $(256,256,256)$]{\label{}\includegraphics[width=0.3\textwidth]{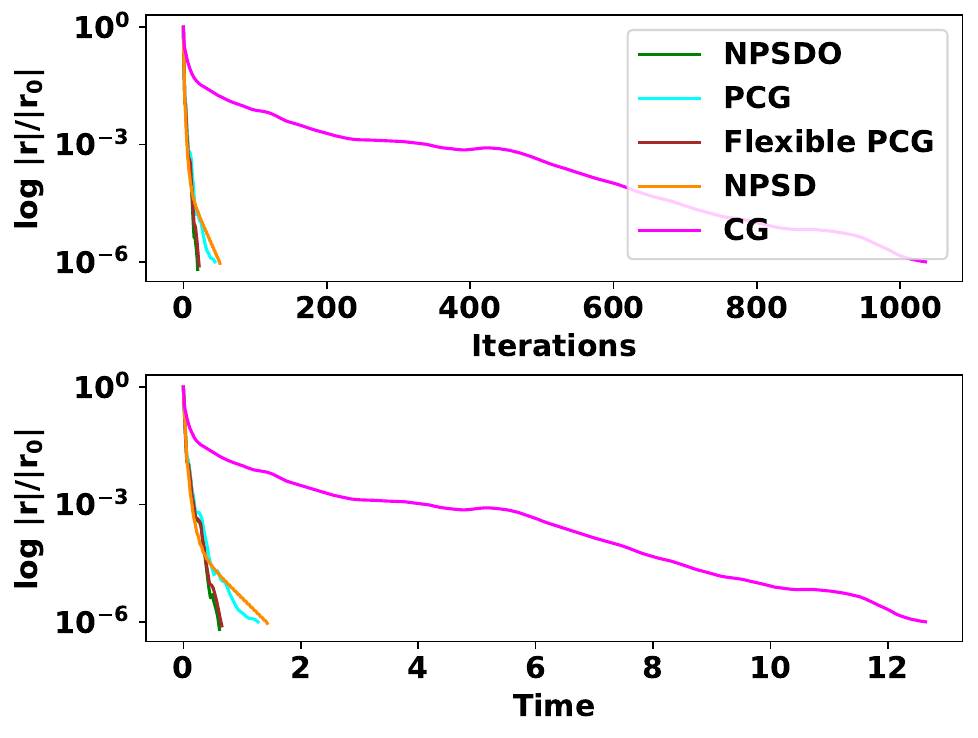}}
    \subfigure[Scooping, $(256,256,256)$]{\label{}\includegraphics[width=0.3\textwidth]{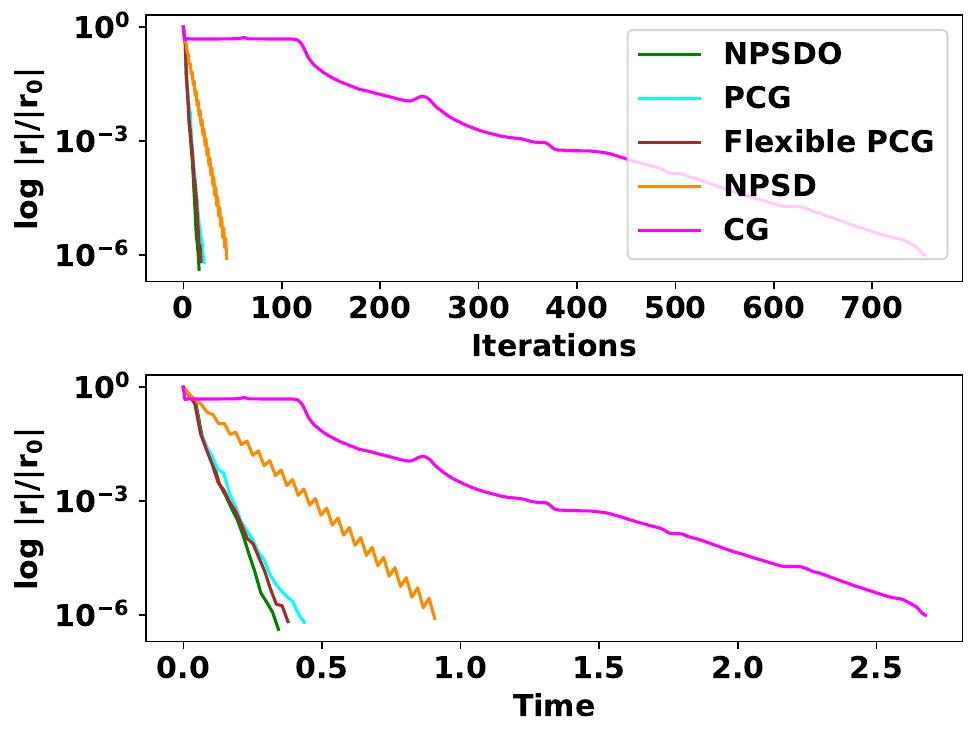}}
    \subfigure[Waterflow torus, $(256,256,256)$]{\label{}\includegraphics[width=0.3\textwidth]{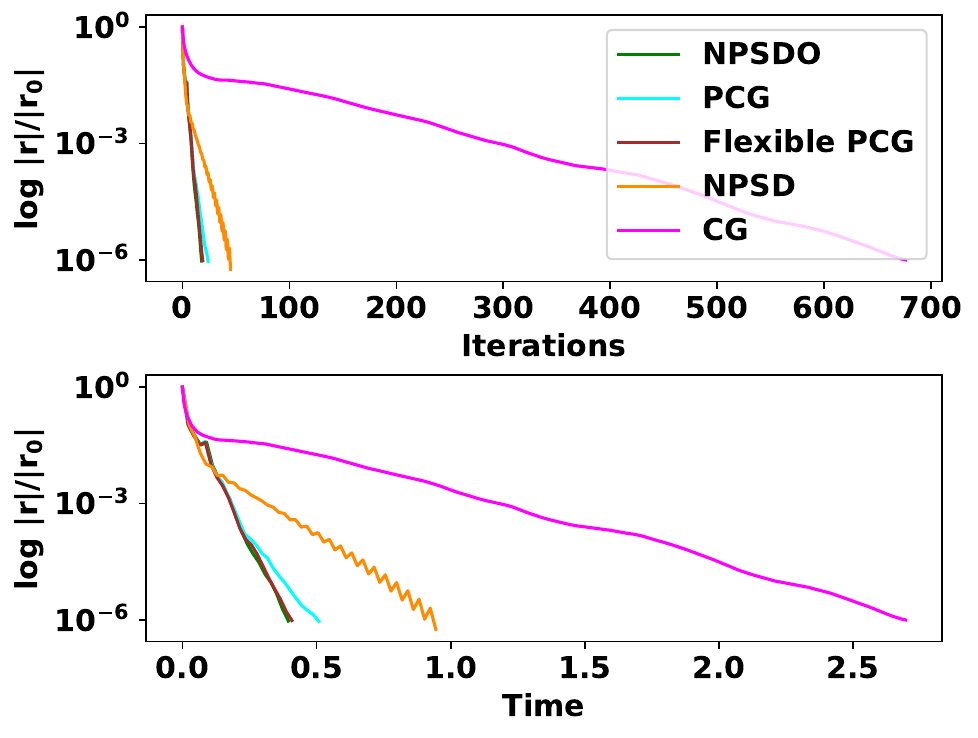}}
    \subfigure[Waterflow ball, $(256,256,256)$]{\label{}\includegraphics[width=0.3\textwidth]{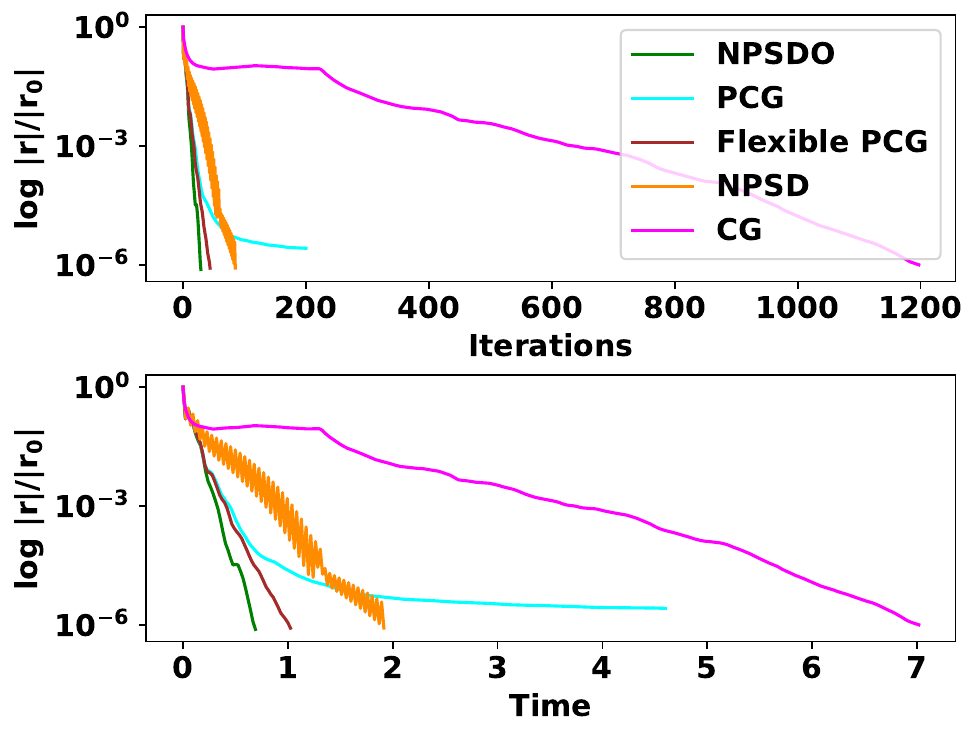}}
  \label{fig:comparing_pcgs}
   \caption{Performance statistics for several solver variants on 12 frames from all scenes}
\end{figure}

\newpage
\subsection{Number of Orthogonalizations}\label{apx:ortho}
The following plots compare the performances of our NPSDO solver with different numbers of orthogonalization steps ($\numOrtho$ in \pr{alg:npsdo}).

\begin{figure}[h!]
    \centering
    \subfigure[Smoke solid, $(128,128,128)$]{\label{}\includegraphics[width=0.3\textwidth]{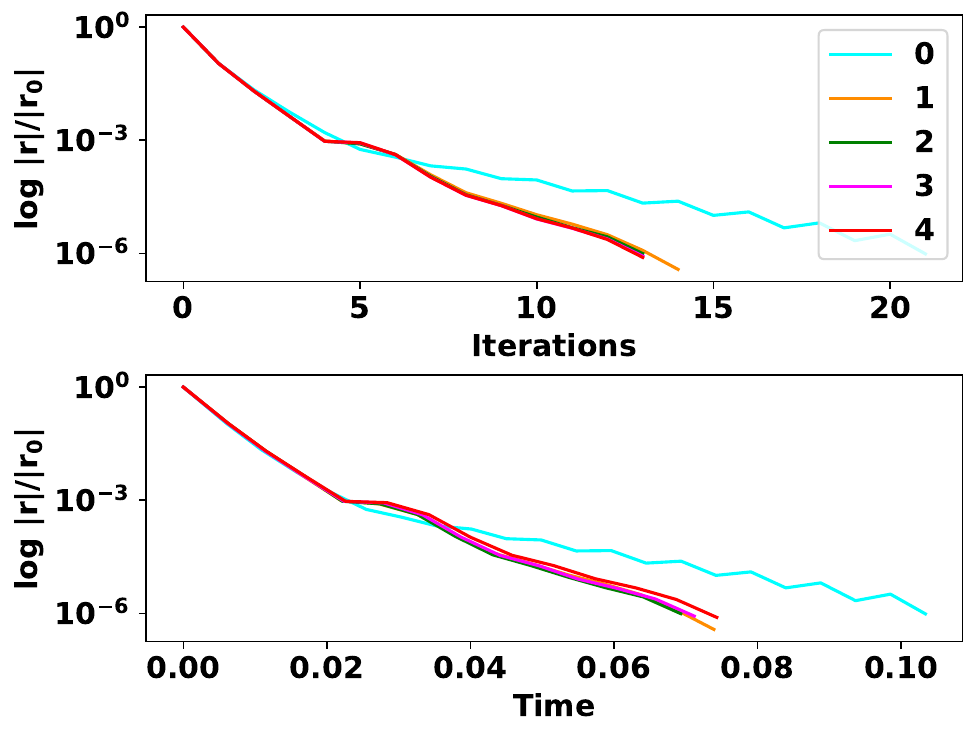}}
    \subfigure[smoke bunny, $(128,128,128)$]{\label{}\includegraphics[width=0.3\textwidth]{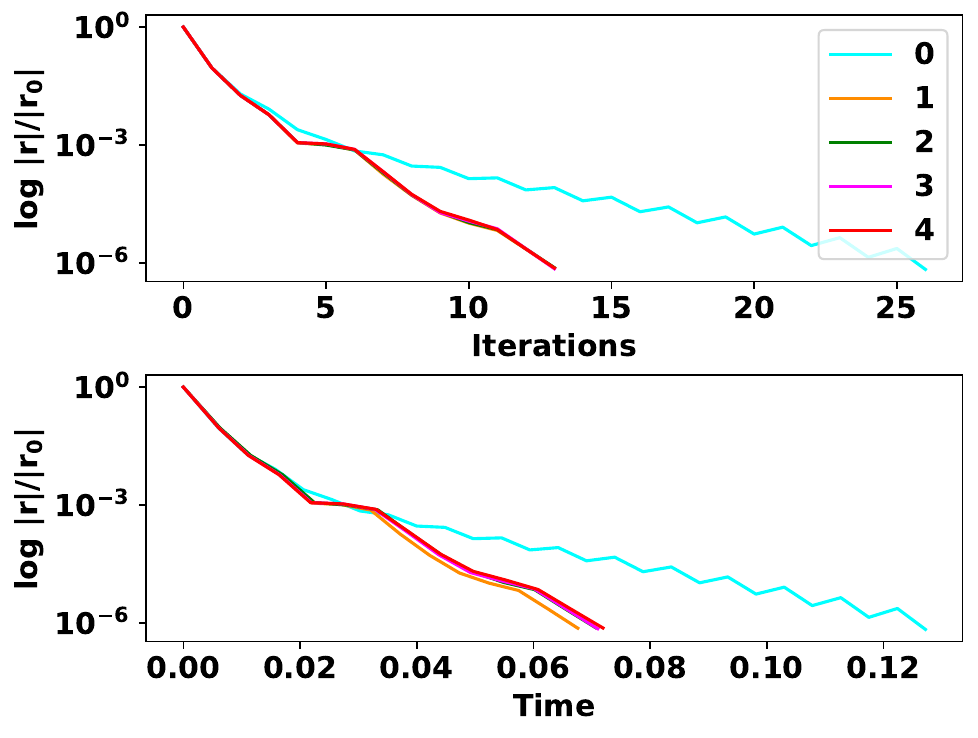}}
    \subfigure[Scooping, $(128,128,128)$]{\label{}\includegraphics[width=0.3\textwidth]{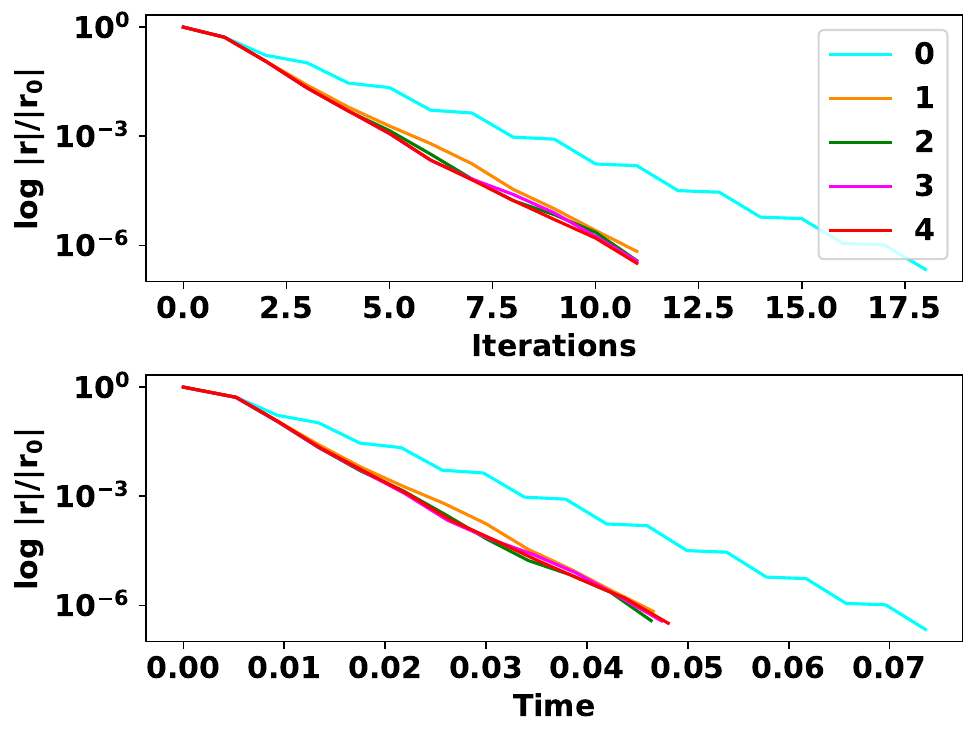}}
    \subfigure[Waterflow torus, $(128,128,128)$]{\label{}\includegraphics[width=0.3\textwidth]{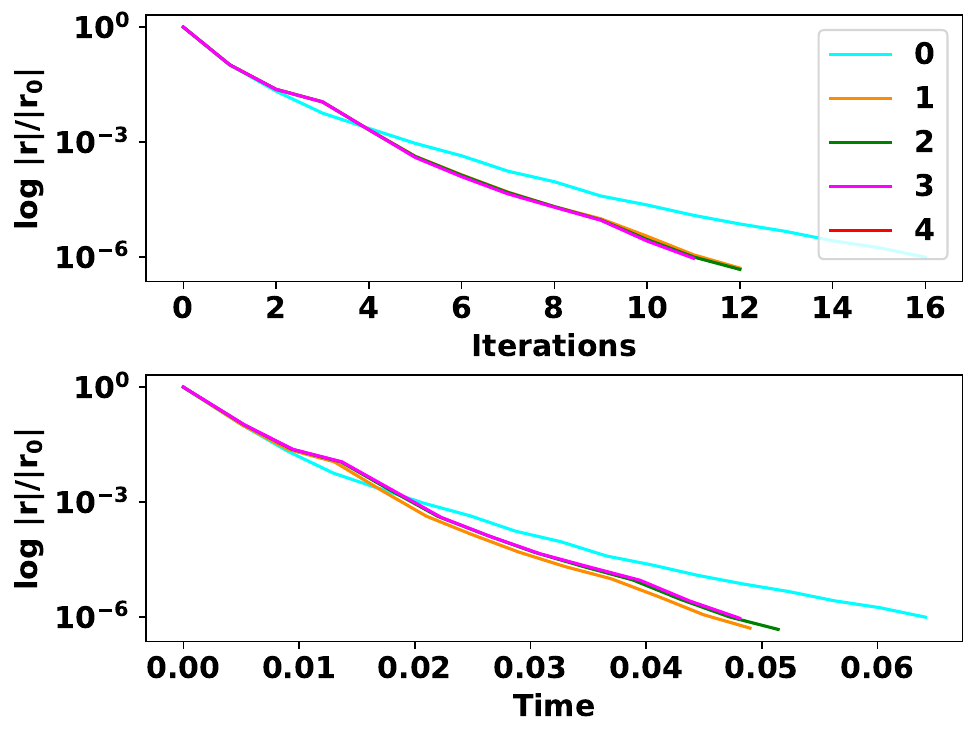}}
    \subfigure[Waterflow ball, $(128,128,128)$]{\label{}\includegraphics[width=0.3\textwidth]{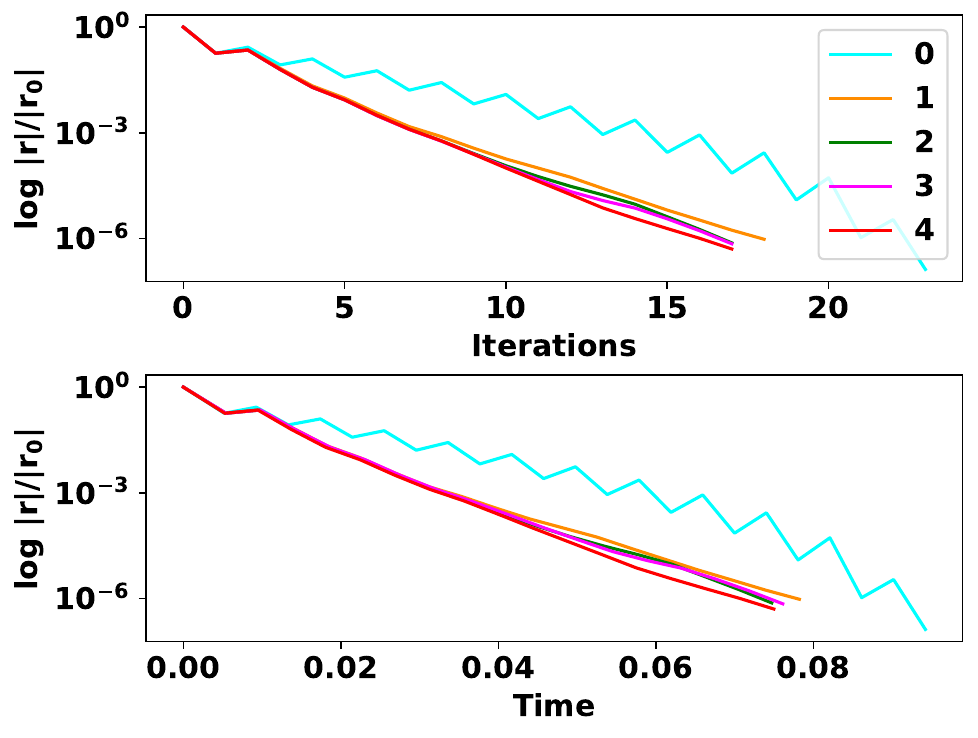}}
    \subfigure[Dambreak pillars, $(256,128,128)$]{\label{}\includegraphics[width=0.3\textwidth]{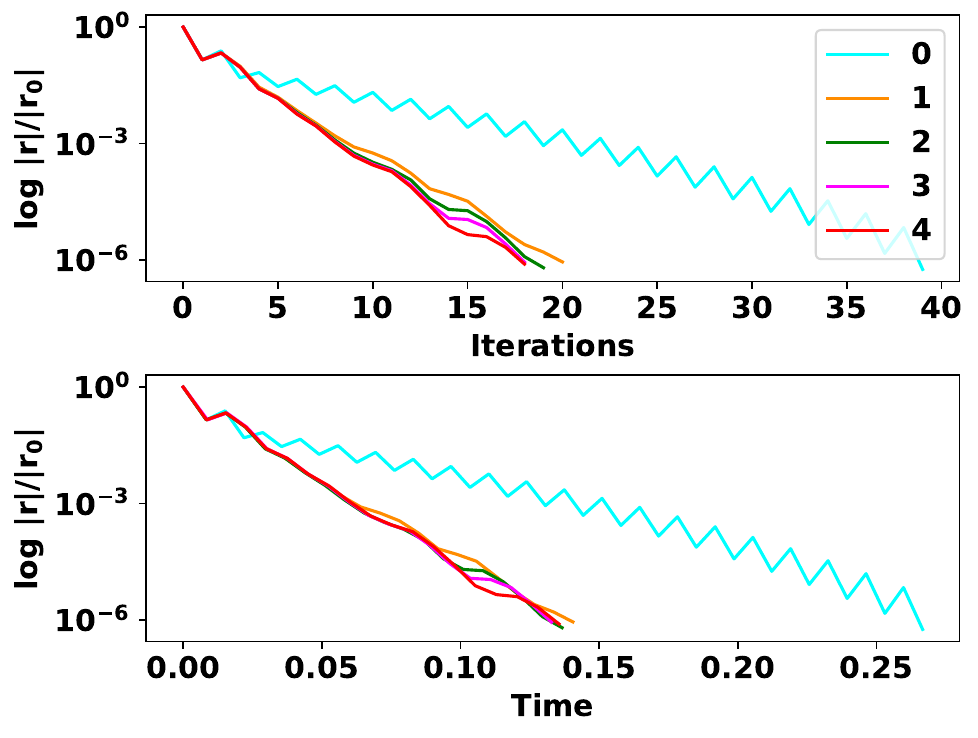}}
    \subfigure[Dambreak bunny, $(256,128,128)$]{\label{}\includegraphics[width=0.3\textwidth]{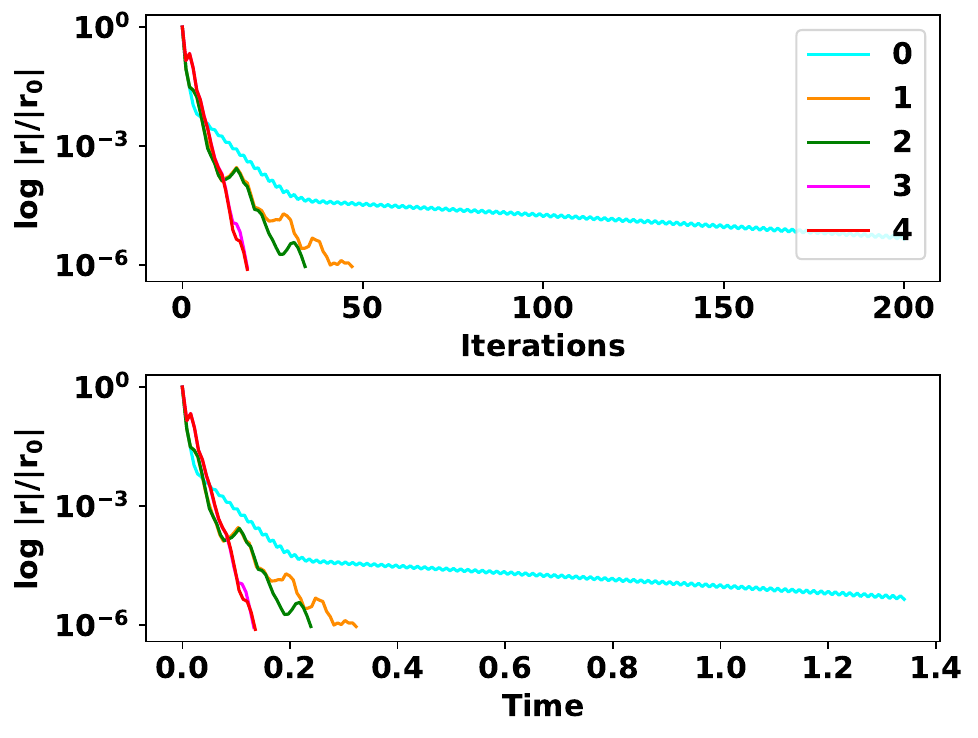}}
    \subfigure[Smoke solid, $(256,256,256)$]{\label{}\includegraphics[width=0.3\textwidth]{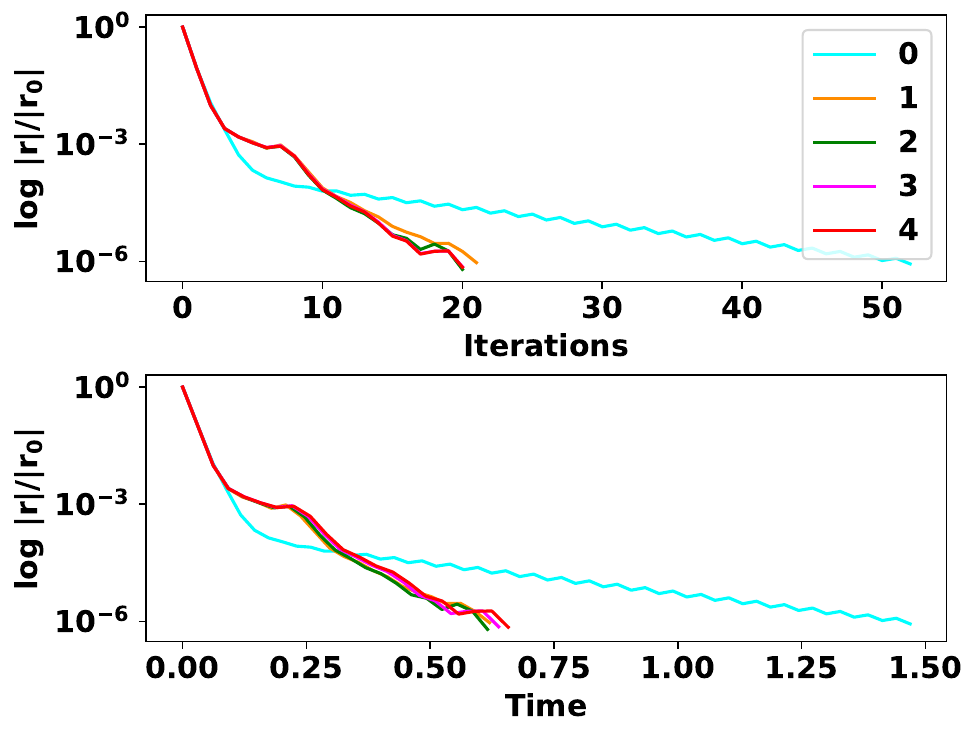}}
    \subfigure[Smoke bunny, $(256,256,256)$]{\label{}\includegraphics[width=0.3\textwidth]{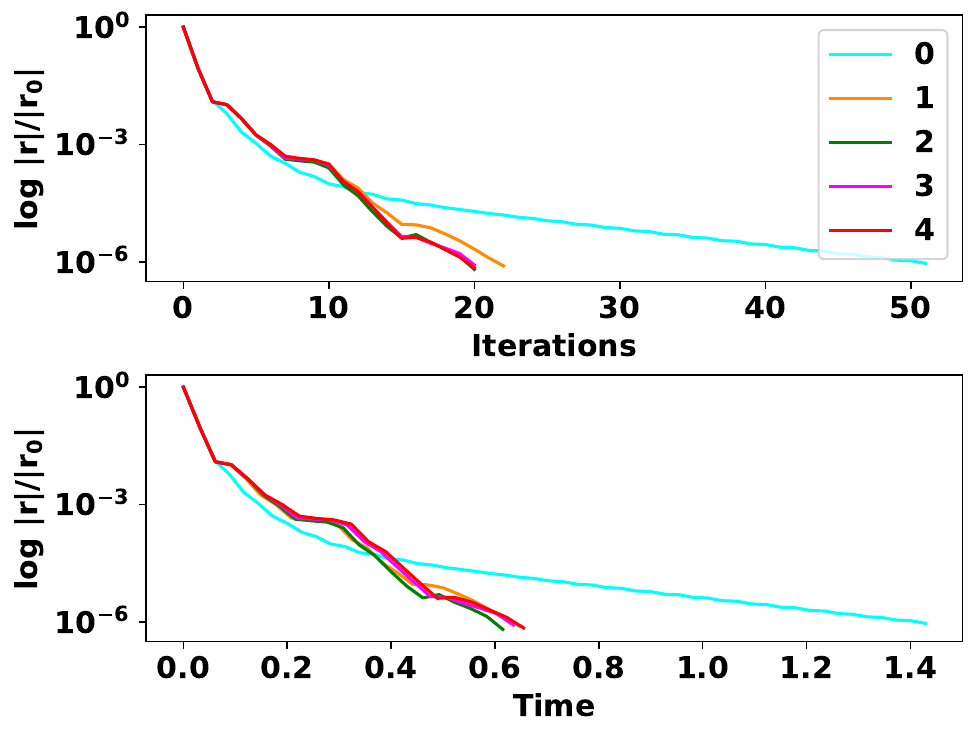}}
    \subfigure[Scooping, $(256,256,256)$]{\label{}\includegraphics[width=0.3\textwidth]{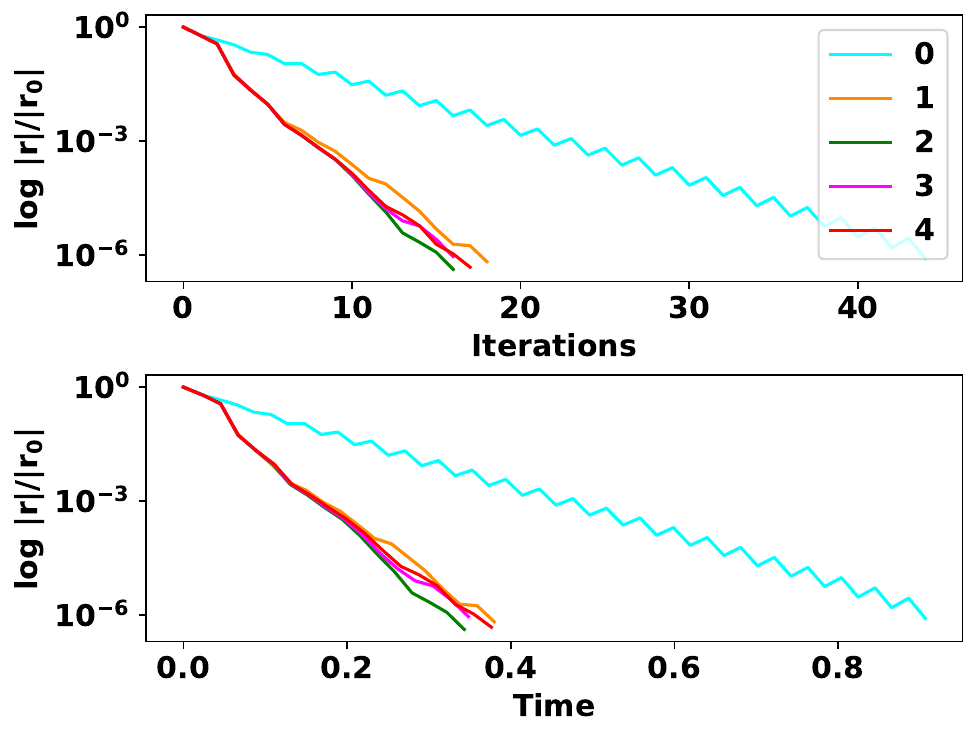}}
    \subfigure[Waterflow torus, $(256,256,256)$]{\label{}\includegraphics[width=0.3\textwidth]{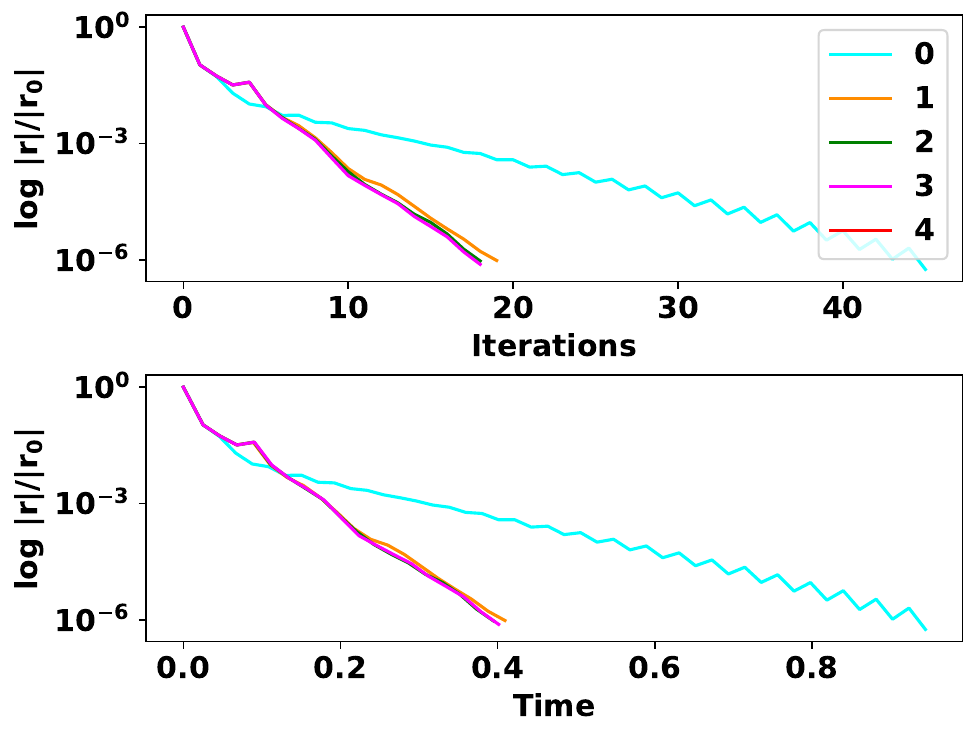}}
    \subfigure[Waterflow ball, $(256,256,256)$]{\label{}\includegraphics[width=0.3\textwidth]{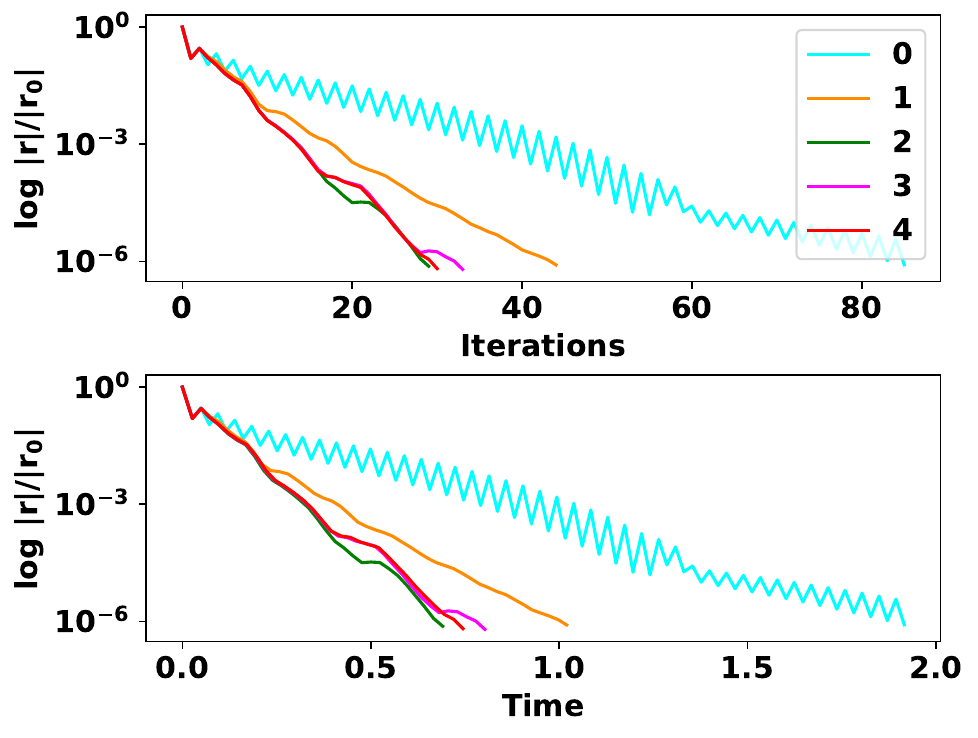}}
  \label{fig:comparing_ortho}
   \caption{Performance statistics for different $\numOrtho$ on 12 frames from all scenes}
\end{figure}

\newpage
\subsection{Timing Benchmarks}\label{apx:benchmarking}
   This section lists the average iteration count and runtime for each
   test case.
   \vskip 0.15in
   \begin{table}[h!]
   \caption{Average iteration count and runtime across all frames for each test suite} \label{tbl:benchmark}
   \begin{center}
      \scriptsize
   \begin{tabular}{l |*{5}{*{2}{c}}}
         & \multicolumn{2}{c}{\textbf{AMGCL}} &  \multicolumn{2}{c}{\textbf{AMGX}}& \multicolumn{2}{c}{\textbf{IC}} & \multicolumn{2}{c}{\textbf{CG}} & \multicolumn{2}{c}{\textbf{Ours}} \\
   Examples & \textbf{Iteration} & \textbf{Time} & \textbf{Iteration} & \textbf{Time} & \textbf{Iteration} & \textbf{Time} & \textbf{Iteration} & \textbf{Time} & \textbf{Iteration} & \textbf{Time} \\
   \hline \\

   Smoke solid $128^3$ & 11.2 &	0.696 & 44.9 & 0.324 & 196.1 &	0.852	& 612.5	& 1.050	& 13.3 & 0.072\\\\
   Smoke solid $256^3$ & 14.3 &	6.433 & 64.8 & 3.156 & 343.7 &	11.1  & 1076.8	& 13.245	& 20.5 & 0.636       \\\\
   Smoke bunny $128^3$ & 11.3 &  0.713 & 43.7 & 0.319 & 200.5 & 0.862   & 615.8  & 1.048  & 14.0 & 0.075        \\\\
   Smoke bunny $256^3$ & 14.3 &	6.046 & 63.7 & 3.119 & 345.7 & 10.9    & 1069.7 & 13.042	& 21.0 & 0.648      \\\\
   Scooping $128^3$ &    12.2 &  0.170 & 34.6 & 0.097 & 128.1 & 0.270   & 392.8  & 0.234  & 11.8&	0.051     \\\\
   Scooping $256^3$ &    12.2 &	1.643 & 50.5&0.738 &246.0 & 1.968   & 749.2	& 2.593  & 18.1 &0.388      \\\\
   Waterflow torus $128^3$ & 9.4  & 0.086 & 20.6&0.046 & 66.9 & 0.129 & 208.6  & 0.095  & 12.0&	0.048 \\\\
   Waterflow torus $256^3$ & 10.4 & 1.024 &30.1&0.346 & 130.0 & 0.946  & 401.2  & 1.121  & 16.6 & 0.338 \\\\
   Waterflow ball $128^3$ & 11.2  & 0.173 & 32.7 & 0.098 & 136.5 & 0.275  & 405.4  & 0.257  & 16.1 & 0.068      \\\\
   Waterflow ball $256^3$ & 12.1  & 2.442 & 57.6 & 1.167  & 328.4 & 3.559  & 969.7  & 4.875  & 28.5 & 0.650   \\\\
   Dambreak pillars $256\cdot 128^2$ & 11.2 & 0.476 & 39.9 & 0.229 & 167.7 & 0.625 & 523.2 & 0.704 & 14.9 & 0.103   \\\\
   Dambreak bunny $256\cdot 128^2$ & 11.1   & 0.395 & 36.7 &	0.183 & 148.2 & 0.514 & 452.1 & 0.522 & 19.3 & 0.129     \\
   \hline \\
   Average & 11.7 & 1.691 & 43.3& 0.819 & 203.2 &	2.667 & 623.1 & 3.232 &	17.2 & 0.267
   \end{tabular}
   \end{center}
   \end{table}

\subsection{Memory Usage}\label{apx:memory}
   The following peak memory usage (in MiB) statistics were recorded with the command \verb|nvidia-smi| on GeForce RTX 3080. Examples that threw out-of-memory error are marked as \texttt{NA}.

   \begin{table}[h!]
      \caption{Peak memory usage on a smoke simulation.}
      \centering
      \vskip 0.1in
      \begin{tabular}{l |*{7}{c}}
      \textbf{Resolution}   & \textbf{AMGCL}  & \textbf{AMGX} & \textbf{IC} & \textbf{CG} & \textbf{NPSDO} & \textbf{DCDM} & \textbf{FN} \\
      \hline\\
      $128^3$ & 1248  & 2150 & 1668  & 1418  & 1548  & 8532  & 5170  \\\\
      $256^3$ & 5032  & 7370 & 8214  & 3716  & 4776  & NA & NA

      \end{tabular}

   \end{table}

\newpage
\subsection{New Grid Sizes}\label{apx:additional_domain_size}

The following plots demonstrate the ability of our solver to generalize to
domains of different resolutions \emph{without retraining the
preconditioner}: it maintains consistently strong
convergence behavior when applied to simulations with grid dimensions
not seen during training. We note that $256^3$ resolution simulations
reported on in the main paper also were not present in the training set.

\begin{figure}[h!]
   \centering
    \includegraphics[width=0.35\textwidth]{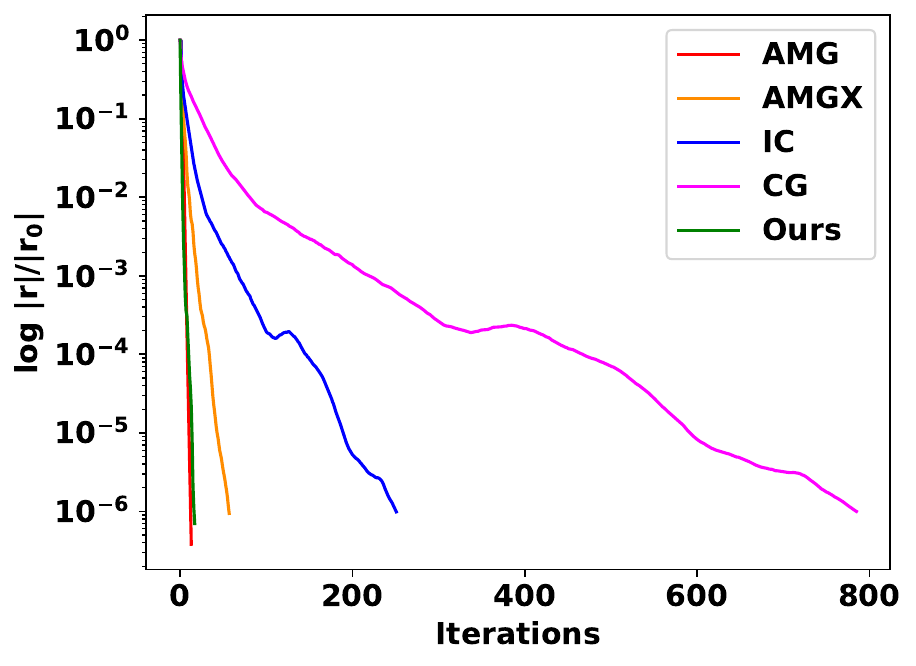}
    \includegraphics[width=0.35\textwidth]{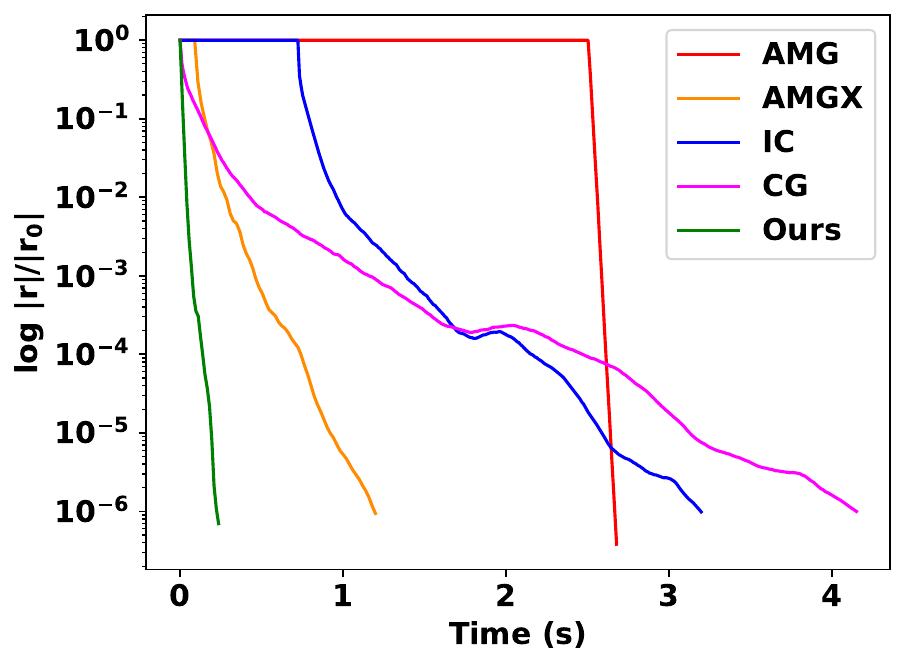}
    \includegraphics[width=0.35\textwidth]{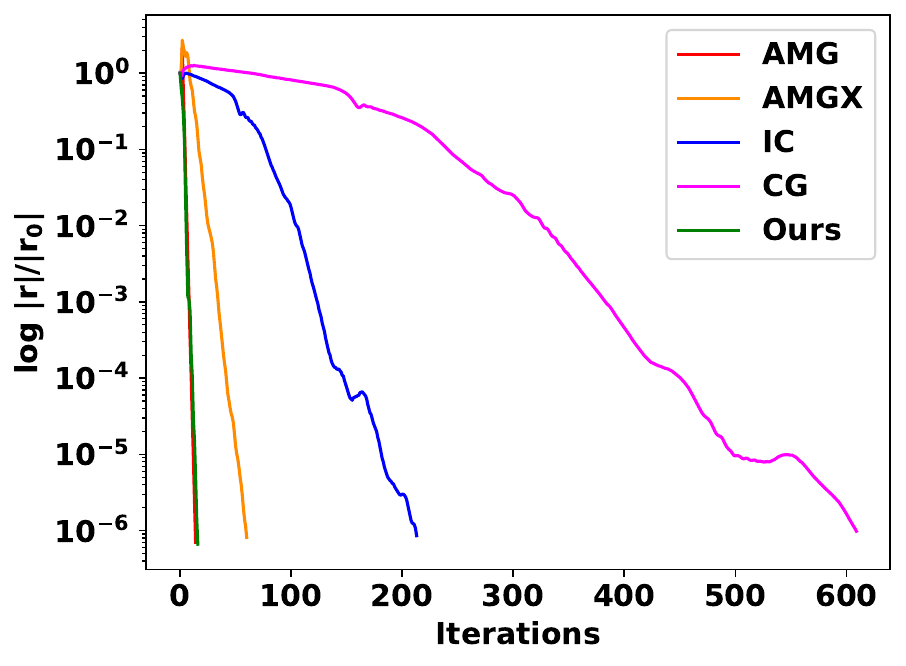}
    \includegraphics[width=0.35\textwidth]{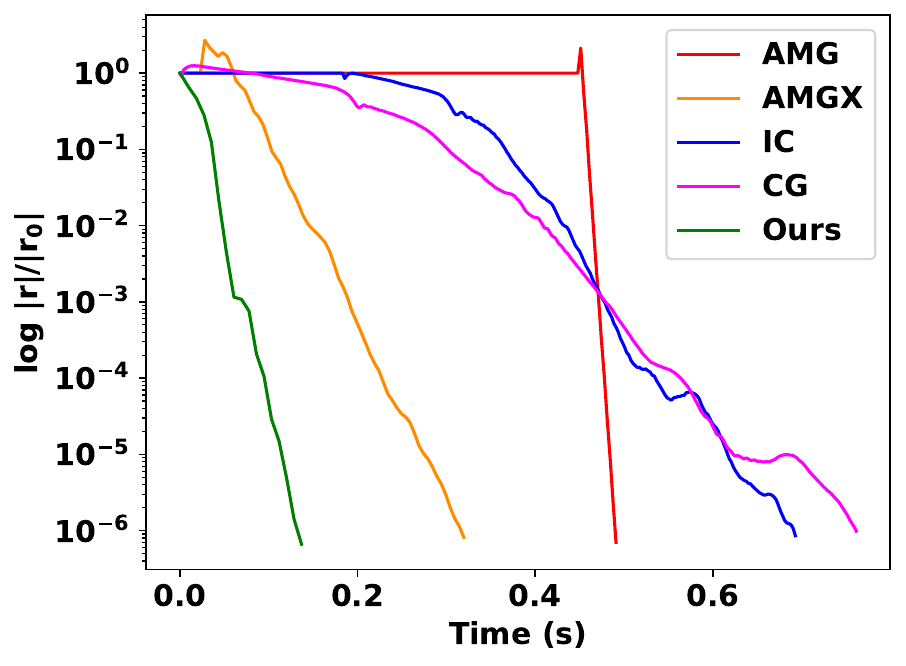}
    \caption{Convergence on a representative frame of a $192^3$-resolution Neumann-only simulation (top pair) and a $(384, 128, 128)$-resolution mixed BC simulation (bottom pair).}
\end{figure}

\newpage
\subsection{Analysis of and Enforcement of Symmetry} \label{apx:symmetry_analysis}
To analyze the deviation of our preconditioner from symmetry,
we applied our network to pairs of test vectors $\mathbf{x}$ and $\mathbf{y}$ and evaluated the following
relative asymmetry measure:
$$
\mathcal{S}(\x, \y; \image) \defeq \frac{\left| \x\cdot \network(\image, \y) - \y\cdot \network(\image, \x) \right|}{\sqrt{\left|\x\cdot \network(\image, \x)\right| \left|\y\cdot \network(\image, \y) \right|}}.
$$
We did this using 100 test vector pairs for each system (image) in our training set,
reporting the aggregate statistics in \pr{fig:symmetry} and \pr{fig:symmetry_rand} for two different types of test vectors:
(a) random pairs of vectors from the training set, and (b) fully random vectors from a normal distribution.
Not only is the trained network's asymmetry small, it is several orders of magnitude smaller than that of the initial network,
suggesting that training the network to approximate the inverse of the discrete Laplace operator (a symmetric matrix)
naturally promotes symmetry.

\begin{figure}[h!]
   \centering
   \subfigure[Random pairs of training RHS vectors.]{\label{fig:symmetry}\includegraphics[width=0.35\textwidth]{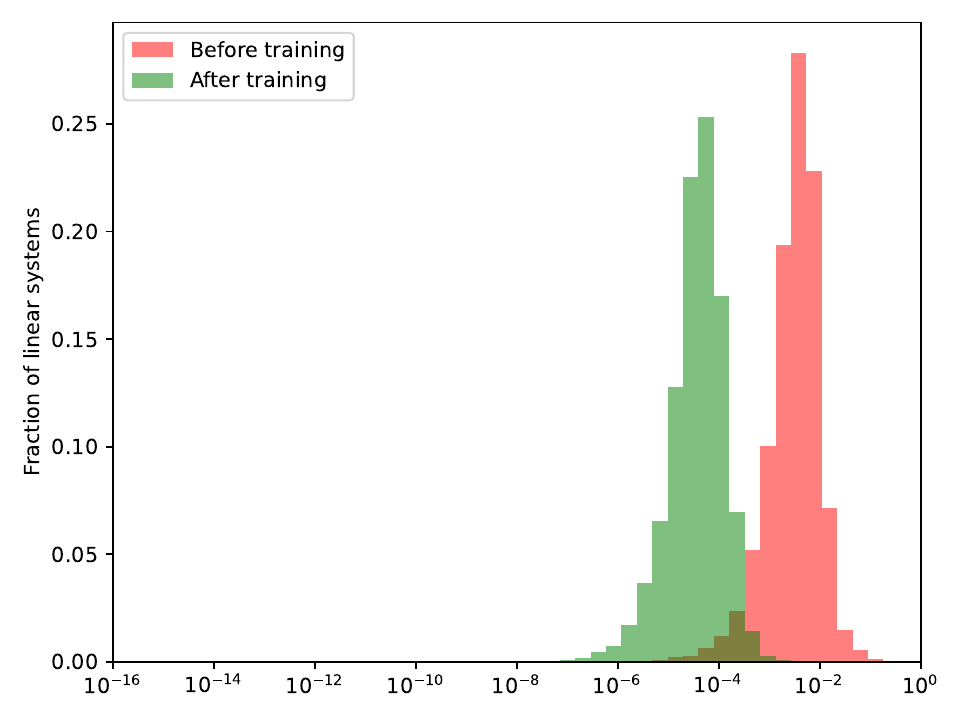}}
    \subfigure[Random vectors from normal distribution.]{\label{fig:symmetry_rand}\includegraphics[width=0.35\textwidth]{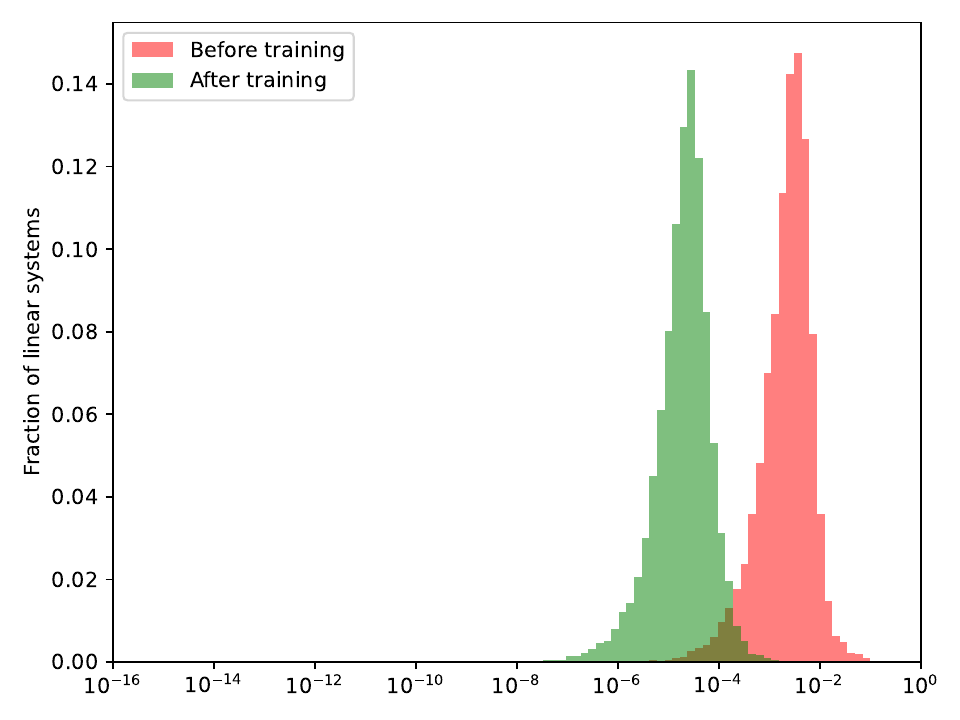}}
    \caption{Histograms of relative symmetry metric $\mathcal{S}(\x, \y; \image)$ across 100 test vector pairs $(\x, \y)$ for each image $\image$ in the training set.}

\end{figure}

We furthermore constructed two variants of our network architecture that enforce symmetry and positive definiteness by construction. The first concatenates the $\network$ architecture with
its transpose ${\network}^\top$, obtaining $\network_\mathrm{sym,1}(\image, \r) \defeq {\network}^\top(\image, \network(\image, \r))$, while
the second constrains the ``pre-smoothing`` and ``post-smoothing`` blocks to be transposes of each other via weight sharing, obtaining $\network_\mathrm{sym,2}$.
However, even when training on a single frame, both of these variants exhibit suboptimal performance:
the training loss is not sufficiently reduced compared to our proposed model, and both symmetric networks perform worse in terms of wallclock time than unpreconditioned CG.

\begin{figure}[h!]
   \centering
   \subfigure[Training loss]{\label{}\includegraphics*[width=0.3\textwidth]{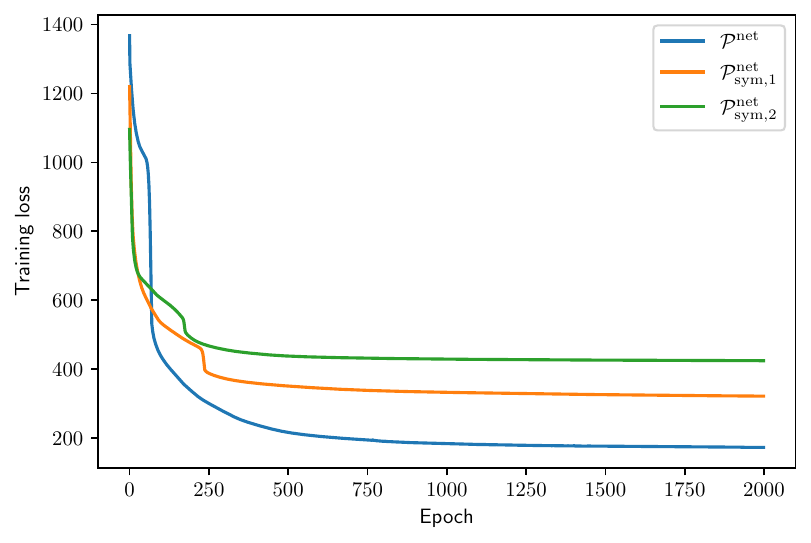}}
   \subfigure[Residual vs. iterations]{\label{fig:symmetry}\includegraphics[width=0.3\textwidth]{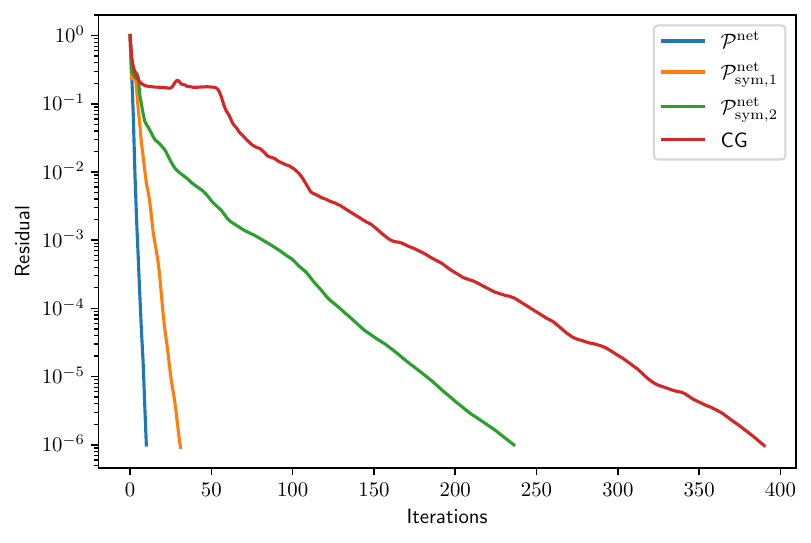}}
    \subfigure[Residual vs. time]{\label{fig:symmetry_rand}\includegraphics[width=0.3\textwidth]{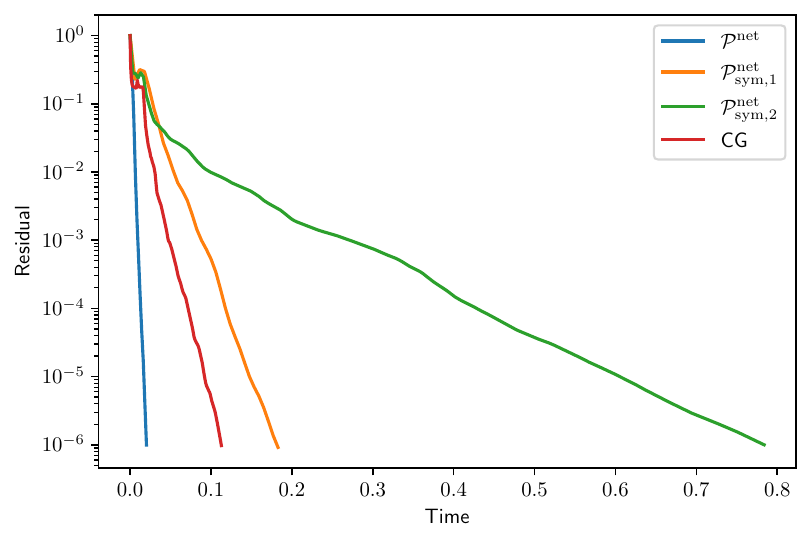}}
    \caption{Comparison between our proposed model and two symmetric models, with baseline method unpreconditioned CG.}
\end{figure}

\newpage
\subsection{Capacity for Residual Reduction}\label{apx:converge_test}

To confirm that our solver can achieve high accuracy despite the single-precision
arithmetic
used in evaluating $\network$, we disabled the convergence test and
ran for a fixed number of iterations (100), recording the final relative
residual achieved for each system in our test set in \pr{fig:final_relative_residual}.
The median relative residual is $4.01\times 10^{-15}$.

\begin{figure}[h!]
   \centering
   \includegraphics[width=0.55\textwidth]{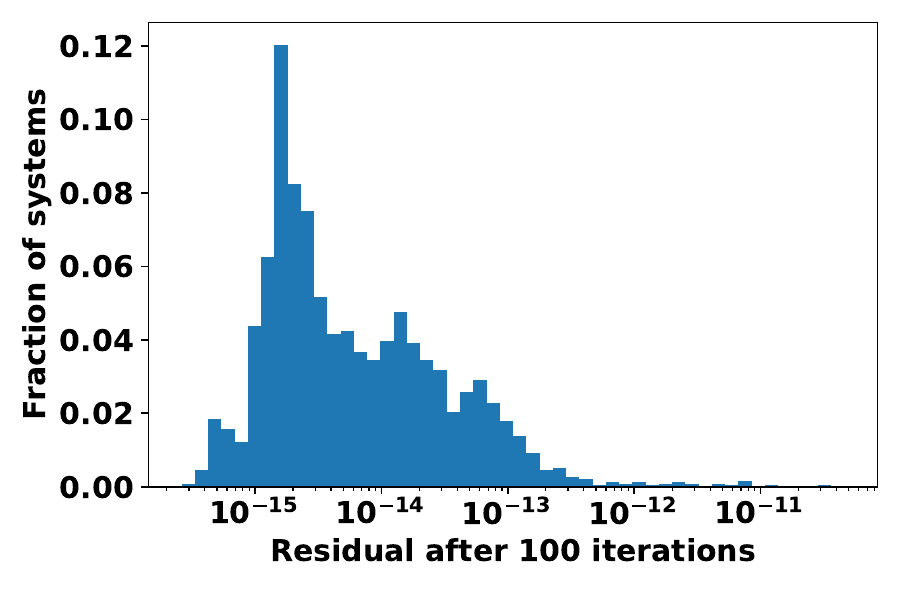}
   \caption{Histogram of the final relative residual achieved after running 100 iterations without a convergence test.}
   \label{fig:final_relative_residual}

\end{figure}

\end{document}

%% file: math_commands.tex
\usepackage{amsmath,amsfonts,bm}

\usepackage{prettyref}
\newrefformat{sec}{Section~\ref{#1}}
\newrefformat{tbl}{Table~\ref{#1}}
\newrefformat{fig}{Figure~\ref{#1}}
\newrefformat{chp}{Chapter~\ref{#1}}
\newrefformat{eqn}{\eqref{#1}}
\newrefformat{set}{\eqref{#1}}
\newrefformat{alg}{Algorithm~\ref{#1}}
\newrefformat{apx}{Appendix~\ref{#1}}
\newrefformat{prop}{Proposition~\ref{#1}}
\newcommand\pr[1]{\prettyref{#1}}

\def\1{\bm{1}}

\DeclareMathAlphabet{\mathsfit}{\encodingdefault}{\sfdefault}{m}{sl}
\SetMathAlphabet{\mathsfit}{bold}{\encodingdefault}{\sfdefault}{bx}{n}

\newcommand{\R}{\mathbb{R}}

\renewcommand{\div}{\nabla \cdot}

\providecommand{\grad}{\nabla}

\providecommand{\pder}[2]{\frac{\partial #1}{\partial #2}}

\renewcommand{\vec}[1]{{\bf #1}}